\pdfoutput=1
\documentclass[11pt]{article}
\usepackage[margin=1.1in]{geometry}
\usepackage[T1]{fontenc}
\usepackage[utf8]{inputenc}
\usepackage{lmodern}
\usepackage{amsmath,amssymb,amsthm,mathtools}
\usepackage{bm}
\usepackage{booktabs}
\usepackage{array}
\usepackage{graphicx}
\usepackage{enumitem}
\usepackage{xcolor}
\usepackage[round,authoryear]{natbib}
\usepackage[colorlinks=true,linkcolor=blue!60!black,citecolor=blue!60!black,urlcolor=blue!60!black]{hyperref}
\usepackage{microtype}

\theoremstyle{plain}
\newtheorem{theorem}{Theorem}[section]
\newtheorem{proposition}[theorem]{Proposition}
\newtheorem{corollary}[theorem]{Corollary}
\newtheorem{lemma}[theorem]{Lemma}
\theoremstyle{definition}
\newtheorem{assumption}[theorem]{Assumption}
\newtheorem{example}[theorem]{Example}

\theoremstyle{remark}
\newtheorem{remark}[theorem]{Remark}

\newcommand{\R}{\mathbb{R}}
\newcommand{\Rbar}{\overline{\mathbb{R}}}
\newcommand{\E}{\mathbb{E}}
\newcommand{\Prob}{\mathbb{P}}
\newcommand{\Pn}{\mathbb{P}_m}
\newcommand{\Gm}{\mathbb{G}_m}
\newcommand{\GQ}{\mathbb{G}_Q}
\newcommand{\Zm}{\mathbb{Z}_m}
\newcommand{\Lb}{\mathbb{L}}
\newcommand{\W}{W}
\newcommand{\ind}{\mathbf{1}}
\newcommand{\dd}{\,\mathrm{d}}
\newcommand{\Fhat}{\hat F_m}
\newcommand{\Fchk}{\check F_m}
\newcommand{\tauhat}{\hat\tau}
\newcommand{\tauchk}{\check\tau}
\newcommand{\Cm}{\mathcal{C}_m}
\newcommand{\linf}{\ell^\infty}
\newcommand{\weak}{\rightsquigarrow}
\newcommand{\toP}{\xrightarrow{\;p\;}}
\newcommand{\eqd}{\overset{d}{=}}
\newcommand{\Var}{\operatorname{Var}}
\newcommand{\Cov}{\operatorname{Cov}}

\newcommand{\Hcal}{\mathcal{H}}
\newcommand{\Vcal}{\mathcal{V}}
\newcommand{\Bcal}{\mathcal{B}}
\newcommand{\Scal}{\mathcal{S}}

\newcommand{\Fcal}{\mathcal{F}}

\defcitealias{pournaderi2026training}{PX26}
\defcitealias{pournaderi2026sharp}{P26}
\newcommand{\PX}{\textup{\citetalias{pournaderi2026training}}}
\newcommand{\Pt}{\textup{\citetalias{pournaderi2026sharp}}}

\DeclareRobustCommand{\codeurl}{\begingroup
  \expandafter\def\expandafter\UrlSpecials\expandafter{\UrlSpecials
    \do\-{\mathchar"002D\nobreak\hskip 0pt plus 1fil\penalty 0\hskip 0pt plus -1fil\relax}}%
  \url{https://github.com/mehrdad-pournaderi/Weighted-conformal-prediction-under-covariate-shift-process-limits-estimated-weights-and-efficiency}\endgroup}

\title{Weighted conformal prediction under covariate shift: process limits, estimated weights and efficiency}
\author{Mehrdad Pournaderi\thanks{Mofid Securities. E-mail: \texttt{m.pournaderi@emofid.com}.}}
\date{November 2026}

\begin{document}
\maketitle

\begin{abstract}
Weighted split conformal prediction handles covariate shift by reweighting calibration scores with the likelihood ratio $v=dQ_X/dP_X$; its training-conditional coverage is a functional of the self-normalized weighted empirical distribution function $\Fhat$, and we study the process $\sqrt m(\Fhat-F_Q)$. It converges in $\linf$ to the Gaussian process $\W(G(t))-F_Q(t)\W(K^2)$, with $G(t)=\E_Q[v\ind\{S\le t\}]$ and $K^2=1+\chi^2(Q\,\|\,P)$, if and only if $\chi^2(Q\,\|\,P)<\infty$ (for non-degenerate scores); $\sup v$ plays no role. The limit gives the joint law of the coverage across levels, consistent bootstraps, inflations asymptotically exact uniformly in the level, and Berry--Esseen-type rates. With estimated ratios, a uniform expansion remains valid when the estimate uses the calibration covariates; the estimation error enters through a covariance term that can reduce the variance below the oracle's, and a pairs bootstrap reproduces it. When the ratio is known, the weighted conformal threshold is semiparametrically efficient if and only if the ratio and the coverage indicator are uncorrelated under the test distribution, with an explicit efficient correction otherwise. When $\chi^2(Q\,\|\,P)=\infty$, the ratio has a regularly varying tail of index $p\in(1,2)$ and the score decouples from it in the far tail, the finite-dimensional limits are $p$-stable at the rate $m^{1-1/p}$ up to a slowly varying factor.
\end{abstract}


\section{Introduction}\label{sec:intro}

Split conformal prediction fits a predictor on a training split, computes nonconformity scores $S_i=s(X_i,Y_i)$ on a calibration split of size $m$, and outputs the prediction set $\{y:s(x,y)\le\hat\tau\}$ with $\hat\tau$ an empirical quantile of the calibration scores \citep{vovk2005algorithmic,angelopoulos2023gentle}. Under exchangeability the set has marginal coverage at least $1-\alpha$. Its \emph{training-conditional} coverage---the coverage of the one set a practitioner actually deploys, as a random variable of the calibration data---has been studied by \citet{vovk2012conditional}, \citet{bian2023training}, \citet{liang2025algorithmic} and \citet{duchi2025sample}, and is the object of PAC-type guarantees.

Under covariate shift \citep{shimodaira2000improving,sugiyama2007covariate}, where the test covariate law $Q_X$ differs from the training law $P_X$ while $P_{Y|X}$ is unchanged, \citet{tibshirani2019conformal} showed that marginal coverage is restored by weighting the calibration scores with the likelihood ratio $v=dQ_X/dP_X$, the test point's own weight being placed at $+\infty$. The object of this paper is the self-normalized importance-weighted empirical distribution function of the calibration scores,
\[
\Fhat(t)=\frac{\sum_{i=1}^m v(X_i)\ind\{S_i\le t\}}{\sum_{i=1}^m v(X_i)},
\]
and its quantile $\hat\tau(\beta)$, whose training-conditional miscoverage at level $\beta$ is $1-F_Q(\hat\tau(\beta))$, where $F_Q$ is the law of the score under the test distribution. This quantile alone does not carry the exact marginal guarantee (for $v\equiv1$, $m=10$ and $\beta=0.9$ it is the ninth order statistic, with marginal coverage $9/11$); the mass at $+\infty$ supplies it, and Proposition~\ref{prop:testweight} shows that, under the second-moment condition of the Gaussian regime, the two training-conditional miscoverage processes differ by $o_p(m^{-1/2})$ uniformly in the level, so the Gaussian-regime results of Sections~\ref{sec:coverage}--\ref{sec:rates} hold for both; the efficient threshold of Section~\ref{sec:efficiency} and the stable regime of Section~\ref{sec:stable} are treated on their own terms. Three papers have studied this random variable, and Table~\ref{tab:compare} summarizes what each covers. \citet{pournaderi2026training} (henceforth \PX{}) proved weighted Dvoretzky--Kiefer--Wolfowitz inequalities with an unspecified constant and derived training-conditional bounds for split, full and jackknife+ methods from them; at confidence $1-\delta$ the deviation of the coverage is of order $B\sqrt{\log(1/\delta)/m}$ when $v\le B$, and of order $B/(\delta\sqrt m)$ when only $\|v\|_{P,2}\le B$ is assumed. \citet{pournaderi2026sharp} (henceforth \Pt{}) observed that at a fixed level no uniformity is needed, which gives explicit finite-sample bounds at the scale of a variance proxy built from $K^2=1+\chi^2(Q\,\|\,P)$ rather than from $B$, a fixed-level central limit theorem (its Proposition~10), lower bounds, and an analysis of clipping and of estimated ratios independent of the calibration sample. \citet{gazin2024asymptotics} obtained, by the functional delta method, the process limit of the false coverage proportion of weighted split conformal sets over a growing test sample, for bounded deterministic weights that are functions of the score; when the test sample dominates, this is the training-conditional coverage process for such weights (Remark~\ref{rem:gazin}), and that paper names the exact asymptotic distribution under estimated weights as an open problem.

This paper studies the \emph{process}
\[
\Zm(t)=\sqrt m\,\bigl(\Fhat(t)-F_Q(t)\bigr),\qquad t\in\Rbar,
\]
which we call the weighted score process, in the training-conditional regime, for covariate weights, and without boundedness.

\paragraph{Contributions.}
The weak convergence of $\Zm$ itself is an application of Donsker theory, and we present it as the foundation rather than as the advance. The advance is what the process view makes possible for two questions that neither the finite-sample bounds nor the pointwise limit can address, \emph{estimated weights} and \emph{efficiency}, which are the center of the paper. The results fall into three groups.
\begin{enumerate}[leftmargin=*,itemsep=3pt]
\item \emph{The Gaussian regime and uniform inference} (Sections~\ref{sec:fclt}--\ref{sec:rates}). If $K^2=\E_Pv^2<\infty$ then $\Zm\weak\GQ$ in $\linf(\Rbar)$, where
\[
\GQ(t)\eqd \W\bigl(G(t)\bigr)-F_Q(t)\,\W(K^2),\qquad G(t)=\E_Q\bigl[v\ind\{S\le t\}\bigr],
\]
for a standard Brownian motion $\W$: a Brownian motion run on the ``second-moment clock'' $G$ and pinned by $F_Q$, which reduces to the Brownian bridge when $v\equiv1$. The condition is necessary, by a converse central limit theorem (Proposition~\ref{prop:necessity}), so $\chi^2(Q\,\|\,P)<\infty$ is exactly the condition for the $\sqrt m$ rate and $\sup v$ is irrelevant. The training-conditional coverage error converges to $\GQ\circ q_Q$ as a process in the level (Theorem~\ref{thm:coverage}), which gives an asymptotic weighted DKW inequality at the scale $K$, consistent multiplier and nonparametric bootstraps, asymptotically exact inflations that certify coverage \emph{simultaneously} over a set or an interval of levels, and a two-sample diagnostic for the covariate-shift assumption; these are consequences of the limit by standard arguments, stated because they are what one uses. For bounded ratios, the duality between the weighted quantile and the weighted distribution function gives a Berry--Esseen bound at the parametric rate for a single level, and a weighted Bahadur--Kiefer bound gives the maximal coverage error over a grid of levels its Gaussian and bootstrap approximations at the rate $m^{-1/4}$ up to logarithms, so that the simultaneous guarantee has a finite-sample form.
\item \emph{Estimated ratios and efficiency} (Sections~\ref{sec:estimated}--\ref{sec:efficiency}). The ratio is fitted from unlabeled covariate samples that almost always include the calibration covariates, so the weights and the scores they weight are dependent. If the estimate $\hat v$ ranges over a Donsker-type class and is $L_2(P)$-consistent, then $\sqrt m(\Fchk-F_{\hat Q})=\Zm+o_p(1)$ uniformly, where $\Fchk$ uses $\hat v$ and $\hat Q$ is the $\hat v$-tilted law, \emph{whether or not} $\hat v$ depends on the calibration covariates (Theorem~\ref{thm:estimated}): the natural centering makes the dependence harmless at first order. For parametric ratio models the estimation error enters through the coupling gradient $\gamma(t)=\E_P[\dot v(X)\{F_{S|X}(t)-F_Q(t)\}]$ (Theorem~\ref{thm:parametric}), and when the calibration covariates are part of the sample used to fit the ratio the asymptotic variance can be \emph{smaller} than with the true ratio (Example~\ref{ex:meanshift}), the conformal analogue of \citet{hirano2003efficient}. A pairs bootstrap that recomputes the weights reproduces this law; the usual bootstrap with the weights held fixed does not (Theorem~\ref{thm:estboot}). The efficiency analysis explains the phenomenon: with a known ratio, so that $\E_Pv=1$ is a known constraint, the self-normalized weighted distribution function is efficient if and only if $G(t)=K^2F_Q(t)$, that is, if and only if the ratio and the coverage indicator are uncorrelated under $Q$; otherwise a control-variate correction computed from the calibration scores and weights alone attains the bound, at the price of the finite-sample marginal guarantee, and the loss of the standard threshold is the squared covariance between the ratio and the coverage indicator divided by the variance of the ratio (Theorems~\ref{thm:eff}--\ref{thm:effthreshold}). Both results rest on the same mechanism, a centering that makes a self-consistency identity exact (Lemma~\ref{lem:identity}), and both are invisible at the resolution of the finite-sample bounds.
\item \emph{Beyond finite $\chi^2$} (Section~\ref{sec:stable}). If $\Prob(v>t)$ is regularly varying with index $-p$, $p\in(1,2)$, and the score decouples from the ratio in the far tail of the shift, the weighted score process converges in finite-dimensional distributions at the rate $m^{1-1/p}$, up to a slowly varying factor, to a $p$-stable process, and the coverage error at a fixed level converges to an explicit $p$-stable law with skewness
\[
\frac{\pi_\beta(1-\beta)^p-(1-\pi_\beta)\beta^p}{\pi_\beta(1-\beta)^p+(1-\pi_\beta)\beta^p},
\]
where $\pi_\beta$ is the probability that a score from the far tail of the shift falls inside the prediction set: the skewness records where the far tail of the shift lands relative to the set. Gaussian inflations are then invalid and the required inflation scales as $m^{-(1-1/p)}$.
\end{enumerate}
Section~\ref{sec:experiments} checks each statement numerically and Section~\ref{sec:discussion} lists open problems. Proofs are in the appendix.

\paragraph{Relation to the literature.}
Empirical processes weighted by a fixed square-integrable function are classical \citep{vandervaart1996weak}; weighted empirical measures in importance sampling have been studied from the large-deviation side by \citet{hult2016large}, under exponential-moment conditions far stronger than the second moment used here. The inverse-probability-weighted empirical process with estimated weights is the object of \citet{breslow2007weighted}, in two-phase sampling, and the variance-reduction phenomenon for estimated weights goes back to \citet{hirano2003efficient}, with its two-sample, nonparametric form in \citet{chen2008semiparametric}. On the conformal side, process limits for the exchangeable false coverage proportion over a test sample were obtained by \citet{gazin2024asymptotics}, building on the finite-sample results of \citet{marques2025universal} and \citet{gazin2024transductive}. Among the finite-sample results under covariate shift, \citet{park2022pac} construct PAC sets by rejection sampling at the scale $B$; \citet{qiu2023prediction} and \citet{yang2024doubly} obtain asymptotic validity under unknown shift through doubly robust constructions; \citet{lei2021conformal} and \citet{barber2023conformal} bound the marginal coverage loss from estimated or fixed weights by total-variation quantities; \citet{wang2026weight} clip unbounded ratios. The regime $\chi^2(Q\,\|\,P)=\infty$ is the regime in which importance sampling is known to be fragile \citep{chatterjee2018sample,vehtari2024pareto}, and Section~\ref{sec:stable} makes this precise for conformal calibration.

\begin{table}[t]
\centering\small
\caption{What each paper on the training-conditional coverage of weighted conformal prediction provides ($\checkmark$: yes; --: no). G24 is \citet{gazin2024asymptotics}. Throughout, $v$ is the likelihood ratio, $m$ the calibration size and $1-\delta$ the confidence level.}
\label{tab:compare}
\begin{tabular}{@{}l@{\hskip 16pt}c@{\hskip 14pt}c@{\hskip 14pt}c@{\hskip 14pt}c@{}}
\toprule
 & \PX{} & \Pt{} & G24 & This paper\\
\midrule
\multicolumn{5}{@{}l}{\emph{Kind of result}}\\
\quad Guarantee at finite $m$ & $\checkmark$ & $\checkmark$ & -- & --\\
\quad Exact limit law of the coverage & -- & one level & $\checkmark$ & $\checkmark$\\
\quad All levels simultaneously & $\checkmark$ & -- & $\checkmark$ & $\checkmark$\\
\quad Computable critical values & --$^{\,\rm a}$ & explicit & plug-in & bootstrap\\
\quad Rates of the approximations & -- & -- & -- & $\checkmark$\\
\addlinespace[5pt]
\multicolumn{5}{@{}l}{\emph{Likelihood ratio}}\\
\quad Weights that depend on the covariates & $\checkmark$ & $\checkmark$ & --$^{\,\rm b}$ & $\checkmark$\\
\quad Unbounded ratio with $\E_Pv^2<\infty$ & $\checkmark^{\,\rm c}$ & truncation & -- & $\checkmark^{\,\rm d}$\\
\quad Ratio with $\E_Pv^2=\infty$ & consistency & truncation & -- & stable limits\\
\quad Ratio estimated on separate data & --$^{\,\rm e}$ & $\checkmark$ & -- & $\checkmark$\\
\quad Ratio estimated on the calibration covariates & --$^{\,\rm e}$ & -- & --$^{\,\rm e}$ & $\checkmark$\\
\quad Bootstrap with an estimated ratio & -- & -- & -- & $\checkmark$\\
\addlinespace[5pt]
\multicolumn{5}{@{}l}{\emph{Scope}}\\
\quad Efficiency bound and efficient threshold & -- & -- & -- & $\checkmark$\\
\quad Lower bound or necessity result & -- & $\checkmark$ & -- & $\checkmark$\\
\quad Full conformal and jackknife+ & $\checkmark$ & -- & -- & --\\
\bottomrule
\end{tabular}

\vspace{5pt}
\begin{minipage}{0.92\linewidth}\footnotesize\raggedright
$^{\rm a}$\,The bounds contain an unspecified universal constant.\\[1pt]
$^{\rm b}$\,Bounded weights that are deterministic functions of the score.\\[1pt]
$^{\rm c}$\,At confidence $1-\delta$ the deviation of the coverage is of order $B/(\delta\sqrt m)$ under $\|v\|_{P,2}\le B$, against $B\sqrt{\log(1/\delta)/m}$ under $v\le B$.\\[1pt]
$^{\rm d}$\,The condition is also necessary for the $\sqrt m$ rate; at a level $\beta$ the deviation is $z_{1-\delta}\,\sigma(\beta)/\sqrt m$, with $\sigma(\beta)\le\|v\|_{P,2}$.\\[1pt]
$^{\rm e}$\,\PX{} gives the $\ell_1$ error term to be added to its bounds for any estimate, without rates; G24 names the limit law under estimated weights as an open problem. The results of this paper are for parametric ratio models.
\end{minipage}
\end{table}%

\section{Setting and a basic identity}\label{sec:setting}

Throughout we condition on the training split, so that the score function $s:\mathcal X\times\mathcal Y\to\R$ is fixed. The calibration sample $Z_i=(X_i,Y_i)$, $i=1,\dots,m$, is i.i.d.\ from $P=P_X\times P_{Y|X}$, and the test point $Z=(X,Y)$ is drawn from $Q=Q_X\times P_{Y|X}$, independently of the calibration sample. We assume $Q\ll P$ and write
\[
v(z)=\frac{dQ}{dP}(z)=\frac{dQ_X}{dP_X}(x),\qquad v_i=v(X_i),\qquad S_i=s(Z_i),\qquad S=s(Z),
\]
so that
\[
\E_P\bigl[v(X)\,g(Z)\bigr]=\E_Q\bigl[g(Z)\bigr]\quad\text{for every integrable }g,\qquad\text{in particular}\qquad \E_Pv=1 .
\]
Let
\begin{gather*}
F_Q(t)=Q(S\le t),\qquad F_{S|X}(t)=\Prob(S\le t\mid X),\\
q_Q(\beta)=\inf\{t:F_Q(t)\ge\beta\},\qquad P_e(\tau)=1-F_Q(\tau),
\end{gather*}
where $F_{S|X}$ is the same under $P$ and under $Q$ because $P_{Y|X}$ is shared, and
\[
F_Q(t)=\E_Q\bigl[F_{S|X}(t)\bigr]=\E_P\bigl[v(X)\,F_{S|X}(t)\bigr].
\]
The self-normalized weighted empirical distribution function of the calibration scores and its $\beta$-quantile are
\begin{equation}\label{eq:Fhat}
\Fhat(t)=\sum_{i=1}^m\hat w_i\ind\{S_i\le t\},\qquad \hat w_i=\frac{v_i}{\sum_{j=1}^m v_j},\qquad
\tauhat(\beta)=\inf\{t:\Fhat(t)\ge\beta\}\in\R\cup\{+\infty\},
\end{equation}
with $\Fhat\equiv0$ on the event $\sum_jv_j=0$, whose probability $\Prob_P(v=0)^m$ tends to zero since $\E_Pv=1$ forces $\Prob_P(v=0)<1$ (it is at most $(1-1/B)^m\le e^{-m/B}$ when $v\le B$). On this event $\tauhat\equiv+\infty$, the prediction set is the whole space, and every statistic built from the weights below (multiplier and bootstrap processes, the variance estimator \eqref{eq:sigmahat}, studentized statistics) is set to zero; the same conventions apply to estimated weights and to bootstrap samples. The event is irrelevant for the limit statements, and the finite-sample statements of Section~\ref{sec:rates} account for its probability separately. The prediction set $\{y:s(x,y)\le\tauhat(\beta)\}$ has training-conditional miscoverage exactly $P_e(\tauhat(\beta))$, and we write
\begin{equation}\label{eq:Cm}
\Cm(\beta)=P_e\bigl(\tauhat(\beta)\bigr)-(1-\beta)=\beta-F_Q\bigl(\tauhat(\beta)\bigr),\qquad\beta\in(0,1),
\end{equation}
for the coverage error at level $\beta$ (positive values are under-coverage). The weighted split conformal set of \citet{tibshirani2019conformal} also gives the test point's weight $v(X)$, placed at $+\infty$, a share of the quantile; its miscoverage $P_e^{\rm split}(D;\beta)$ is treated in Proposition~\ref{prop:testweight}.

The two population quantities that govern everything are
\begin{equation}\label{eq:KG}
K^2:=\E_Pv^2=\E_Qv=1+\chi^2(Q\,\|\,P)\in[1,\infty],\qquad G(t):=\E_Q\bigl[v\ind\{S\le t\}\bigr]=\E_P\bigl[v^2\ind\{S\le t\}\bigr],
\end{equation}
so that $G$ is nondecreasing with $G(-\infty)=0$ and $G(+\infty)=K^2$. When $v\le B$ one has
\[
K^2\le B\qquad\text{and}\qquad G(t)-G(s)\le B\,\{F_Q(t)-F_Q(s)\}\quad\text{for }s<t .
\]
Finally, for $t\in\Rbar$ put
\begin{equation}\label{eq:h}
h_t(z):=v(z)\bigl(\ind\{s(z)\le t\}-F_Q(t)\bigr),\qquad\text{so that}\qquad Ph_t=F_Q(t)-F_Q(t)=0,
\end{equation}
with $h_{\pm\infty}\equiv0$, and write $\Pn$ for the empirical measure of $Z_1,\dots,Z_m$ and
\[
\Gm=\sqrt m\,(\Pn-P)
\]
for the empirical process. The whole paper rests on the following observation.

\begin{lemma}[basic identity]\label{lem:identity}
On the event $\{\sum_jv_j>0\}$, for every $t\in\Rbar$,
\[
\Fhat(t)-F_Q(t)=\frac{\Pn h_t}{\Pn v},\qquad\text{hence}\qquad
\Zm(t):=\sqrt m\bigl(\Fhat(t)-F_Q(t)\bigr)=\frac{\Gm h_t}{\Pn v}.
\]
\end{lemma}

\begin{proof}
Writing out the definitions,
\[
\Fhat(t)-F_Q(t)=\frac{\Pn[v\ind\{S\le t\}]-F_Q(t)\,\Pn v}{\Pn v}=\frac{\Pn h_t}{\Pn v},
\qquad\text{and}\qquad
\Pn h_t=(\Pn-P)h_t
\]
because $Ph_t=0$.
\end{proof}

The self-normalization is thus a control variate: the numerator is the centered process $\Gm h_t$ of a single function class, with no separate denominator term to control. In \PX{} the two were separated by a triangle inequality and bounded one at a time; in \Pt{} the identity was used at a single fixed $t$. The remainder of the paper is the uniform-in-$t$ theory of $\Gm h_t$ and what it implies. Suprema over $t$ coincide with suprema over the rationals, by right-continuity of $t\mapsto h_t(z)$, so no measurability issue arises; we use outer probabilities where \citet{vandervaart1996weak} do.%

\section{Weak convergence of the weighted score process}\label{sec:fclt}

\subsection{The functional central limit theorem}

Let $\Hcal=\{h_t:t\in\Rbar\}$ with $h_t$ as in \eqref{eq:h}. The following covariance function will be the covariance of the limit:
\begin{equation}\label{eq:C}
\begin{aligned}
C(s,t):=P(h_sh_t)&=\E_P\bigl[v^2(\ind\{S\le s\}-F_Q(s))(\ind\{S\le t\}-F_Q(t))\bigr]\\
&=G(s\wedge t)-F_Q(s)G(t)-F_Q(t)G(s)+K^2F_Q(s)F_Q(t).
\end{aligned}
\end{equation}

\begin{theorem}[functional CLT]\label{thm:fclt}
Assume $K^2=\E_Pv^2<\infty$. Then $\Hcal$ is $P$-Donsker and
\[
\Zm\weak\GQ\qquad\text{in }\linf(\Rbar),
\]
where $\GQ$ is a tight centered Gaussian process with covariance \eqref{eq:C}. Equivalently,
\begin{equation}\label{eq:brownian}
\{\GQ(t)\}_{t\in\Rbar}\eqd\bigl\{\W(G(t))-F_Q(t)\,\W(K^2)\bigr\}_{t\in\Rbar}
\end{equation}
for a standard Brownian motion $\W$ on $[0,K^2]$. In particular $\GQ(\pm\infty)=0$,
\[
\Var\GQ(t)=\E_Q\bigl[v\,(\ind\{S\le t\}-F_Q(t))^2\bigr]\le\min\bigl\{K^2,\,B\,F_Q(t)(1-F_Q(t))\bigr\}
\]
(the second bound when $v\le B$), and if $F_Q$ is continuous then $G$ is continuous and $\GQ$ has continuous sample paths. When $v\equiv1$, $G=F_Q$, $K^2=1$ and $\GQ=\mathbb B\circ F_Q$ is the Brownian bridge.
\end{theorem}

The proof (Appendix~\ref{app:fclt}) computes the $L_2(P)$-bracketing numbers of $\Hcal$,
\[
N_{[\,]}\bigl(\epsilon,\Hcal,L_2(P)\bigr)\le\frac{8K^2}{\epsilon^2}+\frac{4K}{\epsilon}+3,
\]
which makes the bracketing entropy integral finite, and invokes Ossiander's bracketing central limit theorem \citep{ossiander1987central,vandervaart1996weak}; the denominator $\Pn v\to1$ a.s.\ is removed by Slutsky's lemma. The bracketing computation is the one used in \PX{} to bound the \emph{expectation} of the supremum. The representation \eqref{eq:brownian} follows from $G(s)\wedge G(t)=G(s\wedge t)$ for the nondecreasing $G$.

\begin{remark}[reading the limit]\label{rem:reading}
The clock $G$ accumulates $\E_Q[v\ind\{S\in dt\}]$: it runs fast where the test distribution puts mass \emph{and} the likelihood ratio is large. The process is pinned at both ends because $\Fhat(\pm\infty)$ is exactly $0$ and $1$; this is the effect of self-normalization, and it is why the variance
\[
\Var\GQ(t)=\bigl(1-F_Q(t)\bigr)^2\,G(t)+F_Q(t)^2\,\bigl(K^2-G(t)\bigr)
\]
is a linear function of $G(t)\in[0,K^2]$, at most $K^2\max\{F_Q(t),1-F_Q(t)\}^2$; at a quantile it is the fixed-level variance of \Pt{} (Section~\ref{sec:coverage}).
\end{remark}

\subsection{Sharpness: the second moment is necessary}

\begin{proposition}[necessity of $K^2<\infty$]\label{prop:necessity}
Let $t\in\R$ with $0<F_Q(t)<1$. If $K^2=\infty$, then the sequence $\Zm(t)$ is not tight; in particular it does not converge in distribution to any law.
\end{proposition}

The proof (Appendix~\ref{app:necessity}) shows that convergence in distribution of $\sqrt m\,\Pn h_t$ along any subsequence forces $Ph_t^2<\infty$, through the characteristic function of a symmetrized summand, and that
\[
Ph_t^2=\E_Q\bigl[v\,(\ind\{S\le t\}-F_Q(t))^2\bigr]\ge\min\{F_Q(t),1-F_Q(t)\}^2\,K^2 .
\]
Together with Theorem~\ref{thm:fclt}:

\begin{corollary}\label{cor:iff}
Suppose $S$ is nondegenerate under $Q$. The following are equivalent: (i) $\chi^2(Q\,\|\,P)<\infty$; (ii) $\Zm$ converges weakly in $\linf(\Rbar)$ to a tight limit; (iii) $\Zm(t)$ is tight for some $t$ with $0<F_Q(t)<1$. When they hold, the limit is \eqref{eq:brownian}.
\end{corollary}

Thus the $\sqrt m$ rate for the weighted score process is governed exactly by the $\chi^2$-divergence of the shift, and by nothing about the supremum of the ratio; Section~\ref{sec:coverage} transfers the statement to the training-conditional coverage under continuity of $F_Q$. This is the asymptotic counterpart of the message of \Pt{}, where $K^2$ appears as the intrinsic scale through a two-point lower bound; here it is the threshold between the Gaussian regime and the regime of Section~\ref{sec:stable}.

\subsection{A finite-sample counterpart}

The same control-variate identity sharpens the weighted DKW inequality of \PX{} in finite samples. Let
\[
\sigma_*^2:=\sup_t Ph_t^2=\sup_t\Var\GQ(t),\qquad\text{so that}\qquad\sigma_*^2\le\min\{K^2,B/4\}\ \text{ when }v\le B .
\]

\begin{proposition}[weighted DKW at the variance scale]\label{prop:bousquet}
Assume $Q\ll P$ and $v\le B$. There is a universal constant $C_1$ such that for all $u>0$ and $\eta\in(0,1)$, with probability at least $1-e^{-u}-e^{-2m\eta^2/B^2}$,
\[
\sup_{t\in\R}\bigl|\Fhat(t)-F_Q(t)\bigr|\le\frac{1}{1-\eta}\left[\frac{C_1K}{\sqrt m}+\sqrt{\frac{2\bigl(\sigma_*^2+2C_1BK/\sqrt m\bigr)u}{m}}+\frac{Bu}{3m}\right].
\]
\end{proposition}

The proof (Appendix~\ref{app:bousquet}) applies Bousquet's inequality \citep{bousquet2002bennett,boucheron2013concentration} to $\sup_t|\Gm h_t|$, whose summands are bounded by $B$ and have variance at most $\sigma_*^2$, and controls the denominator by Hoeffding's inequality. Compared with Lemma~1 of \PX{}, the deviation term changes from
\[
B\sqrt{\frac{2\log(4/\delta)}{m}}\qquad\text{to}\qquad\sqrt{\frac{2\sigma_*^2\log(1/\delta)}{m}},\qquad\sigma_*^2\le\min\{K^2,B/4\},
\]
at least a factor $2\sqrt B$ smaller, while the supremum $B$ survives only in the lower-order terms. The unknown constant $C_1$ of the maximal inequality remains, as in \PX{}; \Pt{} removes constants altogether, but only at a fixed level. Theorem~\ref{thm:fclt} says that asymptotically the leading term is exactly of the order $\sigma_*$ (Corollary~\ref{cor:asympdkw} below gives the Gaussian tail of the limit), so Proposition~\ref{prop:bousquet} has the right scale.

\section{The coverage process and uniform training-conditional inference}\label{sec:coverage}

\subsection{The coverage error as a process in the level}

Recall $\Cm(\beta)=\beta-F_Q(\tauhat(\beta))$ from \eqref{eq:Cm}, and fix $0<\beta_0<\beta_1<1$. We write
\begin{equation}\label{eq:Gtilde}
\tilde G(\beta):=G\bigl(q_Q(\beta)\bigr)=\E_Q\bigl[v\ind\{F_Q(S)\le\beta\}\bigr]
\end{equation}
for the clock $G$ expressed in the level rather than in the threshold.

\begin{theorem}[coverage process]\label{thm:coverage}
Assume $K^2<\infty$ and that $F_Q$ is continuous. Then
\[
\sup_{\beta\in[\beta_0,\beta_1]}\bigl|\sqrt m\,\Cm(\beta)-\Zm(q_Q(\beta))\bigr|\toP0,
\qquad\text{hence}\qquad
\sqrt m\,\Cm\weak\Lb:=\GQ\circ q_Q\quad\text{in }\linf[\beta_0,\beta_1],
\]
where
\[
\{\Lb(\beta)\}_{\beta}\eqd\bigl\{\W(\tilde G(\beta))-\beta\,\W(K^2)\bigr\}_{\beta},\qquad
\sigma^2(\beta):=\Var\Lb(\beta)=(1-\beta)^2\tilde G(\beta)+\beta^2\bigl(K^2-\tilde G(\beta)\bigr).
\]
\end{theorem}

No density of the scores is required and no strict monotonicity of $F_Q$: the statement is about $F_Q(\tauhat(\beta))$, which under continuity of $F_Q$ is the self-normalized weighted $\beta$-quantile of the $Q$-uniform variables $F_Q(S_i)$. The proof (Appendix~\ref{app:coverage}) is a weighted Bahadur--Kiefer argument: it establishes the exact representation
\[
\sqrt m\,\Cm(\beta)=\Zm\bigl(\tauhat(\beta)\bigr)-\sqrt m\,r_m(\beta),\qquad 0\le r_m(\beta)\le\max_i\hat w_i=o_p(m^{-1/2}),
\]
and the asymptotic equicontinuity of $\Zm$ transfers $\tauhat(\beta)$ to $q_Q(\beta)$ uniformly in the level. At a single level the theorem gives
\[
\sqrt m\,\Cm(\beta)\weak N\bigl(0,\sigma^2(\beta)\bigr),\qquad\sigma^2(1-\alpha)=\alpha^2K^2+(1-2\alpha)\,\E_Q\bigl[v\ind\{S>q_Q(1-\alpha)\}\bigr],
\]
which is Proposition~10 of \Pt{}, obtained there by a Lindeberg--Feller argument at a moving quantile under $v\le B$; Theorem~\ref{thm:coverage} removes the boundedness, and its content is the joint law. For $v\equiv1$, $\Lb$ is the Brownian bridge in $\beta$ and $\sigma^2(\beta)=\beta(1-\beta)$, matching the Beta law of exchangeable split conformal coverage \citep{vovk2012conditional}. If $v(X)$ and $\ind\{S\le q_Q(\beta)\}$ are independent under $Q$, then $\tilde G(\beta)=\beta K^2$ and
\[
\sigma^2(\beta)=K^2\,\beta(1-\beta):
\]
$m$ weighted calibration points behave like $m/K^2$ unweighted ones, which is Kish's effective sample size $(\sum v_i)^2/\sum v_i^2\approx m/K^2$ \citep{kish1965survey}. In general the effective sample size is not the whole story: in Table~\ref{tab:coupling} of Section~\ref{sec:experiments}, three shifts with the same $K^2=16/7$---weight on large scores, on small scores, or on a covariate independent of the score---have $\sigma^2(0.9)=0.246$, $0.114$ and $0.206$, a factor two apart with the same $m/K^2$, in the directions given by the sign of $\beta K^2-\tilde G(\beta)$.

\begin{remark}[relation to the false coverage proportion limit of Gazin]\label{rem:gazin}
Decomposing the Brownian motion on $[0,K^2]$ into a bridge and its endpoint, $\W(u)=K\,\mathbb B(u/K^2)+(u/K^2)\W(K^2)$ with $\mathbb B$ a standard Brownian bridge independent of $N:=\W(K^2)/K$, the limit of Theorem~\ref{thm:coverage} becomes
\[
\Lb(\beta)\eqd K\Bigl\{\mathbb B\bigl(\tilde G(\beta)/K^2\bigr)+\bigl(\tilde G(\beta)/K^2-\beta\bigr)N\Bigr\}.
\]
This is the structure of the calibration-side term of Theorem~3.4 of \citet{gazin2024asymptotics}, who obtains, by the functional delta method and under bounded weights and continuously differentiable score distributions, the limit of the false coverage proportion of weighted split conformal sets over a test sample of size $n_{\mathrm{test}}$ growing with $m$, for weights that are deterministic functions \emph{of the score}: with $\rho^w=\|w\|_{P,2}/\|w\|_{P,1}$ and $I^w$ the normalized $w^2$-clock, under the normalization $\sqrt{mn_{\mathrm{test}}/(m+n_{\mathrm{test}})}$ his limit is
\[
\Bigl(1+\frac{m}{n_{\mathrm{test}}}\Bigr)^{-1/2}\rho^w\bigl(\mathbb B\circ I^w+(I-I^w)N\bigr)\;+\;\Bigl(1+\frac{n_{\mathrm{test}}}{m}\Bigr)^{-1/2}\,\mathbb U,
\]
with the two sample-size factors taken at their limits and $\mathbb U$ a second Brownian bridge, independent of $(\mathbb B,N)$: a calibration-side term and a test-side term. When the test sample dominates, $n_{\mathrm{test}}/m\to\infty$, the normalization is asymptotically $\sqrt m$, the first factor tends to one and the test-side term vanishes. In that regime the two statements coincide when $v(X)$ is a function of $S$: then $\rho^w=K$ and, with his miscoverage level $\alpha=1-\beta$, $I^w(\alpha)=1-\tilde G(\beta)/K^2$. In general they differ: weighting by $v(X)$ and by the score-based ratio $\E_P[v\mid S]=dQ_S/dP_S$ are different procedures, and $(\rho^w)^2=\E_P[\E_P(v\mid S)^2]\le K^2$ by Jensen's inequality, with equality if and only if $v$ is a function of the score; when $X$ and $S$ are independent, for instance, $\rho^w=1<K$ while the covariate-weighted procedure still pays the full price $K^2$. His Theorem~3.3 also covers a fixed weight function that is not the ratio, centred at a non-nominal curve; weights estimated from the calibration sample (Section~\ref{sec:estimated}) are outside its scope. Theorem~\ref{thm:coverage} is the statement for covariate weights under $K^2<\infty$ and continuity of $F_Q$ alone, with boundedness replaced by the bracketing bound of Lemma~\ref{lem:bracketH} and the test-point weight handled by Proposition~\ref{prop:testweight}. The bridge form shows that the law is that of an unweighted bridge at the clock $\tilde G/K^2$, tilted by the random level shift $(\tilde G/K^2-\beta)N$; the Brownian form of Theorem~\ref{thm:fclt} is the one that composes with $q_Q$ and that the bootstrap of Section~\ref{sec:bootstrap} reproduces.
\end{remark}

The weighted split conformal set of \citet{tibshirani2019conformal} attaches the test point's weight to $+\infty$, which shifts the quantile level by the random amount $\beta v(X)/\sum_jv_j$; its miscoverage is
\[
P_e^{\rm split}(D;\beta)=Q\Bigl(S>\tauhat\bigl(\beta\,(1+v(X)/\textstyle\sum_jv_j)\bigr)\Bigr),\qquad\tauhat(u):=+\infty\ \text{ for }u>1,
\]
which is at most $P_e(\tauhat(\beta))$ \citep[eq.\ (3)]{pournaderi2026sharp}. The difference is negligible at the $\sqrt m$ scale:

\begin{proposition}[the test weight is negligible]\label{prop:testweight}
Under the assumptions of Theorem~\ref{thm:coverage},
\[
\sqrt m\sup_{\beta\in[\beta_0,\beta_1]}\bigl|P_e^{\rm split}(D;\beta)-P_e(\tauhat(\beta))\bigr|\toP0 .
\]
Consequently Theorem~\ref{thm:coverage} holds verbatim for $P_e^{\rm split}(D;\beta)-(1-\beta)$.
\end{proposition}

The proof (Appendix~\ref{app:coverage}) truncates $v(X)$ at a level $a_m$ with
\[
\frac{a_m}{\sqrt m}\to0\qquad\text{and}\qquad\frac{\sqrt m}{a_m}\,\E_Q\bigl[v\ind\{v>a_m\}\bigr]\to0,
\]
which is possible because $\E_Qv=K^2<\infty$, and uses the equicontinuity of Theorem~\ref{thm:coverage} for the remainder.

\subsection{An asymptotic weighted DKW inequality at the scale $K$}

\begin{corollary}\label{cor:asympdkw}
Assume $K^2<\infty$. Then $\sqrt m\,\|\Fhat-F_Q\|_\infty\weak\|\GQ\|_\infty$, and there is a universal constant $c_0$ such that $\E\|\GQ\|_\infty\le c_0K$ and, for all $u>0$,
\[
\Prob\bigl(\|\GQ\|_\infty\ge\E\|\GQ\|_\infty+u\bigr)\le e^{-u^2/(2\sigma_*^2)},\qquad\sigma_*^2=\sup_t\Var\GQ(t)\le K^2 .
\]
Consequently, for every $\delta\in(0,1)$,
\[
\limsup_{m\to\infty}\Prob\Bigl(\sqrt m\,\|\Fhat-F_Q\|_\infty>c_0K+K\sqrt{2\log(1/\delta)}\Bigr)\le\delta,
\]
and the same bound holds for $\sqrt m\sup_{\beta\in[\beta_0,\beta_1]}|\Cm(\beta)|$ when $F_Q$ is continuous.
\end{corollary}

The first claim is the continuous mapping theorem, the tail bound is the Borell--TIS inequality \citep{adler2007random}, and the expectation bound is Dudley's entropy bound with the covering numbers $O(K^2/\epsilon^2)$ of $\Rbar$ in the canonical metric of $\GQ$ (Appendix~\ref{app:coverage}). The statement is Lemma~1 of \PX{} with $B$ replaced by $K\le\sqrt B$ in both terms, at the price of being asymptotic (Lemma~2 of \PX{} has the scale $K$ at finite $m$, with a factor $1/\delta$ in place of $\sqrt{2\log(1/\delta)}$); for $v\equiv1$ it is the Kolmogorov--Smirnov asymptotics. Unlike the exchangeable case, the limit law of $\|\GQ\|_\infty$ depends on $(F_Q,G)$, which is why the bootstrap of the next subsection is needed for exact asymptotic constants.

\subsection{Bootstrap}\label{sec:bootstrap}

Let $\xi_1,\dots,\xi_m$ be i.i.d.\ multipliers, independent of the data, with
\[
\E\xi=0,\qquad\Var\xi=1,\qquad\|\xi\|_{2,1}=\int_0^\infty\sqrt{\Prob(|\xi|>x)}\,\dd x<\infty,
\]
and define the multiplier process and its composition with the estimated quantile,
\begin{equation}\label{eq:multiplier}
\hat{\mathbb Z}^\xi_m(t)=\frac{1}{\sqrt m\,\Pn v}\sum_{i=1}^m\xi_i\,v_i\bigl(\ind\{S_i\le t\}-\Fhat(t)\bigr),\qquad
\hat\Lb^\xi_m(\beta)=\hat{\mathbb Z}^\xi_m\bigl(\tauhat(\beta)\bigr).
\end{equation}
Centering at $\Fhat(t)$ rather than at the unknown $F_Q(t)$ is what makes \eqref{eq:multiplier} computable. Write $\weak_\xi$ for conditional weak convergence given the data, in outer probability, in the sense of \citet[Section~2.9]{vandervaart1996weak}.

\begin{theorem}[multiplier bootstrap]\label{thm:bootstrap}
Assume $K^2<\infty$. Then $\hat{\mathbb Z}^\xi_m\weak_\xi\GQ$ in $\linf(\Rbar)$. If moreover $F_Q$ is continuous, then $\hat\Lb^\xi_m\weak_\xi\Lb$ in $\linf[\beta_0,\beta_1]$. The same conclusions hold for Efron's nonparametric bootstrap, i.e.\ for $\sqrt m(\Fhat^*-\Fhat)$ computed from a sample drawn with replacement from $\{(S_i,v_i)\}$, and for its composition with $\tauhat$.
\end{theorem}

The proof (Appendix~\ref{app:bootstrap}) is the conditional multiplier central limit theorem for the Donsker class $\Hcal$ with square-integrable envelope, plus the observation that the centering error $-(\Fhat(t)-F_Q(t))\,m^{-1/2}\sum_i\xi_iv_i=o_p(1)\cdot O_p(1)$ is uniformly negligible, and the equicontinuity argument of Theorem~\ref{thm:coverage} for the composition. The identity of Lemma~\ref{lem:identity} holds for the resampled sample with $\Fhat$ in the role of $F_Q$, which is what makes Efron's bootstrap work without a delta method.

\subsection{Simultaneous training-conditional guarantees}

For a single level, Theorem~\ref{thm:coverage} gives the oracle inflation $z_{1-\delta}\,\sigma(\beta)/\sqrt m$ that \Pt{} uses as a benchmark. What the process limit adds is the joint behaviour across levels: a practitioner who reports sets at several nominal levels, or wants a guarantee for every level in a range, needs the law of the supremum of the coverage error over the levels, which Bonferroni over a grid overstates and cannot extend to an interval, and which the weighted DKW inequalities of \PX{} bound only up to an unknown constant.

The plug-in estimator of the single-level variance $\sigma^2(\beta)$ of Theorem~\ref{thm:coverage} is
\begin{equation}\label{eq:sigmahat}
\hat\sigma_m^2(\beta):=\frac{1}{m\,\bar v^{\,2}}\sum_{i=1}^m v_i^2\bigl(\ind\{S_i\le\tauhat(\beta)\}-\beta\bigr)^2,\qquad\bar v:=\Pn v ;
\end{equation}
it is consistent uniformly in $\beta\in[\beta_0,\beta_1]$ under $K^2<\infty$ (Lemma~\ref{lem:sigmahatunif}).

\begin{corollary}[simultaneous inflation]\label{cor:simultaneous}
Assume $K^2<\infty$ and that $F_Q$ is continuous. Let $\Bcal\subset[\beta_0,\beta_1]$ be closed and $\delta\in(0,1/2)$.
\begin{enumerate}[label=(\alph*),leftmargin=*,itemsep=1pt]
\item Let $c_{1-\delta}$ be the $(1-\delta)$-quantile of $\sup_{\beta\in\Bcal}\Lb(\beta)$ and $\hat c$ the conditional $(1-\delta)$-quantile of $\sup_{\beta\in\Bcal}\hat\Lb^\xi_m(\beta)$. Then $\hat c\toP c_{1-\delta}$ and, running the weighted quantile at the levels $\beta+\hat c/\sqrt m$ (a level above one returns the whole space, $\tauhat=+\infty$),
\[
\Prob\Bigl(\exists\,\beta\in\Bcal:\ P_e\bigl(\tauhat(\beta+\hat c/\sqrt m)\bigr)>1-\beta\Bigr)\longrightarrow\delta .
\]
\item (Studentized version.) Let $\hat c^{\,s}$ be the conditional $(1-\delta)$-quantile of $\sup_{\beta\in\Bcal}\hat\Lb^\xi_m(\beta)/\hat\sigma_m(\beta)$. Then $\hat c^{\,s}$ converges in probability to the $(1-\delta)$-quantile of $\sup_{\beta\in\Bcal}\Lb(\beta)/\sigma(\beta)$, and the same conclusion holds for the level-dependent inflation $\beta+\hat c^{\,s}\hat\sigma_m(\beta)/\sqrt m$.
\item If the suprema in (a) or (b) are computed over finite sets $\Bcal_m\subset\Bcal$ whose mesh (the largest distance from a point of $\Bcal$ to $\Bcal_m$) tends to zero, the conclusions hold unchanged for the failure event over all of $\Bcal$.
\end{enumerate}
The same holds for the sets of \citet{tibshirani2019conformal} by Proposition~\ref{prop:testweight}, and with the one-sided supremum replaced by $\sup_\Bcal|\Lb|$ one obtains two-sided simultaneous bands.
\end{corollary}

Taking $\Bcal=\{\beta\}$ recovers the single-level oracle. The gain over Bonferroni for a finite grid is the correlation of $\Lb$ across levels,
\[
\Cov\bigl(\Lb(\beta),\Lb(\beta')\bigr)=\tilde G(\beta\wedge\beta')-\beta\,\tilde G(\beta')-\beta'\,\tilde G(\beta)+K^2\beta\beta',
\]
which is high for nearby levels. For an interval $\Bcal=[\beta_0,\beta_1]$ the corollary gives an asymptotically exact uniform guarantee whose cost---the $(1-\delta)$-quantile of the bounded Gaussian supremum $\sup_\Bcal\Lb$, at most $c_0K+K\sqrt{2\log(1/\delta)}$ by Corollary~\ref{cor:asympdkw}---does not grow with the number of levels. The conditional quantiles are approximated by $R$ Monte Carlo draws; the conclusions hold with $R=R_m\to\infty$, and for fixed $R$ the statistic and the draws are asymptotically exchangeable, so that at a single level the $k$-th order statistic of the draws gives limiting failure probability $1-k/(R+1)$, whence $k=\lceil(1-\delta)(R+1)\rceil$ in the procedure below. Section~\ref{sec:experiments} quantifies{} the gain over Bonferroni.

\begin{quote}\small
\textbf{Procedure 1} (joint inflation by the multiplier bootstrap). \emph{Input:} calibration pairs $(S_i,v_i)_{i\le m}$, a finite set of levels $\Bcal\subset[\beta_0,\beta_1]$ (for an interval, a grid whose mesh decreases with $m$), $\delta\in(0,1/2)$, a number of draws $R$.
\begin{enumerate}[leftmargin=*,itemsep=0pt,topsep=2pt]
\item Compute $\tauhat(\beta)$, $\Fhat(\tauhat(\beta))$ and $\hat\sigma_m(\beta)$ for $\beta\in\Bcal$.
\item For $r=1,\dots,R$ draw $\xi^{(r)}\sim N(0,I_m)$ and set $\hat\Lb^{(r)}(\beta)=(\sqrt m\,\Pn v)^{-1}\sum_i\xi^{(r)}_iv_i\bigl(\ind\{S_i\le\tauhat(\beta)\}-\Fhat(\tauhat(\beta))\bigr)$.
\item Let $\hat c$ be the $k$-th order statistic of $\max_{\beta\in\Bcal}\hat\Lb^{(r)}(\beta)/\hat\sigma_m(\beta)$ over $r=1,\dots,R$, with $k=\lceil(1-\delta)(R+1)\rceil$; if $k>R$ (that is, $R<(1-\delta)/\delta$), set $\hat c=+\infty$.
\item Report the prediction sets at the levels $\beta+\hat c\,\hat\sigma_m(\beta)/\sqrt m$, $\beta\in\Bcal$ (a level above one returns the whole space).
\end{enumerate}
By Corollary~\ref{cor:simultaneous}(b)--(c) the probability that some level in $\Bcal$ under-covers tends to $\delta$ as $m,R\to\infty$; for bounded ratios Corollary~\ref{cor:finitesimul} bounds the error at finite $m$ (for the exact conditional quantile). With estimated weights, step~2 is replaced by the pairs bootstrap of Section~\ref{sec:estboot}, which recomputes the weights.
\end{quote}

\subsection{A diagnostic for the covariate-shift assumption}
The guarantee of weighted conformal prediction rests on two unverifiable premises: that the shift is a pure covariate shift ($P_{Y|X}$ unchanged) and that the ratio $v$ is correct. If a small labeled sample $S'_1,\dots,S'_n$ from the test distribution is available, a necessary consequence of the two premises can be tested, because under them the weighted source distribution $\Fhat$ and the plain target empirical distribution
\[
\hat F^{\rm tgt}_n(t)=\frac1n\sum_{j=1}^n\ind\{S'_j\le t\}
\]
estimate the same $F_Q$.

\begin{proposition}[two-sample weighted Kolmogorov--Smirnov test]\label{prop:twosample}
Let the target sample be independent of the calibration sample, $m,n\to\infty$ with $m/(m+n)\to\lambda\in(0,1)$, and $K^2<\infty$. Under the covariate-shift premises,
\[
T_{m,n}:=\sqrt{\frac{mn}{m+n}}\,\bigl\|\Fhat-\hat F^{\rm tgt}_n\bigr\|_\infty\weak\bigl\|\sqrt{1-\lambda}\,\GQ-\sqrt{\lambda}\,\mathbb B\circ F_Q\bigr\|_\infty,
\]
with $\mathbb B$ a Brownian bridge independent of $\GQ$, and the critical value is consistently estimated by the multiplier bootstrap applied to both samples. Let $F^{\rm tgt}$ denote the law of the score under the actual test distribution. If the calibration scores are weighted by a misspecified ratio $\tilde v$ with $\E_P\tilde v=1$, or if the conditional law of $Y$ given $X$ differs between the two populations, then $\Fhat$ converges uniformly to
\[
F_{\tilde Q}(t)=\E_P\bigl[\tilde v(X)\ind\{S\le t\}\bigr]
\]
(with $\tilde v=v$ in the second case), which is computed under the \emph{training} conditional law, while $\hat F^{\rm tgt}_n\to F^{\rm tgt}$; hence $T_{m,n}\to\infty$ in probability whenever $F_{\tilde Q}\ne F^{\rm tgt}$.
\end{proposition}

The proof is immediate from Theorem~\ref{thm:fclt}, the classical Donsker theorem for the target sample, independence, and the continuous mapping theorem. The labeled target points are used only for the test, not to calibrate. Rejection refutes at least one of the two premises; failure to reject does not validate either, since a change of $P_{Y|X}$ that leaves the law of the score unchanged is invisible to the test (two outcome laws symmetric about zero give the same absolute residuals). Under a fixed alternative the statistic grows like
\[
\sqrt{\frac{mn}{m+n}}\,\bigl\|F_{\tilde Q}-F^{\rm tgt}\bigr\|_\infty ,
\]
so with $n\ll m$ the power is governed by the size of the target sample.

\section{Rates: a finite-sample form of the simultaneous guarantee}\label{sec:rates}

Everything in Section~\ref{sec:coverage} is a limit statement. This section quantifies two of them for bounded ratios: a Berry--Esseen bound at the parametric rate for the single-level central limit theorem, and a Gaussian and bootstrap approximation of the law of the maximal coverage error over a finite grid of levels with an error of order $m^{-1/4}$ up to logarithms, which turns the asymptotically exact inflation of Corollary~\ref{cor:simultaneous} into a guarantee at every finite $m$ with a slack of known order but unspecified constant. Throughout this section
\[
v\le B\qquad\text{and}\qquad F_Q\ \text{is continuous},
\]
and we recall that $B\ge1$ because $\E_Pv=1$, and that $\hat\sigma_m^2(\beta)$ of \eqref{eq:sigmahat} is the plug-in estimator of the single-level variance $\sigma^2(\beta)$.

\subsection{A single level: the parametric rate}

\begin{proposition}[Berry--Esseen bound at a fixed level]\label{prop:berryesseen}
There is a universal constant $C_0$ such that for every $\beta\in(0,1)$ and every $m\ge1$,
\[
\sup_{x\in\R}\Bigl|\Prob\bigl(\sqrt m\,\Cm(\beta)\le x\,\sigma(\beta)\bigr)-\Phi(x)\Bigr|\le\frac{C_0\,B}{\sigma(\beta)\sqrt m}.
\]
\end{proposition}

The proof (Appendix~\ref{app:rates}) does not go through the Bahadur--Kiefer representation. It uses the exact duality between the quantile and the distribution function,
\[
\bigl\{F_Q(\tauhat(\beta))<w\bigr\}=\bigl\{\Fhat(q_Q(w)-)\ge\beta\bigr\}\qquad\text{a.s., }w\in(0,1),
\]
which turns each event about the coverage into an event about a \emph{single} sum of bounded i.i.d.\ variables, to which the classical Berry--Esseen theorem applies; the only additional work is that the variance of that sum depends on $w$, and a mean-value argument controls the resulting change of scale. When the variance is of its maximal order $\sigma^2(\beta)\asymp B$, which is the adversarial coupling behind the variance proxy of \Pt{}, the bound reads $C_0\sqrt{B/m}$: the one-dimensional central limit theorem of \Pt{} (Proposition~10 there) holds with an error at the parametric rate, which is the rate at which the finite-sample bounds of that paper operate.

\subsection{A grid of levels}

Fix $0<\beta_0\le b_1<\dots<b_L\le\beta_1<1$ and write $\sigma_j:=\sigma(b_j)$, $\hat\sigma_j:=\hat\sigma_m(b_j)$. The one-sided maximal statistics, their feasible studentized versions, and the Gaussian limit are
\begin{equation}\label{eq:maxstats}
T_m:=\max_{j\le L}\frac{\sqrt m\,\Cm(b_j)}{\sigma_j},\qquad
\hat T_m:=\max_{j\le L}\frac{\sqrt m\,\Cm(b_j)}{\hat\sigma_j},\qquad
T:=\max_{j\le L}\frac{\Lb(b_j)}{\sigma_j},
\end{equation}
and the multiplier bootstrap statistic, computed with \emph{Gaussian} multipliers $\xi_i\sim N(0,1)$ in \eqref{eq:multiplier}, is
\begin{equation}\label{eq:bootmax}
\hat T^\xi_m:=\max_{j\le L}\frac{\hat\Lb^\xi_m(b_j)}{\hat\sigma_j},\qquad
\hat c_{1-\delta}:=\inf\bigl\{x:\Prob_\xi(\hat T^\xi_m\le x)\ge1-\delta\bigr\},
\end{equation}
where $\Prob_\xi$ is the conditional probability given the data. Two-sided versions, with $|\Cm(b_j)|$, $|\Lb(b_j)|$ and $|\hat\Lb^\xi_m(b_j)|$ in place of the signed quantities, are denoted by a superscript $|\cdot|$.

\begin{theorem}[Gaussian and bootstrap approximation on a grid]\label{thm:rates}
There is a constant $C$ depending only on $(\beta_0,\beta_1)$ such that, for all $m\ge2$ and $L\ge1$,
\begin{align}
\sup_{x\in\R}\bigl|\Prob(T_m\le x)-\Prob(T\le x)\bigr|&\le C\,B\,m^{-1/4}\log^{5/4}(Lm),\label{eq:rate-gauss}\\
\sup_{\delta\in(0,1)}\bigl|\Prob(\hat T_m>\hat c_{1-\delta})-\delta\bigr|&\le C\,B\,m^{-1/4}\log^{5/4}(Lm),\label{eq:rate-boot}
\end{align}
and the same bounds hold for the two-sided statistics $T_m^{|\cdot|}$, $T^{|\cdot|}$, $\hat T_m^{|\cdot|}$ and $\hat c^{\,|\cdot|}_{1-\delta}$.
\end{theorem}

The proof (Appendix~\ref{app:rates}) has three ingredients. The exact representation of Theorem~\ref{thm:coverage} and a quantitative modulus of continuity for the weighted score process reduce $\sqrt m\,\Cm(b_j)$ to the plain sum $\Gm h_{q_Q(b_j)}$ up to an error of order $B^{3/4}m^{-1/4}\log^{3/4}m+B\log(Lm)/\sqrt m$ with probability $1-O(1/m)$: this is a weighted Bahadur--Kiefer bound, and the exponent $1/4$ is the one of \citet{kiefer1967bahadur} for the unweighted sample quantile process, where it is exact up to iterated logarithms. The high-dimensional central limit theorem of \citet{chernozhukov2022improved} approximates the law of the maximum of the $L$-dimensional vector sum by that of the Gaussian maximum with an error of order $(B^2\log^5(Lm)/m)^{1/4}$, uniformly in $L$; its bootstrap counterpart, in the same paper, handles the multiplier critical value. Nazarov's anti-concentration inequality for Gaussian maxima converts the Bahadur--Kiefer error into an error in probability. The logarithmic factor is inherited from the high-dimensional central limit theorem and the Bahadur--Kiefer step, and the dependence on the number $L$ of levels is only through $\log L$: the guarantee below costs nothing in rate when the grid is refined.

The comparison between Proposition~\ref{prop:berryesseen} and Theorem~\ref{thm:rates} is the quantitative form of the message of Section~\ref{sec:coverage}: at a single level the duality gives the parametric rate because no uniformity in the threshold is needed, whereas the joint law of several levels requires $\Zm$ at the \emph{random} points $\tauhat(b_j)$, and the Bahadur--Kiefer remainder is of exact order $m^{-1/4}$, up to iterated logarithms, already for $v\equiv1$ \citep{kiefer1967bahadur}. Whether the grid maximum can nevertheless be approximated beyond $m^{-1/4}$ is open.

\begin{corollary}[simultaneous inflation at finite $m$]\label{cor:finitesimul}
In the setting of Theorem~\ref{thm:rates}, let $\delta\in(0,1/2)$ and let the weighted quantile be run at the inflated levels
\[
b_j'=b_j+\frac{\hat c_{1-\delta}\,\hat\sigma_j}{\sqrt m},\qquad j=1,\dots,L .
\]
Then, with $C$ depending only on $(\beta_0,\beta_1)$,
\[
\Prob\Bigl(\exists\,j\le L:\ P_e\bigl(\tauhat(b_j')\bigr)>1-b_j\Bigr)\le\delta+C\,B\,m^{-1/4}\log^{5/4}(Lm),
\]
and the same bound holds for the sets of \citet{tibshirani2019conformal}, whose test-point weight can only raise the level.
\end{corollary}

Corollary~\ref{cor:finitesimul} is to Corollary~\ref{cor:simultaneous} what Proposition~\ref{prop:bousquet} is to Corollary~\ref{cor:asympdkw}: the same statement with an explicit finite-sample slack. Its critical value is the exact Gaussian one up to the bootstrap error, so for a fine grid it is smaller than the Bonferroni value by the factor documented in Section~\ref{sec:experiments}{} for the unstudentized version (the studentized critical values behave in the same way), and smaller than the DKW-type inflation $B\sqrt{2\log(4/\delta)/m}$ of \PX{} by a factor of order $\sqrt B$. The constant $C$ is not explicit---it contains the constants of \citet{chernozhukov2022improved} and of the bracketing maximal inequality---so the corollary is a rate statement, not a certificate in the sense of \Pt{}; a certificate with explicit constants that is uniform over levels and operates at the scale $\sigma(\beta)$ rather than $B$ remains open.%

\section{Estimated likelihood ratios}\label{sec:estimated}

In practice $v$ is estimated, typically by a probabilistic classifier of ``target versus source'' covariates trained on unlabeled covariate samples from both populations, and the source sample used for this purpose usually \emph{includes} the calibration covariates. \Pt{} treats a fixed estimate $\hat v$ independent of the calibration sample, through an oracle inequality whose statistical term is evaluated at the pseudo-target $\hat Q=(\hat v/\E_P\hat v)\cdot P$ and whose bias is $\sup_t|F_Q(t)-F_{\hat Q}(t)|$, and notes that a dependent $\hat v$ ``would require uniformity over the class of possible $\hat v$, i.e.\ exactly the empirical-process machinery'' its argument avoids. This section supplies that machinery.

Let $\hat v\ge0$ be a measurable function of the covariate, possibly depending on the calibration sample itself (in practice through the covariates $X_1,\dots,X_m$) and on auxiliary data (unlabeled source and target covariates). Define the estimated-weight distribution function, its quantile and the pseudo-target
\begin{equation}\label{eq:Fchk}
\Fchk(t)=\frac{\sum_i\hat v(X_i)\ind\{S_i\le t\}}{\sum_i\hat v(X_i)},\qquad
\tauchk(\beta)=\inf\{t:\Fchk(t)\ge\beta\},\qquad
F_{\hat Q}(t)=\frac{P[\hat v\ind\{S\le t\}]}{P\hat v},
\end{equation}
where $P[\cdot]$ integrates over $Z\sim P$ with $\hat v$ held fixed (so $F_{\hat Q}$ is random when $\hat v$ is).

\begin{assumption}[Donsker weight class]\label{ass:donsker}
There is a class $\Vcal$ of nonnegative measurable functions on $\mathcal X$ with $v\in\Vcal$ and an envelope $V\ge\sup_{u\in\Vcal}u$ with $\|V\|_{P,2}<\infty$, such that
\[
\Hcal_\Vcal:=\bigl\{z\mapsto u(x)\bigl(\ind\{s(z)\le t\}-c\bigr):\ u\in\Vcal,\ t\in\Rbar,\ c\in[0,1]\bigr\}
\]
is $P$-Donsker; and $\Prob(\hat v\in\Vcal)\to1$ and
\[
\|\hat v-v\|_{P,2}^2=\int\bigl(\hat v(x)-v(x)\bigr)^2\,\dd P_X(x)\toP0,
\]
the norm being computed for the realized $\hat v$.
\end{assumption}

\begin{theorem}[uniform expansion with estimated weights]\label{thm:estimated}
Under Assumption~\ref{ass:donsker},
\[
\sup_{t\in\Rbar}\Bigl|\sqrt m\bigl(\Fchk(t)-F_{\hat Q}(t)\bigr)-\Zm(t)\Bigr|\toP0,
\]
that is, uniformly in $t$,
\[
\sqrt m\,(\Fchk-F_Q)=\Zm+\sqrt m\,(F_{\hat Q}-F_Q)+o_p(1).
\]
\end{theorem}

The theorem (Appendix~\ref{app:estimated}) says that, centered at the pseudo-target, the estimated-weight process is asymptotically the oracle process, however $\hat v$ was obtained inside $\Vcal$. The proof is the identity of Lemma~\ref{lem:identity} for $\hat v$,
\[
P\bigl[\hat v\,(\ind\{S\le t\}-F_{\hat Q}(t))\bigr]=0\qquad\text{exactly},
\]
together with the asymptotic equicontinuity of $\Gm$ over $\Hcal_\Vcal$, which is insensitive to the dependence between $\hat v$ and the sample. The bias $F_{\hat Q}-F_Q$ is the quantity bounded by $\Delta$ in Theorem~12 of \Pt{}; here its exact first-order contribution is identified. Assumption~\ref{ass:donsker} holds for the usual parametric ratio models:

\begin{lemma}[parametric classes]\label{lem:parametric}
Let $\Theta\subset\R^d$ be bounded and $\Vcal=\{v_\theta:\theta\in\Theta\}$ with $v_\theta\ge0$,
\[
|v_\theta(x)-v_{\theta'}(x)|\le\dot V(x)\,\|\theta-\theta'\|,\qquad V=\sup_\theta v_\theta,\qquad V,\dot V\in L_2(P).
\]
Then $N_{[\,]}(\epsilon,\Hcal_\Vcal,L_2(P))\le C\epsilon^{-(d+3)}$ for $\epsilon\le1$ and a constant $C$ depending on $d$, $\operatorname{diam}\Theta$, $\|V\|_{P,2}$ and $\|\dot V\|_{P,2}$; in particular $\Hcal_\Vcal$ is $P$-Donsker.
\end{lemma}

Exponential tilts $v_\theta(x)=\exp(\theta^\top x-\psi(\theta))$ when $\E_Pe^{c\|X\|}<\infty$ for every $c$, and the ratios
\[
v_\theta(x)=\frac{(1-\pi)\,\sigma(\theta^\top\phi(x))}{\pi\,\bigl(1-\sigma(\theta^\top\phi(x))\bigr)}
\]
obtained from a logistic ``target versus source'' classifier with bounded features $\phi$ and target fraction $\pi$ in the pooled sample, satisfy the lemma on a bounded parameter set.

\subsection{Parametric ratio models: the coupling gradient}

Suppose $\hat v=v_{\hat\theta}$ for a parametric model with $v_{\theta_0}=v$, and that $\theta\mapsto v_\theta$ is Fr\'echet differentiable in $L_1(P)$ at $\theta_0$ with derivative $\dot v\in L_1(P)^d$:
\[
\bigl\|v_\theta-v_{\theta_0}-(\theta-\theta_0)^\top\dot v\bigr\|_{P,1}=o(\|\theta-\theta_0\|).
\]
Define the \emph{coupling gradient}
\begin{equation}\label{eq:gamma}
\gamma(t):=P\bigl[\dot v\,(\ind\{S\le t\}-F_Q(t))\bigr]=\E_P\bigl[\dot v(X)\{F_{S|X}(t)-F_Q(t)\}\bigr]\in\R^d .
\end{equation}
Its name reflects that the estimation error of the ratio affects the weighted distribution function only through the covariance between the score function of the weight model and the conditional distribution of the score given the covariate; a direction of parameter error along which $F_{S|X}(t)$ does not vary is harmless, and $\gamma\equiv0$ when the score is independent of the covariate (in which case $F_Q=F_P$ and no weighting is needed). This is the exact first-order version of the observation in \Pt{} that only the score-projected part of the weight error enters the bias. For the estimator we allow a general influence-function expansion in which the calibration covariates may appear:
\begin{equation}\label{eq:thetaexp}
\sqrt m\,(\hat\theta-\theta_0)=\Gm\psi+\xi_m+o_p(1),\qquad P\psi=0,\quad\psi=\psi(X)\in L_2(P_X)^d,
\end{equation}
where $\xi_m$ is independent of the calibration sample and $\xi_m\weak N(0,\Sigma_{\rm aux})$; the influence function depends on the covariate only, as it does for every ratio estimator fitted from covariates, and this is what gives the covariance formula \eqref{eq:covpsi} below its conditional form (for a general $\psi(X,S)$ one would keep $P[h_t\psi]$). If $\hat\theta$ is computed from an unlabeled source sample of size $N_s$ and an unlabeled target sample of size $N_t$, both independent of the calibration data, then $\psi\equiv0$ and $\xi_m$ collects the two auxiliary terms; if the calibration covariates are the source sample (or part of it), $\psi$ is their influence function.

\begin{theorem}[parametric estimated weights]\label{thm:parametric}
Let Assumption~\ref{ass:donsker} hold with $\Vcal=\{v_\theta\}$, let $\theta\mapsto v_\theta$ be $L_1(P)$-differentiable at $\theta_0$, and let \eqref{eq:thetaexp} hold. Then, in $\linf(\Rbar)$,
\[
\sqrt m\,(\Fchk-F_Q)\weak\GQ+\gamma^\top\bigl(\mathbb G\psi+\mathcal N\bigr),
\]
where $(\GQ,\mathbb G\psi)$ is the jointly Gaussian limit of $(\Gm h_\cdot,\Gm\psi)$, with
\begin{equation}\label{eq:covpsi}
\Cov\bigl(\GQ(t),\mathbb G\psi\bigr)=P[h_t\psi]=\E_Q\bigl[\{F_{S|X}(t)-F_Q(t)\}\,\psi(X)\bigr],
\end{equation}
and $\mathcal N\sim N(0,\Sigma_{\rm aux})$ is independent of $(\GQ,\mathbb G\psi)$.
\end{theorem}

\begin{corollary}[coverage with estimated weights]\label{cor:estcoverage}
In the setting of Theorem~\ref{thm:parametric}, let $F_P$ and $F_Q$ be continuous, let $F_Q$ be strictly increasing in a neighbourhood of $q=q_Q(\beta)$, and let $\gamma$ be continuous at $q$. Then
\[
\sqrt m\bigl(P_e(\tauchk(\beta))-(1-\beta)\bigr)\weak N\bigl(0,\tau^2(\beta)\bigr),
\]
where
\[
\tau^2(\beta)=\sigma^2(\beta)+2\,\gamma(q)^\top c(q)+\gamma(q)^\top\bigl(P\psi\psi^\top+\Sigma_{\rm aux}\bigr)\gamma(q),
\qquad
c(q)=\E_Q\bigl[\{F_{S|X}(q)-\beta\}\,\psi(X)\bigr],
\]
with $\sigma^2(\beta)$ as in Theorem~\ref{thm:coverage}.
\end{corollary}

The cross term $2\gamma(q)^\top c(q)$ is the new element: it is absent when the ratio is estimated from data independent of the calibration sample, and it can be negative. In that case using the calibration covariates in the ratio estimate lowers the asymptotic variance of the training-conditional coverage, possibly below that of the oracle weights. The mechanism is the one of \citet{hirano2003efficient} for inverse-probability weighting: the estimated weights absorb the part of the sampling fluctuation of the calibration covariates that is aligned with the model's score function.

\begin{example}[Gaussian mean shift]\label{ex:meanshift}
Let $P_X=N(0,1)$ and $Q_X=N(\theta_0,1)$, so that the ratio is the exponential tilt
\[
v_\theta(x)=\exp\bigl(\theta x-\theta^2/2\bigr),\qquad \dot v(x)=v(x)\,(x-\theta_0),\qquad \gamma(t)=\Cov_Q\bigl(X,F_{S|X}(t)\bigr)
\]
(the normalizing constant cancels in $\Fchk$). Estimate $\theta$ by $\hat\theta=\bar X^{\rm tgt}_{N_t}-\bar X^{\rm src}_{N_s}$ with $m/N_t\to\lambda_t$ and $m/N_s\to\lambda_s$.
\begin{enumerate}[leftmargin=*,itemsep=1pt]
\item[(a)] If the source sample is independent of the calibration data: $\psi\equiv0$, $\Sigma_{\rm aux}=\lambda_t+\lambda_s$, and
\[
\tau^2(\beta)=\sigma^2(\beta)+\gamma(q)^2\,(\lambda_t+\lambda_s).
\]
\item[(b)] If the calibration covariates are the source sample ($N_s=m$): $\psi(x)=-x$, $c(q)=-\gamma(q)$, $\Sigma_{\rm aux}=\lambda_t$, and
\[
\tau^2(\beta)=\sigma^2(\beta)+\gamma(q)^2\,(\lambda_t-1).
\]
\end{enumerate}
In case (b) the estimated weights reproduce the oracle variance when $\lambda_t=1$ and beat it when $\lambda_t<1$ (a target sample asymptotically larger than the calibration sample); in case (a) the estimation always costs variance. Section~\ref{sec:experiments}{} confirms the three predictions.
\end{example}

\subsection{Feasible inference with estimated weights}\label{sec:estboot}

Corollary~\ref{cor:estcoverage} gives the asymptotic variance $\tau^2(\beta)$, but its ingredients $\gamma$, $\psi$ and $\Sigma_{\rm aux}$ depend on the ratio model and on how the auxiliary data were used, and the bootstrap of Section~\ref{sec:bootstrap} does not reproduce them: if the estimated weights are held fixed and only the calibration scores are resampled or reweighted, the bootstrap process converges to the oracle limit $\GQ$ and misses the term $\gamma^\top(\mathbb G\psi+\mathcal N)$ together with its covariance with $\GQ$. The remedy is to resample everything that was used to compute the weights and to recompute them. Let $\mathcal A$ denote the auxiliary unlabeled samples (target covariates $X'_1,\dots,X'_{N_t}$, and an independent source sample if one is used), so that $\hat\theta=\hat\theta(X_{1:m},\mathcal A)$, and let $\hat v=v_{\hat\theta}$. A \emph{pairs bootstrap} draws $(X_i^*,S_i^*)$, $i\le m$, with replacement from the calibration pairs and, independently, a bootstrap resample $\mathcal A^*$ of each auxiliary sample, recomputes $\hat\theta^*=\hat\theta(X^*_{1:m},\mathcal A^*)$ and $\hat v^*_i=v_{\hat\theta^*}(X_i^*)$, and forms
\begin{equation}\label{eq:pairsboot}
\Fchk^*(t)=\frac{\sum_i\hat v^*_i\ind\{S^*_i\le t\}}{\sum_i\hat v^*_i},\qquad
\Lambda^*_m(\beta):=\sqrt m\bigl(\Fchk^*(\tauchk(\beta))-\Fchk(\tauchk(\beta))\bigr).
\end{equation}
Nothing about $F_Q$, $\gamma$ or $\psi$ needs to be known or estimated. Two regularity conditions on the ratio model and its estimator are needed beyond Theorem~\ref{thm:parametric}.

\begin{assumption}[bootstrappable ratio model]\label{ass:boot}
\begin{enumerate}[label=(\roman*),leftmargin=*,itemsep=1pt]
\item\label{ass:boot-smooth} The class $\Vcal=\{v_\theta:\theta\in\Theta\}$ satisfies the conditions of Lemma~\ref{lem:parametric}, $\theta_0$ is an interior point of $\Theta$, and on a neighbourhood $U\subset\Theta$ of $\theta_0$ the map $\theta\mapsto v_\theta(x)$ is differentiable, with derivative $\dot v_\theta(x)$ such that
\[
\bigl|v_\theta(x)-v_{\theta'}(x)-(\theta-\theta')^\top\dot v_{\theta'}(x)\bigr|\le\ddot H(x)\,\|\theta-\theta'\|^2,\qquad \ddot H\in L_1(P),
\]
$\|\dot v_\theta-\dot v_{\theta_0}\|_{P,1}\to0$ as $\theta\to\theta_0$, and the class $\{\dot v_\theta\ind\{S\le t\}\}$ over $\theta\in U$ and $t\in\Rbar$ is $P$-Glivenko--Cantelli. (The second-order bound implies the $L_1(P)$-differentiability required in Theorem~\ref{thm:parametric}, with $\dot v=\dot v_{\theta_0}$.)
\item\label{ass:boot-est} Conditionally on the data, in probability, the bootstrap estimator admits the expansion $\sqrt m(\hat\theta^*-\hat\theta)=\mathbb G_m^*\psi+\xi^*_m+o_{p^*}(1)$, where $\mathbb G^*_m=\sqrt m(\Pn^*-\Pn)$ is the bootstrap empirical process of the calibration pairs and $\xi^*_m$, computed from $\mathcal A^*$ alone, satisfies $\xi^*_m\weak_*N(0,\Sigma_{\rm aux})$.
\end{enumerate}
\end{assumption}

Condition~\ref{ass:boot-est} is more than the consistency of the bootstrap for $\hat\theta$ alone: it is a bootstrap influence-function expansion \emph{jointly} with the calibration empirical process, the term $\mathbb G^*_m\psi$ being computed from the same resampled pairs as $\Fchk^*$, which is what preserves the covariance between the weight estimate and the scores in the limit. It is immediate for estimators that are smooth functions of sample means, such as the moment estimator of the exponential tilt in Example~\ref{ex:meanshift}, and it holds under the standard conditions for bootstrapping $M$-estimators \citep[Chapter~10]{kosorok2008introduction} for the maximum likelihood estimator of a logistic ``target versus source'' classifier, the two covariate samples being resampled separately. Condition~\ref{ass:boot-smooth} holds for exponential tilts with $\E_Pe^{c\|X\|}<\infty$ and for bounded-feature logistic ratios.

\begin{theorem}[pairs bootstrap with estimated weights]\label{thm:estboot}
Let the assumptions of Theorem~\ref{thm:parametric} and Assumption~\ref{ass:boot} hold, let $F_P$ and $F_Q$ be continuous, and let $\gamma$ be continuous on $\R$ (which holds whenever $F_{S|X=x}$ is continuous for $P_X$-almost every $x$). Then, in $\linf[\beta_0,\beta_1]$,
\[
\sqrt m\bigl(\beta-F_Q(\tauchk(\beta))\bigr)\weak\Lb^{\rm est}(\beta):=\Lb(\beta)+\gamma(q_Q(\beta))^\top\bigl(\mathbb G\psi+\mathcal N\bigr),
\qquad\text{and}\qquad
\Lambda^*_m\weak_*\Lb^{\rm est},
\]
the second convergence being conditional on the data, in probability. Consequently, if $\tau^2(\beta)>0$ on a closed $\Bcal\subset[\beta_0,\beta_1]$ (for which $\E_P[v^2F_{S|X}(q)\{1-F_{S|X}(q)\}]>0$ at $q=q_Q(\beta)$ suffices), $\hat c^*$ is the conditional $(1-\delta)$-quantile of $\sup_{\beta\in\Bcal}\Lambda^*_m(\beta)$ and $\delta\in(0,1/2)$, then
\[
\Prob\Bigl(\exists\,\beta\in\Bcal:\ P_e\bigl(\tauchk(\beta+\hat c^*/\sqrt m)\bigr)>1-\beta\Bigr)\longrightarrow\delta ,
\]
whereas the bootstrap of Theorem~\ref{thm:bootstrap} applied with $\hat v$ held fixed converges to $\Lb$ and calibrates to $\sigma^2(\beta)$ in place of $\tau^2(\beta)$. At a single level its limiting failure probability is $1-\Phi(z_{1-\delta}\,\sigma(\beta)/\tau(\beta))$: it is anti-conservative when the ratio is estimated from data independent of the calibration sample and $\gamma(q)^\top\Sigma_{\rm aux}\gamma(q)>0$ (then $\tau^2=\sigma^2+\gamma(q)^\top\Sigma_{\rm aux}\gamma(q)$; equality is possible with $\gamma(q)\ne0$ if $\Sigma_{\rm aux}$ is singular or the auxiliary samples grow faster than $m$), and conservative in case~(b) of Example~\ref{ex:meanshift} when $\lambda_t<1$ and $\gamma(q)\neq0$.
\end{theorem}

The proof (Appendix~\ref{app:estimated}) has the structure of Theorems~\ref{thm:estimated} and~\ref{thm:parametric} with the empirical measure in the role of the population: the identity of Lemma~\ref{lem:identity} holds exactly for the bootstrap sample when $\Fchk^*$ is centered at the $\hat v^*$-tilt of $\Pn$, the bootstrap empirical process is asymptotically equicontinuous over the Donsker class $\Hcal_\Vcal$ \citep[Theorem~3.6.1]{vandervaart1996weak}, and the difference between the $\hat v^*$-tilt and the $\hat v$-tilt of $\Pn$ is $\gamma^\top(\hat\theta^*-\hat\theta)$ to first order by condition~\ref{ass:boot-smooth}. Section~\ref{sec:experiments} shows{} the two bootstraps side by side.

\begin{remark}[what the fixed-level argument cannot see]
The decomposition behind Theorem~12 of \Pt{} is
\[
P_e(\tauchk)\le\bigl\{1-F_{\hat Q}(\tauchk)\bigr\}+\Delta,\qquad\Delta=\|F_Q-F_{\hat Q}\|_\infty\le\tfrac12\,\|v-\tilde v\|_{P,1},
\]
with $\tilde v=\hat v/\E_P\hat v$: a worst-case bias plus a statistical term evaluated under $\hat Q$. Theorem~\ref{thm:estimated} shows that this statistical term is asymptotically exact even when $\hat v$ depends on the calibration covariates, and Corollary~\ref{cor:estcoverage} shows that the bias term, once $\hat v$ is $\sqrt m$-consistent, is not a worst case but a Gaussian fluctuation correlated with the statistical term. A finite-sample PAC certificate that exploits the negative cross term is an open problem; the asymptotic statement is what the process limit delivers.
\end{remark}

\section{Is the weighted conformal threshold efficient?}\label{sec:efficiency}

Corollary~\ref{cor:estcoverage} and Example~\ref{ex:meanshift} showed that estimating the likelihood ratio from a sample that includes the calibration covariates can give a \emph{smaller} asymptotic variance than using the true ratio. A phenomenon of this kind signals that the oracle procedure is not efficient in the model in which the estimated procedure is evaluated \citep{hirano2003efficient}. This section makes the statement precise for the training-conditional coverage of a weighted quantile, in three models that differ by what is known about the shift, and identifies the efficient threshold when the ratio is known. The results are stated for the score $S$ and the functional
\[
\theta_t(P):=F_Q(t)=\E_P\bigl[v(X)\ind\{S\le t\}\bigr],
\]
whose estimation is what calibration at a threshold amounts to; the transfer to the quantile and to the coverage is the delta method.

\subsection{Three models}

Let $m_t(x):=\Prob(S\le t\mid X=x)$, so that $F_Q(t)=\E_Q[m_t(X)]$ and $G(t)=\E_P[v^2m_t(X)]$.
\begin{itemize}[leftmargin=*,itemsep=2pt]
\item $\mathcal M_1$ (\emph{known ratio}). The calibration sample is i.i.d.\ from $P=P_X\otimes P_{S|X}$, where $P_{S|X}$ is unrestricted and $P_X$ ranges over the laws with
\[
\int v\,\dd P_X=1,\qquad\int v^2\,\dd P_X<\infty ,
\]
the function $v$ being known: the model of weighted conformal prediction with an oracle ratio, in which $P_X$ is unknown but the constraint $\int v\,\dd P_X=1$ is known.
\item $\mathcal M_2$ (\emph{known marginals}). $P_X$ and $Q_X$ are both known, and $P_{S|X}$ is unrestricted.
\item $\mathcal M_3$ (\emph{unknown ratio, two samples}). Nothing is known about $P_X$, $Q_X$ or $P_{S|X}$; besides the $m$ labeled source points, $N$ unlabeled target covariates $X_1',\dots,X_N'\sim Q_X$ are observed, with $m/(m+N)\to\pi\in(0,1)$.
\end{itemize}
Throughout, $K^2>1$, that is, $v$ is not almost surely constant; otherwise there is no shift.

\subsection{The known-ratio model}

\begin{theorem}[efficiency bound with a known ratio]\label{thm:eff}
In $\mathcal M_1$, the functional $\theta_t$ is pathwise differentiable at every $P$ with $K^2<\infty$, with efficient influence function
\[
\psi_t:=v\,\ind\{S\le t\}-F_Q(t)-c_t\,(v-1),\qquad c_t:=\frac{G(t)-F_Q(t)}{K^2-1},
\]
and efficiency bound
\[
V_1(t):=\E_P\psi_t^2=G(t)-F_Q(t)^2-\frac{\bigl(G(t)-F_Q(t)\bigr)^2}{K^2-1}.
\]
The self-normalized weighted empirical distribution function $\Fhat(t)$ is regular, with influence function $h_t$ and asymptotic variance
\[
C(t,t)=V_1(t)+\frac{\bigl(K^2F_Q(t)-G(t)\bigr)^2}{K^2-1};
\]
it is efficient at $t$ if and only if $G(t)=K^2F_Q(t)$, that is, if and only if $v(X)$ and $\ind\{S\le t\}$ are uncorrelated under $Q$.
\end{theorem}

The proof (Appendix~\ref{app:efficiency}) is the usual projection argument \citep{bickel1993efficient,vandervaart1998asymptotic}. The tangent space of $\mathcal M_1$ is the sum of the conditional tangent space $\{b(S,X):\E[b\mid X]=0\}$ and of the marginal tangent space $\{a(X):\E a=0,\ \E[va]=0\}$; the second constraint is the derivative of $\int v\,\dd P_X=1$. The influence function $v\ind\{S\le t\}-F_Q(t)$ of the unnormalized weighted mean has an $X$-component $vm_t(X)-F_Q(t)$ that is not orthogonal to $v-1$, and the projection removes the component along $v-1$. Self-normalization, $\Fhat=\Pn(v\ind\{S\le t\})/\Pn v$, subtracts $F_Q(t)(v-1)$ instead of $c_t(v-1)$ from the influence function: it is a control-variate correction with the coefficient $F_Q(t)$, which is the right one exactly when $c_t=F_Q(t)$. The unnormalized, the self-normalized and the efficient estimators are thus the control-variate estimators with coefficients $0$, $F_Q(t)$ and $c_t$, and the second is efficient only by coincidence.

Control variates built from the weights are classical in importance sampling: the regression estimator of \citet{hesterberg1995weighted}, which corrects the weighted mean by a least-squares multiple of $\Pn v-1$, and the control variates of \citet{owen2000safe} are estimators of this form, with a known variance reduction over the self-normalized estimator. What Theorem~\ref{thm:eff} adds is the model: in $\mathcal M_1$ the regression estimator with coefficient $c_t$ is not only better but efficient, and the self-normalized estimator is efficient exactly when $G=K^2F_Q$; Theorem~\ref{thm:effthreshold} transfers this to the threshold and the coverage, where the correction has a price that importance sampling does not face, the loss of the finite-sample marginal guarantee.

The efficient estimator is computable, and the corresponding threshold has an exact process limit.

\begin{theorem}[the efficient threshold]\label{thm:effthreshold}
Assume $K^2<\infty$, $K^2>1$ and that $F_Q$ is continuous. Let
\begin{gather*}
\hat c_m(t):=\frac{\Pn(v^2\ind\{S\le t\})-\Fhat(t)}{\Pn v^2-1},\qquad
\tilde F_m(t):=\Pn\bigl(v\ind\{S\le t\}\bigr)-\hat c_m(t)\,(\Pn v-1),\\
\tilde\tau(\beta):=\inf\{t:\tilde F_m(t)\ge\beta\},
\end{gather*}
with the conventions $\hat c_m\equiv0$ on the event $\Pn v^2=1$ and $\tilde\tau(\beta)=+\infty$ when the infimum is over the empty set (events of vanishing probability). Then $\sqrt m(\tilde F_m-F_Q)$ is asymptotically linear with influence function $\psi_t$, uniformly in $t$, and, in $\linf[\beta_0,\beta_1]$,
\[
\sqrt m\,\bigl(\beta-F_Q(\tilde\tau(\beta))\bigr)\weak\tilde\Lb(\beta):=\Lb(\beta)+\bigl(\beta-\tilde c(\beta)\bigr)\mathcal N,
\]
where $\tilde c(\beta)=c_{q_Q(\beta)}$ and $\mathcal N$ is a centered Gaussian variable with
\[
\Var\mathcal N=K^2-1,\qquad\Cov\bigl(\Lb(\beta),\mathcal N\bigr)=\tilde G(\beta)-K^2\beta .
\]
In particular
\[
\Var\tilde\Lb(\beta)=\sigma^2_{\rm eff}(\beta):=\sigma^2(\beta)-\frac{\bigl(\beta K^2-\tilde G(\beta)\bigr)^2}{K^2-1}\le\sigma^2(\beta),
\]
and for every regular estimator $\tau_m(\beta)$ of $q_Q(\beta)$ in $\mathcal M_1$, at a point where $F_Q$ has a positive derivative, the asymptotic variance of $\sqrt m(\beta-F_Q(\tau_m(\beta)))$ is at least $\sigma^2_{\rm eff}(\beta)$.
\end{theorem}

The gain $(\beta K^2-\tilde G(\beta))^2/(K^2-1)$ is the squared covariance $\Cov_Q(v,\ind\{S\le q_Q(\beta)\})^2$ divided by $\Var_P(v)$: it vanishes when the shift is unrelated to the size of the score and is large when the test distribution concentrates where scores are large or small. The price is finite-sample validity: $\tilde F_m$ is not monotone and $\tilde\tau$ is not a weighted quantile, so the marginal guarantee of \citet{tibshirani2019conformal} is lost, and only the asymptotic training-conditional statement remains. At a level with $\sigma^2_{\rm eff}(\beta)>0$, the efficient threshold can be run at the inflated level $\beta+z_{1-\delta}\hat\sigma_{\rm eff}(\beta)/\sqrt m$ as in Section~\ref{sec:coverage}, where
\[
\hat\sigma^2_{\rm eff}(\beta):=\Pn\bigl\{v\,\ind\{S\le\tilde\tau(\beta)\}-\beta-\hat c_m(\tilde\tau(\beta))\,(v-1)\bigr\}^2\ge0
\]
is consistent for $\sigma^2_{\rm eff}(\beta)$; the asymptotic inflation is no larger than the conformal one, and strictly smaller when $\tilde G(\beta)\ne\beta K^2$. The condition $\sigma^2_{\rm eff}(\beta)>0$ cannot be dropped: $\sigma^2_{\rm eff}(\beta)=0$ exactly when $v$ is almost surely constant on each of the events $\{S\le q_Q(\beta)\}$ and $\{S>q_Q(\beta)\}$; the limit $\tilde\Lb(\beta)$ is then degenerate and the inflation carries no guarantee (Proposition~\ref{prop:sigmaeff} and Example~\ref{ex:degenerate}). Section~\ref{sec:experiments} quantifies{} the gain in the examples of this paper; it is modest when the ratio is bounded and the score is only mildly related to the shift, and substantial under a strong tilt.

\subsection{Known marginals and unknown ratio}\label{sec:eff23}
\begin{proposition}[the other two bounds]\label{prop:bounds23}
\begin{enumerate}[label=(\alph*),leftmargin=*,itemsep=1pt]
\item In $\mathcal M_2$ the efficient influence function of $\theta_t$ is $v(X)\{\ind\{S\le t\}-m_t(X)\}$ and the bound is
\[
V_2(t):=\E_P\bigl[v^2m_t(X)\{1-m_t(X)\}\bigr]\le V_1(t).
\]
\item In $\mathcal M_3$, formulated as a pooled i.i.d.\ sample of $m+N$ observations with source indicator $R\sim\mathrm{Bernoulli}(\pi)$ (given $R=1$ the pair $(X,S)$ is drawn from $P$, given $R=0$ only $X\sim Q_X$ is observed), the efficient influence function of $\theta_t$ is
\[
\frac{R}{\pi}\,v(X)\bigl\{\ind\{S\le t\}-m_t(X)\bigr\}+\frac{1-R}{1-\pi}\bigl\{m_t(X)-F_Q(t)\bigr\},
\]
\citep[Lemma~1]{yang2024doubly}; its second moment is the bound in the normalization by $\sqrt{m+N}$, and in the normalization by $\sqrt m$ used in this paper the bound is
\[
V_3(t):=\E_P\bigl[v^2m_t(X)\{1-m_t(X)\}\bigr]+\frac{m}{N}\,\Var_Q\bigl(m_t(X)\bigr)\ \longrightarrow\ V_2(t)\quad(N/m\to\infty).
\]
\end{enumerate}
\end{proposition}

Part~(a) is the projection of Theorem~\ref{thm:eff} with the marginal tangent space removed, and part~(b) is the classical bound for a missing-at-random mean \citep{hahn1998role}, which \citet{yang2024doubly} transported to prediction sets under covariate shift and used to build doubly robust calibration procedures. Three remarks place these bounds next to the conformal procedures.

First, the excess of the oracle weighted conformal variance over the known-marginals bound is
\[
C(t,t)-V_2(t)=\E_P\bigl[v^2\{m_t(X)-F_Q(t)\}^2\bigr],
\]
the weighted variance of the regression of the coverage indicator on the covariates; the control variate of Theorem~\ref{thm:effthreshold} recovers the part of it that is explained linearly by $v$. The rest is recovered either by an estimate of $m_t$, that is, a model for the conditional distribution of the score given the covariates, as in \citet{yang2024doubly}, or, when the ratio is unknown and estimated nonparametrically, by the estimation itself: Theorem~7 of \citet{chen2008semiparametric} shows that inverse probability weighting with a nonparametrically estimated ratio attains the two-sample bound $V_3$ at a fixed $t$, under smoothness and boundedness conditions on the ratio, the analogue for the two-sample setting of \citet{hirano2003efficient}.

Second, Theorem~\ref{thm:parametric} and Example~\ref{ex:meanshift} are about a parametric submodel of $\mathcal M_3$, and are the parametric, uniform-in-level form of the mechanism behind the results just cited. In a parametric model the bound is smaller than $V_3$, and nothing prevents the variance of case~(b) of the example from falling below $V_3$; what the example illustrates is the mechanism---the oracle estimator ignores the unlabeled target sample, and any estimator that uses it correctly can do better, as $V_3<C(t,t)$ whenever $(m/N)\Var_Q(m_t)<\E_P[v^2(m_t-F_Q(t))^2]$. The limit $V_3\to V_2$ as $N/m\to\infty$ says that with abundant unlabeled target data, knowing the marginals or not makes no difference to the bound, and that the gap to the oracle conformal procedure is then exactly the regression variance above.

Third, in $\mathcal M_1$ the efficient threshold uses nothing beyond what weighted conformal prediction uses: the calibration scores, the weights, and the fact that the weights average to one. It is the only one of the three bounds attained without estimating a nuisance function: $V_2$ needs a model for $m_t$ and $V_3$ a nonparametric estimate of the ratio.

\section{Beyond finite $\chi^2$: stable limits}\label{sec:stable}

When $\chi^2(Q\,\|\,P)=\infty$, Proposition~\ref{prop:necessity} rules out the $\sqrt m$ rate, but weighted conformal prediction remains consistent: $\sup_t|\Fhat(t)-F_Q(t)|\to0$ in probability with no further condition, as shown in Lemma~4 of \PX{} by a Feller-type weak law (the condition $x\Prob_P(v>x)\to0$ used there is automatic, since $\E_Pv=1$), and \Pt{} gives finite-sample bounds through truncation of the ratio. This section identifies the rate and the limit law in the standard heavy-tailed model.

\begin{assumption}\label{ass:heavy}
\begin{enumerate}[label=(\roman*),leftmargin=*,itemsep=1pt]
\item (regular variation) For some $p\in(1,2)$ and a slowly varying $\ell$,
\[
\Prob_P(v>t)=t^{-p}\,\ell(t);
\]
let $a_m\to\infty$ satisfy $m\,\Prob_P(v>a_m)\to1$, so that $a_m=m^{1/p}\ell^\#(m)$ for a slowly varying $\ell^\#$.
\item (tail decoupling) There is a probability law $\Pi$ on $\R$ such that, for all $r\ge1$ and all $\Pi$-continuity sets $A$,
\[
\Prob_P\bigl(v>rt,\ S\in A\mid v>t\bigr)\longrightarrow r^{-p}\,\Pi(A)\qquad\text{as }t\to\infty .
\]
\end{enumerate}
\end{assumption}

Condition (i) says $\E_Pv^{p'}<\infty$ for $p'<p$ and $=\infty$ for $p'>p$; in particular $\E_Pv<\infty$ and $K^2=\infty$. Condition (ii) says that, in the far tail of the shift, the magnitude of the ratio and the score become independent, with $\Pi$ the law of the score of a point drawn from the extreme part of the shift; it holds whenever the conditional law of $S$ given $v(X)=u$ converges weakly to $\Pi$ as $u\to\infty$, which is the typical situation when $v$ is large on a region of covariate space where the conditional law of the score stabilizes (Section~\ref{sec:experiments}{} has an example). For a Borel set $A\subset\R$ write
\[
\mathbf c_A(s)=\ind\{s\in A\}-Q(S\in A)
\]
for the centered indicator of $A$.

\begin{theorem}[stable limits for the weighted score process]\label{thm:stablefidi}
Under Assumption~\ref{ass:heavy}, for all $t_1<\dots<t_k$ with $\Pi(\{t_j\})=0$,
\[
\frac{m}{a_m}\Bigl(\Fhat(t_j)-F_Q(t_j)\Bigr)_{j=1}^k\weak\bigl(\Scal(t_j)\bigr)_{j=1}^k,
\]
where $\Scal$ is the mean-zero infinitely divisible random vector with no Gaussian part and L\'evy measure
\[
\nu(\cdot)=\int_0^\infty\!\!\int\ind\bigl\{r\,(\mathbf c_{(-\infty,t_j]}(s))_{j=1}^k\in\cdot\bigr\}\,\Pi(\dd s)\,p\,r^{-p-1}\dd r,
\]
that is,
\[
\E\,e^{i\langle u,\Scal\rangle}=\exp\Bigl\{\int\bigl(e^{i\langle u,y\rangle}-1-i\langle u,y\rangle\bigr)\,\nu(\dd y)\Bigr\}.
\]
\end{theorem}

The rate is $m/a_m=m^{1-1/p}/\ell^\#(m)$, and the limit is a $p$-stable process in $t$ driven by the point process of the largest weights: in the LePage representation \citep{lepage1981convergence,samorodnitsky1994stable} the points $\Gamma_j^{-1/p}$, with $\Gamma_j$ the arrival times of a unit-rate Poisson process, are placed at i.i.d.\ scores $\eta_j\sim\Pi$, and $\Scal(t)$ is the compensated sum of $\Gamma_j^{-1/p}(\ind\{\eta_j\le t\}-F_Q(t))$ (for $p>1$ the series converges only after centering). The coverage at a fixed level inherits an explicit one-dimensional stable law.

\begin{theorem}[stable limit of the coverage error]\label{thm:stablecoverage}
Let Assumption~\ref{ass:heavy} hold, let $F_Q$ be continuous, fix $\beta\in(0,1)$, let $q=q_Q(\beta)$ and assume $\Pi(\{s:F_Q(s)=\beta\})=0$. Put
\[
\pi_\beta:=\Pi\bigl((-\infty,q]\bigr),
\]
the probability that a score from the far tail of the shift falls inside the prediction set. Then
\[
\frac{m}{a_m}\Bigl(P_e(\tauhat(\beta))-(1-\beta)\Bigr)\weak\Scal_\beta\sim S_p(\varsigma_\beta,\kappa_\beta,0),
\]
a $p$-stable law in the parametrization of \citet{samorodnitsky1994stable} with mean zero, skewness and scale
\[
\kappa_\beta=\frac{\pi_\beta(1-\beta)^p-(1-\pi_\beta)\beta^p}{\pi_\beta(1-\beta)^p+(1-\pi_\beta)\beta^p},\qquad
\varsigma_\beta^{\,p}=\Gamma(1-p)\cos\!\Bigl(\frac{\pi p}{2}\Bigr)\bigl\{\pi_\beta(1-\beta)^p+(1-\pi_\beta)\beta^p\bigr\}.
\]
\end{theorem}

Both proofs (Appendix~\ref{app:stable}) verify the classical conditions for convergence of row sums of a triangular array to an infinitely divisible law \citep{gnedenko1954limit,kallenberg2002foundations}: the L\'evy measure condition follows from (i)--(ii), the absence of a Gaussian part from Karamata's theorem, and the centering from $\E_P[v\,\mathbf c_A(S)]=0$. The coverage theorem uses the exact duality
\[
\bigl\{F_Q(\tauhat(\beta))\le w\bigr\}=\bigl\{\Fhat^{W}(w)\ge\beta\bigr\}
\]
for the weighted empirical distribution $\Fhat^W$ of the $Q$-uniform variables $W_i=F_Q(S_i)$---the quantile-crossing device of \Pt{}---at the moving threshold $w_m=\beta-za_m/m$; this is what absorbs the jumps of $\Fhat$, whose largest size $\max_i\hat w_i\asymp a_m/m$ is of the same order as the fluctuation and would defeat the Bahadur--Kiefer argument of Theorem~\ref{thm:coverage}.

\begin{remark}[interpretation and consequences]\label{rem:stable}
(a) The two sides of $\Scal_\beta$ come from the two things a single calibration point of enormous weight can do: if its score falls inside the set it pulls the threshold down and the set under-covers (positive jumps, intensity $\propto\pi_\beta(1-\beta)^p$); if it falls outside it pushes the threshold up (negative jumps, intensity $\propto(1-\pi_\beta)\beta^p$). The skewness $\kappa_\beta$ is the balance of the two; for small miscoverage $1-\beta$ and $\pi_\beta$ not close to one, $\kappa_\beta\approx-1$ and the heavy tail is on the side of over-coverage, while $\pi_\beta\to1$ (the extreme part of the shift lands inside the set) makes under-coverage heavy-tailed. $\pi_\beta$ plays, in this regime, the role that the coupling term $\Lambda(\alpha)$ plays in the variance proxy of \Pt{}.
(b) The inflation needed for an $(\alpha,\delta)$-PAC guarantee is of order $a_m/m=m^{-(1-1/p)}\ell^\#(m)$, not $m^{-1/2}$; a Gaussian plug-in based on an empirical variance is invalid, since the empirical variance of the weights diverges. The truncation bounds of \PX{} (Lemma~4) and \Pt{} (Remark~7), and the clipped ratios of \citet{wang2026weight}, are the finite-sample tools in this regime, and the theorem identifies the price they pay.
(c) At the boundary $p=2$ with $\E_Pv^2=\infty$ the limit is Gaussian again, with a normalization determined by the truncated second moment $\E_P[v^2\ind\{v\le x\}]$, which is slowly varying: for $\Prob_P(v>t)\sim ct^{-2}$ it is $\sqrt{m/\log m}$, for $\Prob_P(v>t)\sim c/(t^2\log t)$ it is $\sqrt{m/\log\log m}$. Since $\E_Pv=1$, the tail index cannot fall below $1$; at $p=1$ consistency persists by Feller's weak law (\PX{}, Lemma~4), but $m/a_m$ is then slowly varying, so no polynomial rate holds, while rates faster than any power of $\log m$ are possible for suitable slowly varying tails. The regime map is therefore: $\chi^2(Q\,\|\,P)<\infty$, Gaussian at $\sqrt m$ (Section~\ref{sec:fclt}); $1<p<2$, stable at $m^{1-1/p}$ up to the slowly varying factor $\ell^\#$; $p=1$, consistency without a polynomial rate. It is the conformal counterpart of the sample-size requirements for importance sampling in \citet{chatterjee2018sample} and of the Pareto-tail diagnostics of \citet{vehtari2024pareto}.
\end{remark}%

\section{Numerical illustration}\label{sec:experiments}

Unless stated otherwise, the experiments use the same regression model,
\[
Y=2X+\varepsilon,\qquad \varepsilon\mid X\sim N\bigl(0,(0.2+|X|)^2\bigr),
\]
the fixed predictor $\hat\mu(x)=2x$ (a learned predictor is used in the last experiment), and the score $S=|Y-\hat\mu(X)|=|\varepsilon|$, so that
\[
F_{S|X=x}(t)=2\,\Phi\bigl(t/(0.2+|x|)\bigr)-1
\]
depends on $x$ and the score is coupled with the covariate. Three covariate shifts are used throughout, listed in Table~\ref{tab:shifts}; two further shifts with the same $K^2$ as shift (A) appear in Table~\ref{tab:coupling}, and the estimated-weights experiment below uses{} a fourth, parametric shift (D) with its own noise scale. In shift (C) the ratio is large only near $x=0$, where the conditional law of the score is that of $0.2\,|N(0,1)|$, which is the law $\Pi$ of Assumption~\ref{ass:heavy}(ii). The population quantities $F_Q$, $G$, $\sigma^2(\beta)$ and $\gamma$ are computed by numerical quadrature, with $F_Q$ and $G$ tabulated on a grid of thresholds and interpolated linearly; the miscoverage $P_e(\tauhat)=1-F_Q(\tauhat)$ is thus computed numerically, without test-sample Monte Carlo error. Code and stored results are available at \codeurl.

\begin{table}[t]
\centering\small
\caption{The three covariate shifts used throughout the experiments.}
\label{tab:shifts}
\begin{tabular}{llllll}
\toprule
shift & $P_X$ & $Q_X$ & $v(x)$ & $B$ & $K^2$\\
\midrule
(A) & $U(0,1)$ & density $4x^3$ & $4x^3$ & $4$ & $16/7\approx2.29$\\
(B) & $N(0,1)$ & $N(1,1)$ & $e^{x-1/2}$ & $\infty$ & $e\approx2.72$\\
(C) & $U(0,1)$ & density $\tfrac13x^{-2/3}$ & $\tfrac13x^{-2/3}$ & $\infty$ & $\infty$ ($p=3/2$)\\
\bottomrule
\end{tabular}
\end{table}

\begin{figure}[t]
\centering
\includegraphics[width=0.95\linewidth]{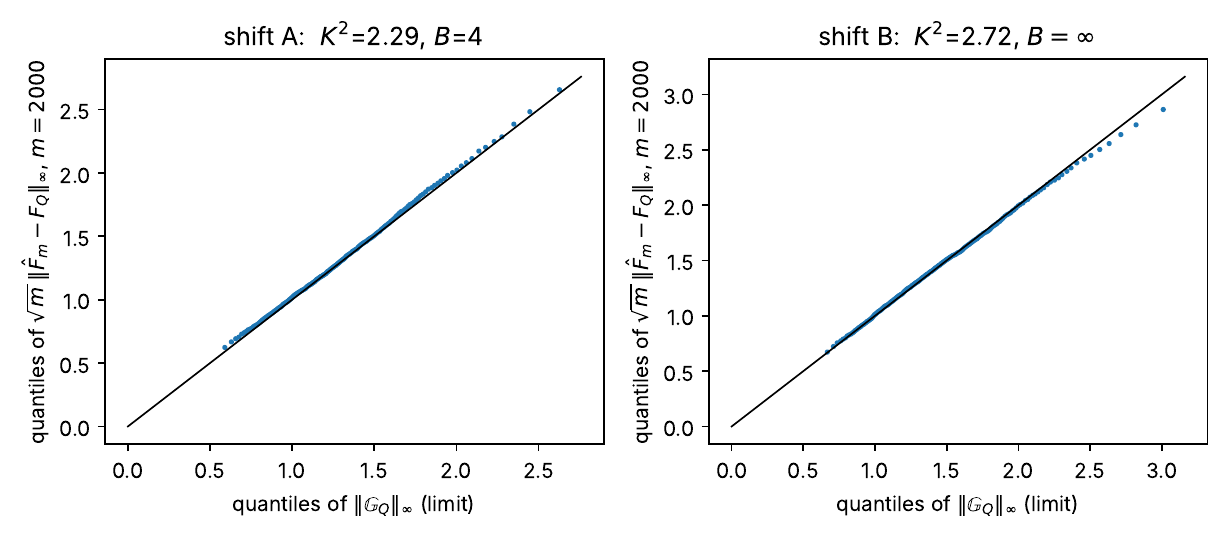}
\caption{Quantiles of $\sqrt m\|\Fhat-F_Q\|_\infty$ at $m=2000$ against quantiles of $\|\GQ\|_\infty$ simulated from \eqref{eq:brownian}, for shifts (A) and (B).}
\label{fig:ks}
\end{figure}

\paragraph{The functional limit (Theorem~\ref{thm:fclt}, Corollary~\ref{cor:asympdkw}).}
Table~\ref{tab:ks} compares quantiles of the weighted Kolmogorov--Smirnov statistic $\sqrt m\|\Fhat-F_Q\|_\infty$ over $3000$ calibration samples with those of $\|\GQ\|_\infty$, simulated from the representation \eqref{eq:brownian} with $20{,}000$ Brownian paths on a grid of $4001$ points. The agreement is within Monte Carlo error already at $m=500$, for the bounded shift (A) and for the unbounded one (B) alike; Figure~\ref{fig:ks} shows the full quantile--quantile comparison at $m=2000$. The $0.95$-quantiles of the limit, $2.03$ and $2.34$, are about $1.3$--$1.4$ times $K$ (the classical Kolmogorov--Smirnov $0.95$-quantile is $1.358$), illustrating the scale $K$ of Corollary~\ref{cor:asympdkw}. By comparison, the deviation term of the weighted DKW bound of \PX{} at $\delta=0.05$ is
\[
B\sqrt{2\log 80}\approx11.8
\]
for shift (A), with an unknown constant on top, and the bound is void for shift (B).

\begin{table}[t]
\centering\small
\caption{Quantiles of $\sqrt m\|\Fhat-F_Q\|_\infty$ (3000 replications) against those of the limit $\|\GQ\|_\infty$ (20{,}000 simulated paths).}
\label{tab:ks}
\begin{tabular}{llccccc}
\toprule
shift & & $0.50$ & $0.90$ & $0.95$ & $0.99$ & mean\\
\midrule
(A) & limit & 1.225 & 1.813 & 2.030 & 2.447 & 1.285\\
 & $m=500$ & 1.236 & 1.817 & 1.988 & 2.412 & 1.291\\
 & $m=2000$ & 1.229 & 1.850 & 2.054 & 2.483 & 1.301\\
\midrule
(B) & limit & 1.375 & 2.101 & 2.338 & 2.821 & 1.457\\
 & $m=500$ & 1.358 & 2.082 & 2.285 & 2.740 & 1.432\\
 & $m=2000$ & 1.386 & 2.090 & 2.306 & 2.727 & 1.456\\
\bottomrule
\end{tabular}
\end{table}

\paragraph{Simultaneous guarantees (Corollary~\ref{cor:simultaneous}).}
For shift (A), $m=2000$ and $\delta=0.1$, Table~\ref{tab:simul} reports the probability that \emph{at least one} of the levels in $\Bcal$ under-covers after inflation, for four procedures: no inflation; the single-level Gaussian plug-in $z_{0.9}\,\hat\sigma(\beta)/\sqrt m$ at each level, which is the asymptotically exact single-level oracle, and a Bonferroni-corrected version of it; and the joint multiplier-bootstrap inflation $\hat c/\sqrt m$ of Corollary~\ref{cor:simultaneous}(a) ($400$ Gaussian multiplier draws, unstudentized). For the finite set $\Bcal=\{0.80,0.85,0.90,0.95\}$ the single-level inflation fails jointly $22\%$ of the time, Bonferroni is conservative ($6\%$), and the joint bootstrap attains the nominal $\delta$ while inflating less than Bonferroni. For the $61$-point grid $\Bcal=\{0.80,0.8025,\dots,0.95\}$, used both for the bootstrap and for the evaluation (the experiment certifies the grid; Corollary~\ref{cor:simultaneous}(c) extends the guarantee to the interval as the mesh decreases), the single-level inflation fails for some level $31\%$ of the time; a Bonferroni correction over the $61$ points would use $z_{1-0.1/61}=2.94$ in place of $z_{0.9}=1.28$, an inflation of about $0.034$; the joint bootstrap again attains $\delta$ with an inflation of $0.024$, barely above that of the four-point set---the supremum of the coverage process over a fine grid is only slightly larger than over a few well-spread points, as the correlation structure of $\Lb$ predicts. Monte Carlo standard errors of the failure rates are about $0.008$.

\begin{table}[t]
\centering\small
\caption{Probability that some level in $\Bcal$ under-covers, and mean inflation, for shift (A), $m=2000$, $\delta=0.1$, $1500$ calibration samples. Pointwise: $z_{1-\delta}\hat\sigma(\beta)/\sqrt m$ at each level; Bonferroni: $z_{1-\delta/|\Bcal|}\hat\sigma(\beta)/\sqrt m$; joint: $\hat c/\sqrt m$ from Corollary~\ref{cor:simultaneous}.}
\label{tab:simul}
\begin{tabular}{lcccc}
\toprule
 & none & pointwise & Bonferroni & joint bootstrap\\
\midrule
$\Bcal=\{0.80,0.85,0.90,0.95\}$: joint failure & 0.741 & 0.223 & 0.061 & 0.110\\
\quad mean inflation & -- & 0.0149 & 0.0228 & 0.0208\\
$61$-point grid on $[0.80,0.95]$: joint failure & 0.817 & 0.315 & -- & 0.107\\
\quad mean inflation & -- & 0.0151 & -- & 0.0238\\
\bottomrule
\end{tabular}
\end{table}

\paragraph{Estimated weights (Example~\ref{ex:meanshift}).}
For shift (B), $\beta=0.8$ and $m=2000$, the tilt parameter is estimated by $\hat\theta=\bar X^{\rm tgt}-\bar X^{\rm src}$ from an unlabeled target sample of size $N_t$ and a source sample that is either independent of the calibration data (of size $N_s=m$) or \emph{is} the calibration covariates. Here $\sigma^2(\beta)=0.736$ and $\gamma(q)=-0.158$, so the predicted effects of weight estimation, $\gamma(q)^2(\lambda_t\pm1)$ and $\gamma(q)^2(\lambda_t+\lambda_s)$, lie between $-0.019$ and $+0.075$: small relative to the Monte Carlo error of a variance estimate ($\pm0.008$ with $20{,}000$ replications, the weights having heavy fourth moments). All variants are therefore computed from the \emph{same} calibration and auxiliary draws in each replication, so that the paired difference between the variance of a variant and that of the oracle is estimated with a standard error of $0.003$--$0.005$. Table~\ref{tab:est} reports both. The five paired differences match the predictions of Corollary~\ref{cor:estcoverage} to within one standard error: estimating the ratio from an independent source sample always costs variance; using the calibration covariates with $N_t=m$ costs nothing; and using them with $N_t=4m$ gives a variance \emph{below} the oracle's, by $0.020\pm0.003$ against the predicted $0.019$---eight standard errors from zero. A worst-case bound in terms of $\|v-\tilde v\|_{P,1}$ cannot distinguish these cases.

\begin{table}[t]
\centering\small
\caption{Variance of $\sqrt m\{P_e(\tauchk(0.8))-0.2\}$ with estimated tilt parameter, shift (B), $m=2000$, $20{,}000$ replications with common random numbers. ``cal'' means the calibration covariates are the source sample; ``indep'' an independent source sample of size $m$. The last two columns give the difference from the oracle variance predicted by Corollary~\ref{cor:estcoverage} and its paired empirical estimate (one standard error in parentheses).}
\label{tab:est}
\begin{tabular}{lcccc}
\toprule
weights & predicted $\tau^2$ & empirical & predicted $\tau^2-\sigma^2$ & empirical, paired\\
\midrule
oracle $v$ & 0.736 & 0.734 & -- & --\\
cal, $N_t=4m$ & 0.718 & 0.714 & $-0.019$ & $-0.020$ (0.003)\\
cal, $N_t=m$ & 0.736 & 0.732 & $0$ & $-0.002$ (0.003)\\
indep, $N_t=4m$ & 0.768 & 0.766 & $+0.031$ & $+0.032$ (0.003)\\
indep, $N_t=m$ & 0.787 & 0.785 & $+0.050$ & $+0.051$ (0.003)\\
cal, $N_t=m/4$ & 0.812 & 0.812 & $+0.075$ & $+0.078$ (0.005)\\
\bottomrule
\end{tabular}
\end{table}

\paragraph{Estimated weights in a bounded parametric model (Theorems~\ref{thm:parametric} and~\ref{thm:estboot}).}
This experiment isolates the two questions at the center of the paper in a bounded model in which every prediction of Theorem~\ref{thm:parametric} can be computed. Shift (D): $P_X=\mathrm{Beta}(\theta_s,1)$ and $Q_X=\mathrm{Beta}(\theta_t,1)$ with $\theta_s=3$, $\theta_t=4$, so $v(x)=(\theta_t/\theta_s)x^{\theta_t-\theta_s}=4x/3$, $B=4/3$, $K^2=16/15$, and the noise scale $2-1.8X$ in place of $0.2+|X|$, so that the predictor is noisy where the training distribution has its mass and accurate where the test distribution has its mass. Both shape parameters are estimated by maximum likelihood, $\hat\theta=-n/\sum_i\log X_i$, from a target sample of size $N_t$ and from a source sample that is either independent of the calibration data ($N_s=m$) or the calibration covariates; $\hat v(x)=(\hat\theta_t/\hat\theta_s)x^{\hat\theta_t-\hat\theta_s}$. With $\gamma:=\E_Q[\log X\{F_{S\mid X}(q)-\beta\}]$ one has $\gamma_t=-\gamma_s=\gamma$, $\Var\psi_s=\theta_s^2$ and $\Var\psi_t=\theta_t^2$, hence
\begin{align*}
\tau^2&=\sigma^2+\gamma^2(\theta_s^2\lambda_s+\theta_t^2\lambda_t)&&\text{(independent source sample)},\\
\tau^2&=\sigma^2+\gamma^2(\theta_t^2\lambda_t-\theta_s^2)&&\text{(calibration covariates as source sample)};
\end{align*}
reuse lowers the variance below the oracle's if and only if $N_t>(\theta_t/\theta_s)^2m$. At $\beta=0.8$, $\sigma^2=0.147$ and $\gamma=0.043$. Table~\ref{tab:estw} shows, for $m=500$ and $2000$, that the empirical variances follow the predictions in both directions (from $+69\%$ at $N_t=m/4$ to $-10\%$ at $N_t=16m$ with the calibration covariates as source), that the fixed-weight multiplier bootstrap estimates $\sigma^2$ in every configuration and fails at the rate $1-\Phi(z_{0.9}\sigma/\tau)$, between $0.08$ and $0.17$, and that the pairs bootstrap tracks $\tau^2$ and is calibrated within Monte Carlo error already at $m=500$.

\begin{table}[t]
\centering\footnotesize\setlength{\tabcolsep}{3.5pt}
\caption{Estimated weights in the bounded model of Section~\ref{sec:experiments}, $\beta=0.8$, $\delta=0.1$, $2000$ replications per row: predicted variance $\tau^2(\beta)$ of Theorem~\ref{thm:parametric} ($\tau^2=\sigma^2=0.147$ for the true ratio), empirical variance of $\sqrt m\{\beta-F_Q(\tauchk(\beta))\}$ (standard errors at most $0.008$), average variance estimates of the fixed-weight multiplier bootstrap (naive) and of the pairs bootstrap ($R=400$), limiting failure probability $1-\Phi(z_{0.9}\sigma/\tau)$ of the naive inflation (the pairs bootstrap has limit $0.100$), and probability that the set run at the inflated level under-covers (standard errors at most $0.008$).}
\label{tab:estw}
\begin{tabular}{llcccccccc}
\toprule
 & & & & & \multicolumn{2}{c}{bootstrap variance} & naive & \multicolumn{2}{c}{failure}\\
$m$ & source sample & $N_t$ & $\tau^2$ & empirical & naive & pairs & limit & naive & pairs\\
\midrule
$500$ & true ratio & -- & 0.147 & 0.147 & 0.147 & -- & 0.100 & 0.095 & --\\
 & independent, $N_s=m$ & $m$ & 0.193 & 0.194 & 0.147 & 0.194 & 0.132 & 0.125 & 0.097\\
 & calibration covariates & $m/4$ & 0.248 & 0.245 & 0.148 & 0.245 & 0.162 & 0.155 & 0.101\\
 & calibration covariates & $m$ & 0.160 & 0.157 & 0.147 & 0.160 & 0.110 & 0.101 & 0.090\\
 & calibration covariates & $4m$ & 0.138 & 0.136 & 0.147 & 0.137 & 0.093 & 0.092 & 0.106\\
 & calibration covariates & $16m$ & 0.133 & 0.127 & 0.147 & 0.132 & 0.088 & 0.081 & 0.098\\
\addlinespace[2pt]
$2000$ & true ratio & -- & 0.147 & 0.146 & 0.147 & -- & 0.100 & 0.089 & --\\
 & independent, $N_s=m$ & $m$ & 0.193 & 0.194 & 0.147 & 0.193 & 0.132 & 0.137 & 0.104\\
 & calibration covariates & $m/4$ & 0.248 & 0.260 & 0.147 & 0.248 & 0.162 & 0.168 & 0.102\\
 & calibration covariates & $m$ & 0.160 & 0.158 & 0.147 & 0.160 & 0.110 & 0.112 & 0.103\\
 & calibration covariates & $4m$ & 0.138 & 0.137 & 0.147 & 0.138 & 0.093 & 0.089 & 0.103\\
 & calibration covariates & $16m$ & 0.133 & 0.131 & 0.147 & 0.132 & 0.088 & 0.086 & 0.101\\
\bottomrule
\end{tabular}
\end{table}

\paragraph{Estimated weights, feasible bootstrap (Theorem~\ref{thm:estboot}).}
Table~\ref{tab:estboot} compares, for shift (B), $m=2000$ and $\delta=0.1$, the multiplier bootstrap of Section~\ref{sec:bootstrap} applied with the estimated weights held fixed and the pairs bootstrap of Section~\ref{sec:estboot}, both with $R=400$ draws, in four configurations of Example~\ref{ex:meanshift} in which $\tau^2(\beta)\ne\sigma^2(\beta)=0.736$: an independent source sample with $N_s=N_t=m$ ($\tau^2=0.787$), the calibration covariates as source sample with $N_t=4m$ ($\tau^2=0.718$), with $N_t=m/4$ ($\tau^2=0.812$) and with $N_t=m/16$ ($\tau^2=1.113$), the last added because in the first three the weight-estimation term is only a few percent of $\sigma^2$. The naive bootstrap estimates $\sigma^2(\beta)$ in every configuration (average estimate $0.72$--$0.75$), while the pairs bootstrap tracks $\tau^2(\beta)$: $0.773$, $0.700$, $0.801$ and $1.126$ against the predicted $0.787$, $0.718$, $0.812$ and $1.113$. The table also gives the limiting failure probability of the naive single-level inflation, $1-\Phi(z_{0.9}\sigma/\tau)$, which is $0.107$, $0.097$, $0.111$ and $0.149$ against $0.100$ for the pairs bootstrap: the asymptotic difference between the two procedures is small in the first three configurations and visible in the fourth, where the naive inflation fails $17\%$ of the time against $13\%$ for the pairs bootstrap (joint failure over the four levels $0.170$ against $0.135$). At $m=8000$ ($1500$ replications, first three configurations) the bootstrap variance estimates are $0.728$--$0.731$ (naive) against $0.778$, $0.711$ and $0.802$ (pairs), and the single-level failure rates are $0.131$, $0.119$, $0.139$ (naive) against $0.123$, $0.118$, $0.125$ (pairs).

Both bootstraps over-reject by a few points at these sample sizes, including the pairs bootstrap whose variance estimate is correct on average. Table~\ref{tab:estbootdecomp} locates the excess. By the quantile duality of Section~\ref{sec:rates}, the set run at the level $\beta+c/\sqrt m$ under-covers if and only if $\sqrt m\{\Fchk(q-)-\beta\}\ge c$, $q=q_Q(\beta)$, so the failure rate can be examined for the \emph{fixed} statistic $\sqrt m\{\Fchk(q-)-\beta\}$ against different critical values. With the population critical value $z_{0.9}\tau(\beta)$ the failure rate is $0.101$ and $0.095$ in the two configurations shown: the statistic is well approximated by its $N(0,\tau^2)$ limit in the relevant tail, so the excess is not in the statistic. With $z_{0.9}\sigma(\beta)$ it is $0.112$ and $0.147$, the asymptotic benchmarks of the naive procedure. With the bootstrap critical values it rises to $0.13$--$0.14$ (pairs) and $0.13$--$0.18$ (naive), and increasing the number of draws from $400$ to $4000$, or replacing the empirical quantile by the rank-corrected order statistic of Procedure~1, changes none of these rates by more than Monte Carlo error. What differs between the fixed and the bootstrap critical values is that the latter are random, with a coefficient of variation of $12\%$--$17\%$, and \emph{negatively} correlated with the statistic (correlation $-0.5$ to $-0.75$): with lognormal weights the bootstrap variance is driven by the few largest weights, and a sample in which those weights fall below $q$ has both a large coverage error and a small estimated variance. The excess is thus a second-order effect of studentizing by a dependent, heavy-tailed variance estimate. It is common to both bootstraps, it lies outside the bounded-weight rates of Section~\ref{sec:rates}, and it is separate from the first-order comparison between the two bootstraps, which is what Theorem~\ref{thm:estboot} is about and which the tables confirm, with a size governed by $\tau^2-\sigma^2$.

\begin{table}[t]
\centering\footnotesize\setlength{\tabcolsep}{4pt}
\caption{Estimated weights, shift (B), $\beta=0.8$, $m=2000$, $\delta=0.1$, $3000$ replications: predicted and empirical variance of $\sqrt m\{\beta-F_Q(\tauchk(\beta))\}$, the average bootstrap variance estimates of the naive (weights fixed) and pairs bootstraps ($R=400$), the limiting failure probability $1-\Phi(z_{0.9}\sigma/\tau)$ of the naive single-level inflation (the pairs bootstrap has limit $0.100$), and the probability that the inflated set under-covers, at the single level and jointly over the levels $\{0.80,0.85,0.90,0.95\}$. Standard errors of the failure rates are about $0.006$.}
\label{tab:estboot}
\begin{tabular}{lccccccccc}
\toprule
 & & & \multicolumn{2}{c}{bootstrap variance} & naive & \multicolumn{2}{c}{failure, single level} & \multicolumn{2}{c}{failure, joint}\\
weights & $\tau^2(0.8)$ & empirical & naive & pairs & limit & naive & pairs & naive & pairs\\
\midrule
indep, $N_s=N_t=m$ & 0.787 & 0.781 & 0.724 & 0.773 & 0.107 & 0.137 & 0.128 & 0.140 & 0.139\\
cal, $N_t=4m$ & 0.718 & 0.670 & 0.722 & 0.700 & 0.097 & 0.111 & 0.114 & 0.115 & 0.127\\
cal, $N_t=m/4$ & 0.812 & 0.788 & 0.725 & 0.801 & 0.111 & 0.124 & 0.114 & 0.135 & 0.130\\
cal, $N_t=m/16$ & 1.113 & 1.117 & 0.749 & 1.126 & 0.149 & 0.169 & 0.129 & 0.170 & 0.135\\
\bottomrule
\end{tabular}
\end{table}

\begin{table}[t]
\centering\footnotesize\setlength{\tabcolsep}{4pt}
\caption{Where the over-rejection comes from: shift (B), $\beta=0.8$, $m=2000$, $\delta=0.1$, $1500$ replications. Failure rate of the set run at $\beta+c/\sqrt m$ for the population critical values $c=z_{0.9}\tau$ and $z_{0.9}\sigma$ and for bootstrap critical values with $R$ draws (empirical $(1-\delta)$-quantile); mean and standard deviation of $c$ and its correlation with the statistic $\sqrt m\{\Fchk(q-)-\beta\}$. Standard errors of the failure rates are about $0.008$.}
\label{tab:estbootdecomp}
\begin{tabular}{lcccccccc}
\toprule
 & \multicolumn{4}{c}{indep, $N_s=N_t=m$ ($\tau^2=0.787$)} & \multicolumn{4}{c}{cal, $N_t=m/16$ ($\tau^2=1.113$)}\\
critical value $c$ & failure & mean $c$ & sd $c$ & corr & failure & mean $c$ & sd $c$ & corr\\
\midrule
$z_{0.9}\tau$ & 0.101 & 1.137 & -- & -- & 0.095 & 1.352 & -- & --\\
$z_{0.9}\sigma$ & 0.112 & 1.100 & -- & -- & 0.147 & 1.100 & -- & --\\
naive, $R=400$ & 0.138 & 1.085 & 0.152 & $-0.48$ & 0.175 & 1.091 & 0.182 & $-0.64$\\
pairs, $R=400$ & 0.129 & 1.116 & 0.145 & $-0.49$ & 0.139 & 1.281 & 0.187 & $-0.70$\\
naive, $R=4000$ & 0.134 & 1.092 & 0.132 & $-0.54$ & 0.179 & 1.095 & 0.166 & $-0.70$\\
pairs, $R=4000$ & 0.129 & 1.122 & 0.130 & $-0.57$ & 0.141 & 1.294 & 0.179 & $-0.76$\\
\bottomrule
\end{tabular}
\end{table}

\paragraph{Same $K^2$, different coupling (Theorem~\ref{thm:coverage}, Section~\ref{sec:efficiency}).}
The effective sample size $m/K^2$ of \citet{kish1965survey} is often used as a summary of the cost of a shift. Table~\ref{tab:coupling} shows that it is not the whole story. Three shifts with $P_X$ uniform and the same $K^2=16/7$ are compared: the shift (A), $v(x)=4x^3$, which puts the test mass where the noise, hence the score, is large; its mirror image $v(x)=4(1-x)^3$, which puts it where the score is small; and a shift through an independent second covariate: $X=(X_1,X_2)$ with independent $U(0,1)$ coordinates, the score built from $X_1$ as before ($Y=2X_1+\varepsilon$, $\varepsilon\mid X\sim N(0,(0.2+X_1)^2)$) and $v(x)=4x_2^3$, so that $Q_X=U(0,1)\otimes\mathrm{Beta}(4,1)$, $X_2\perp S$ and $\tilde G(\beta)=\beta K^2$ exactly. The Kish approximation $K^2\beta(1-\beta)$ is exact in the third case, while in the other two the actual variance at $\beta=0.8$ is $12\%$ above and $27\%$ below it, in the directions predicted by the sign of $\beta K^2-\tilde G(\beta)$: weight on large scores makes the coverage \emph{more} variable than the effective sample size suggests, weight on small scores less. The same quantity determines the gain available to the efficient threshold, which is zero in the third case and largest in the second.

\begin{table}[t]
\centering\small
\caption{Three shifts with the same $K^2=16/7$ (Kish effective sample size $7m/16$) and different coupling between the weights and the scores: the clock $\tilde G(\beta)$, the Kish approximation $K^2\beta(1-\beta)$, the predicted variance $\sigma^2(\beta)$ of Theorem~\ref{thm:coverage}, its empirical value over $20{,}000$ calibration samples of size $m=2000$ (standard errors $0.001$--$0.004$), and the efficiency gain $(\beta K^2-\tilde G(\beta))^2/(K^2-1)$ of Theorem~\ref{thm:effthreshold}.}
\label{tab:coupling}
\begin{tabular}{llcccccc}
\toprule
shift & $\beta$ & $\tilde G(\beta)$ & $\beta K^2$ & Kish & $\sigma^2(\beta)$ & empirical & gain\\
\midrule
$v(x)=4x^3$ (weight on large scores) & 0.8 & 1.755 & 1.829 & 0.366 & 0.410 & 0.409 & 0.004\\
 & 0.9 & 2.007 & 2.057 & 0.206 & 0.246 & 0.245 & 0.002\\
$v(x)=4(1-x)^3$ (weight on small scores) & 0.8 & 1.993 & 1.829 & 0.366 & 0.267 & 0.265 & 0.021\\
 & 0.9 & 2.171 & 2.057 & 0.206 & 0.114 & 0.114 & 0.010\\
$v=4X_2^3$, $X_2\perp S$ & 0.8 & 1.829 & 1.829 & 0.366 & 0.366 & 0.364 & 0.000\\
 & 0.9 & 2.057 & 2.057 & 0.206 & 0.206 & 0.205 & 0.000\\
\bottomrule
\end{tabular}
\end{table}

\paragraph{A learned predictor.}
The experiments above use the fixed predictor $\mu(x)=2x$. To show the theory in the setting it is meant for, a cubic least-squares fit $\hat\mu$ was trained once on an independent training split of $500$ points from $P$ with $Y=2X+0.5\sin(6X)+\varepsilon$, and then held fixed; its root-mean-square distance to the true regression function, $\{\E_Q(\mu(X)-\hat\mu(X))^2\}^{1/2}$, is $0.067$. All quantities of Theorem~\ref{thm:coverage} refer to the deployed model: $F_{Q,\hat\mu}$ and $G_{\hat\mu}$ are those of the score $|Y-\hat\mu(X)|$, computed here by quadrature given $\hat\mu$. For shift (A) and $m=1000$, the empirical variances of $\sqrt m\{\beta-F_{Q,\hat\mu}(\tauhat(\beta))\}$ over $5000$ calibration samples are $0.410, 0.338, 0.249, 0.140$ at $\beta=0.80,0.85,0.90,0.95$, against the predicted $0.410, 0.336, 0.246, 0.136$, and the joint multiplier-bootstrap inflation over these four levels (unstudentized, Corollary~\ref{cor:simultaneous}(a), $400$ draws) under-covers some level with probability $0.100$ at $\delta=0.1$ ($2000$ replications), against $0.219$ for the single-level inflation applied at each level and $0.737$ without inflation.

\begin{figure}[t]
\centering
\includegraphics[width=0.6\linewidth]{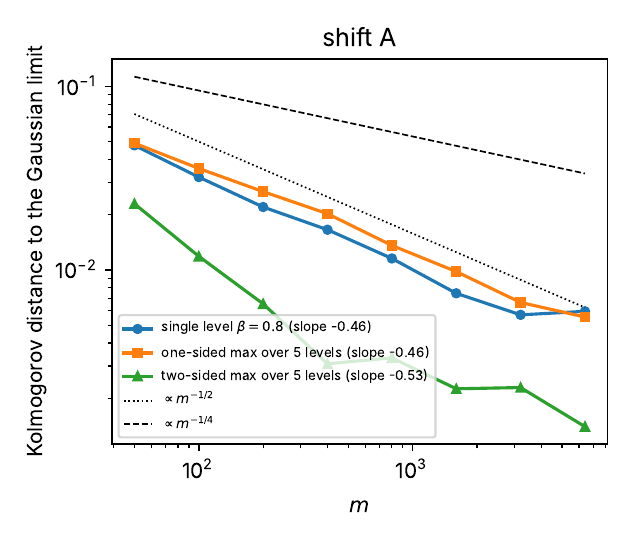}
\caption{Shift (A). Kolmogorov distance between the law of the studentized coverage error at a single level, and of its one- and two-sided maxima over five levels, and the corresponding Gaussian limits, against $m$ (log-log), with the reference slopes $-1/2$ and $-1/4$.}
\label{fig:rates}
\end{figure}

\paragraph{Rates (Proposition~\ref{prop:berryesseen}, Theorem~\ref{thm:rates}).}
For shift (A) and the levels $0.5,0.6,0.7,0.8,0.9$, Figure~\ref{fig:rates} shows the Kolmogorov distance between the Monte Carlo law of three statistics and their Gaussian limits as $m$ grows from $50$ to $6400$ ($200{,}000$ calibration samples per $m$; the Monte Carlo resolution is about $0.003$): the single-level statistic $\sqrt m\,\Cm(0.8)/\sigma(0.8)$, and the one- and two-sided maxima over the five levels, studentized by the true $\sigma(b_j)$. The single-level distance decays at the parametric rate of Proposition~\ref{prop:berryesseen}. The distance for the grid maximum decays at essentially the same rate, visibly faster than the $m^{-1/4}$ of Theorem~\ref{thm:rates}; the Bahadur--Kiefer step of the proof is a device, and the open problem (iv) of Section~\ref{sec:discussion} records the gap. The multiplier bootstrap critical value of the one-sided studentized maximum, at $\delta=0.1$ and with $400$ multiplier draws, rejects with frequency 0.163, 0.127 and 0.114 for $m=100$, $400$ and $1600$ ($10{,}000$ replications, standard error $0.003$): the bootstrap over-rejects at small $m$, where the studentization by $\hat\sigma$ and the skewness of the one-sided maximum both work against it, and converges to the nominal level slowly, in agreement with the modest rate of \eqref{eq:rate-boot}.

\paragraph{Efficiency (Theorems~\ref{thm:eff}--\ref{thm:effthreshold}).}
Table~\ref{tab:eff} reports, for shifts (A) and (B) and $\beta\in\{0.8,0.9\}$, the asymptotic variance $\sigma^2(\beta)$ of the oracle weighted threshold, the bound $\sigma^2_{\rm eff}(\beta)$ of the known-ratio model attained by the control-variate threshold $\tilde\tau$, the known-marginals bound $V_2$, the unknown-ratio bound $V_3$ for $m/N\in\{1,1/4\}$, and the empirical variances of $\sqrt m\{\beta-F_Q(\hat\tau)\}$ and $\sqrt m\{\beta-F_Q(\tilde\tau)\}$ over $20{,}000$ calibration samples of size $m=2000$ computed from the same data. Under the bounded shift (A) the ratio and the coverage indicator are nearly uncorrelated under $Q$ at these levels, $\tilde G(\beta)\approx\beta K^2$, and the four quantities are within two percent of each other: the oracle weighted threshold is practically efficient and nothing is to be gained. Under the tilt (B) the picture changes: at $\beta=0.8$ the oracle variance is $0.736$ (empirically $0.721$ here and $0.734$ in Table~\ref{tab:est}, two independent runs with heavy-tailed weights), the control-variate threshold attains $0.589$ (empirically $0.579$, a paired reduction of $0.141\pm0.01$ against the predicted $0.147$), and the known-marginals bound is $0.540$. The estimated-weight procedures of Table~\ref{tab:est}, which exploit the unlabeled target sample through a correctly specified parametric ratio model, reduce the variance by at most $0.019$ in the same setting; in this example most of the available improvement over the oracle lies in the constraint $\E_Pv=1$ that the self-normalization uses with the wrong coefficient, not in the covariate sample.

\begin{table}[t]
\centering\small
\caption{Asymptotic variances and efficiency bounds for the coverage error (Section~\ref{sec:efficiency}), with empirical variances at $m=2000$ over $20{,}000$ replications (Monte Carlo standard errors $0.003$--$0.007$). $\sigma^2$: oracle weighted threshold (Theorem~\ref{thm:coverage}); $\sigma^2_{\rm eff}$: known-ratio bound, attained by $\tilde\tau$ (Theorem~\ref{thm:effthreshold}); $V_2$: known marginals; $V_3$: unknown ratio with $N$ unlabeled target covariates, reported for $m/N\in\{1,1/4\}$, i.e.\ one and four target covariates per calibration point (Proposition~\ref{prop:bounds23}).}
\label{tab:eff}
\begin{tabular}{llccccccc}
\toprule
shift & $\beta$ & $\sigma^2$ & emp.\ $\hat\tau$ & $\sigma^2_{\rm eff}$ & emp.\ $\tilde\tau$ & $V_2$ & $V_3$, $m/N=1$ & $V_3$, $m/N=1/4$\\
\midrule
(A) & 0.8 & 0.410 & 0.400 & 0.406 & 0.397 & 0.402 & 0.407 & 0.404\\
(A) & 0.9 & 0.246 & 0.245 & 0.244 & 0.243 & 0.242 & 0.244 & 0.242\\
(B) & 0.8 & 0.736 & 0.721 & 0.589 & 0.579 & 0.540 & 0.574 & 0.548\\
(B) & 0.9 & 0.545 & 0.524 & 0.463 & 0.450 & 0.421 & 0.436 & 0.425\\
\bottomrule
\end{tabular}
\end{table}

\begin{figure}[t]
\centering
\includegraphics[width=0.95\linewidth]{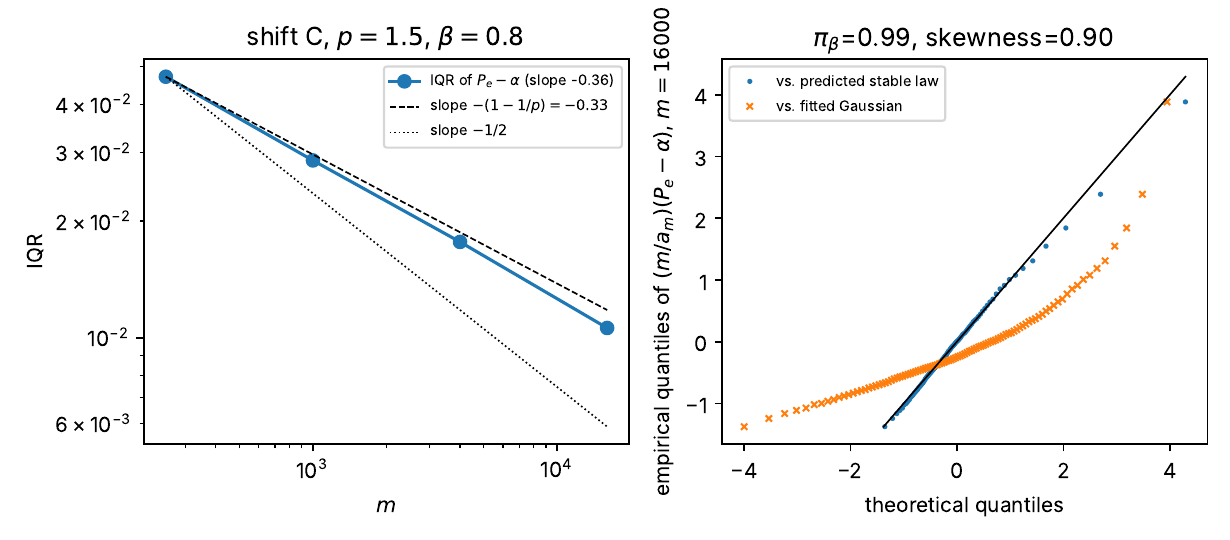}
\caption{Shift (C), $p=3/2$, $\beta=0.8$. Left: interquartile range of the coverage error against $m$ (log-log), with the stable and Gaussian reference slopes. Right: quantiles of $(m/a_m)\{P_e(\tauhat(\beta))-0.2\}$ at $m=16{,}000$ against the predicted stable law of Theorem~\ref{thm:stablecoverage} and against a fitted Gaussian.}
\label{fig:stable}
\end{figure}

\paragraph{The stable regime (Theorems~\ref{thm:stablefidi}--\ref{thm:stablecoverage}).}
For shift (C) and $\beta=0.8$, the left panel of Figure~\ref{fig:stable} shows the interquartile range of the coverage error $P_e(\tauhat(\beta))-0.2$ over $3000$ replications, for $m$ ranging from $250$ to $16{,}000$: the log-log slope is $-0.36$, against $-(1-1/p)=-0.33$ predicted by Theorem~\ref{thm:stablecoverage} and $-0.5$ for the Gaussian regime. Here $q=0.544$ and $\pi_\beta=0.993$: a score from the far tail of the shift almost surely falls inside the set, so the predicted skewness is $\kappa_\beta=0.90$ and the heavy tail is on the side of under-coverage. The right panel compares the empirical quantiles of
\[
\frac{m}{a_m}\bigl\{P_e(\tauhat(\beta))-0.2\bigr\},\qquad m=16{,}000,\quad a_m=\tfrac13\,m^{2/3},
\]
with those of the predicted law $S_{3/2}(0.380,0.90,0)$ and with a Gaussian fit: the stable law is matched across the bulk and the left tail (empirical $0.05$, $0.25$, $0.50$ and $0.75$ quantiles $-1.08$, $-0.60$, $-0.26$, $0.20$ against $-1.03$, $-0.59$, $-0.25$, $0.20$), with the heavy upper tail somewhat lighter at this $m$ than in the limit (empirical $0.95$-quantile $1.31$ against $1.42$), while the Gaussian fit misses both tails and the median.

\section{Discussion}\label{sec:discussion}

\paragraph{Three resolutions of one random variable.}
\PX{}, \Pt{}, \citet{gazin2024asymptotics} and this paper study the same object, the training-conditional coverage of the weighted quantile, at different resolutions (Table~\ref{tab:compare}). The finite-sample bounds of the first two papers are certificates: uniform in the threshold with an unspecified constant, with sub-Gaussian tails at the scale $B=\sup v$ or polynomial tails at the scale $\|v\|_{P,2}$, or at a fixed level with explicit constants at the scale of the variance proxy. This paper says what those bounds are converging to, and uses the limit for what certificates cannot do: joint inference across levels, resampling, estimated weights, efficiency, and the regime in which the $\sqrt m$ rate fails. Apart from Section~\ref{sec:rates}, whose constants are not explicit, nothing here is a certificate at finite $m$.

\paragraph{What the limit says about the finite-sample bounds.}
Two messages are quantitative. First, the variance proxy $\sigma_\alpha^2$ of \Pt{} is the value of $C(q,q)$ at the adversarial coupling, whereas the actual asymptotic variance $\sigma^2(\beta)$ is plug-in estimable: $\hat\sigma_m^2(\beta)$ of \eqref{eq:sigmahat} is consistent uniformly in $\beta$ under $K^2<\infty$, by the Glivenko--Cantelli theorem for the class $\{v^2(\ind\{S\le t\}-c)^2\}$ with integrable envelope $v^2$, and Lemma~\ref{lem:sigmahat} gives its rate for bounded $v$. Whether such an estimate can enter a finite-sample certificate without paying $B$ again in its concentration is the open question of \Pt{}; asymptotically it costs nothing, and Corollary~\ref{cor:finitesimul} gives a finite-sample answer up to an unspecified constant. Second, the bootstrap of Theorem~\ref{thm:bootstrap} replaces the distribution-free constants of concentration inequalities by the exact ones of the limit law, uniformly over levels, at the price of asymptotic validity; Theorem~\ref{thm:rates} prices this at $m^{-1/4}$ up to logarithms, and Proposition~\ref{prop:berryesseen} at the parametric rate $m^{-1/2}$ for a single level.

\paragraph{Open problems.}
(i) \emph{The stable regime at the process level.} The finite-dimensional convergence of Theorem~\ref{thm:stablefidi} should extend to convergence in a Skorokhod topology of the truncated process to the ``stable bridge'' $M((-\infty,t])-F_Q(t)M(\R)$, with $M$ the $p$-stable random measure with control measure $\Pi$, by the point-process methods of \citet{resnick1986point} together with a maximal inequality for the small-jump part; it would give the joint law across levels in that regime. (ii) \emph{Efficiency with an estimated ratio.} Theorem~\ref{thm:eff} settles the known-ratio model. In the two-sample model $\mathcal M_3$ the bound $V_3$ is attained by doubly robust procedures that estimate the conditional score distribution \citep{yang2024doubly} and, at a fixed $t$, by weighting with a nonparametrically estimated ratio, under smoothness and boundedness conditions on the ratio \citep[Theorem~7]{chen2008semiparametric}. Two questions remain. The uniform-in-level version of the latter under the second-moment conditions of this paper, that is, the process limit of $\Fchk$ when the ratio is estimated by a sieve, is open; Theorem~\ref{thm:estimated} gives the expansion once the sieve estimate is in a Donsker-type class, but the efficiency of the resulting process has not been established. And in a correctly specified parametric ratio model the bound lies below $V_3$; Example~\ref{ex:meanshift} shows that the plug-in estimated-weight procedure can fall below the oracle variance, but which estimators of the ratio attain the parametric bound is not known. (iii) \emph{A distribution-free diagnostic.} The diagnostic of Proposition~\ref{prop:twosample}{} has a limit law that depends on $(F_Q,G)$. Even the unnormalized process $\Gm(v\ind\{S\le\cdot\})$ has the bridge-type covariance $G(s\wedge t)-F_Q(s)F_Q(t)$, so no time change alone makes it distribution-free; a transformation in the spirit of Khmaladze's innovation martingale, applied to the mismatched-clock bridge, might remove the dependence on $(F_Q,G)$ and the need for the bootstrap there, and we leave its construction open. (iv) \emph{Sharp rates.} Proposition~\ref{prop:berryesseen} gives the parametric rate at a single level and Theorem~\ref{thm:rates} the rate $m^{-1/4}$ for the grid maximum; the simulations of Section~\ref{sec:experiments}{} suggest that the truth is $m^{-1/2}$ for the grid as well, and a proof would have to avoid the Bahadur--Kiefer step, perhaps through the duality of Proposition~\ref{prop:berryesseen} combined with a Gaussian comparison for covariances that differ by $O(m^{-1/2})$. An explicit-constant version, which would make Corollary~\ref{cor:finitesimul} a certificate in the sense of \Pt{}, is open.

\section*{Reproducibility}
Python code reproducing every table and figure (\texttt{code/experiments.py}) and the stored results (\texttt{results/}) are available at \codeurl.

\section*{Acknowledgments}
The author thanks Javad Alipanah for providing access to the AI tools used in this work.

\paragraph{Statement on AI use.}
AI tools were used in this work, under the direction of the author. Claude Fable 5.1 and 5.2 (Anthropic) were used to assist with writing, the development of the proofs, and running the simulations. GPT-6 Astra (OpenAI) was used to review the manuscript, and the errors it identified were corrected. The author takes full responsibility for the mathematical content and conclusions of the paper.

\bibliographystyle{plainnat}
\bibliography{refs}

\appendix
\section{Proofs for Section~\ref{sec:fclt}}\label{app:fclt}

We use the notation of \citet{vandervaart1996weak}: for a class $\mathcal F$ of square-integrable functions, the bracketing number $N_{[\,]}(\epsilon,\mathcal F,L_2(P))$ is the minimal number of brackets $[l,u]$ with
\[
\|u-l\|_{P,2}\le\epsilon
\]
needed to cover $\mathcal F$, and the bracketing entropy integral is
\[
J_{[\,]}(\eta,\mathcal F,L_2(P))=\int_0^\eta\sqrt{\log N_{[\,]}(\epsilon,\mathcal F,L_2(P))}\,\dd\epsilon .
\]

\begin{lemma}[partition lemma]\label{lem:partition}
Let $\mu$ be a finite Borel measure on $\R$ with total mass $M$ and let $\eta>0$. There are points
\[
-\infty=t_0<t_1<\dots<t_N=+\infty,\qquad N\le M/\eta+1,
\]
such that $\mu((t_i,t_{i+1}))\le\eta$ for all $i$.
\end{lemma}

\begin{proof}
Put $t_0=-\infty$ and, recursively,
\[
t_{i+1}=\sup\bigl\{t\in\R:\ \mu((t_i,t])\le\eta\bigr\},
\]
with $t_{i+1}=+\infty$ if $\mu((t_i,\infty))\le\eta$. Since $\mu((t_i,t])$ decreases to zero as $t$ decreases to $t_i$, the new point is strictly larger than the previous one. By monotone convergence and by right-continuity of $t\mapsto\mu((t_i,t])$,
\[
\mu\bigl((t_i,t_{i+1})\bigr)=\lim_{t\uparrow t_{i+1}}\mu((t_i,t])\le\eta,
\qquad\text{and}\qquad
\mu\bigl((t_i,t_{i+1}]\bigr)\ge\eta\quad\text{whenever }t_{i+1}<\infty,
\]
the latter because $\mu((t_i,t])>\eta$ for every $t>t_{i+1}$. The intervals $(t_i,t_{i+1}]$ with a finite right endpoint are disjoint and each has mass at least $\eta$, so there are at most $M/\eta$ finite points.
\end{proof}

\begin{lemma}[bracketing numbers of $\Hcal$]\label{lem:bracketH}
If $K^2<\infty$ then, for every $\epsilon>0$,
\[
N_{[\,]}\bigl(\epsilon,\Hcal,L_2(P)\bigr)\le\frac{8K^2}{\epsilon^2}+\frac{4K}{\epsilon}+3 .
\]
In particular the bound on $N_{[\,]}(\epsilon K,\Hcal,L_2(P))$ does not depend on $K$, and $J_{[\,]}(1,\Hcal,L_2(P))<\infty$.
\end{lemma}

\begin{proof}
Fix $\epsilon>0$. Apply Lemma~\ref{lem:partition} once to the measure
\[
\mu_1(A)=\E_P\bigl[v^2\ind\{S\in A\}\bigr]\qquad\text{(mass $K^2$)}\qquad\text{with}\qquad\eta=\epsilon^2/4,
\]
obtaining at most $4K^2/\epsilon^2+1$ intervals, and once to the measure
\[
\mu_2(A)=Q(S\in A)\qquad\text{(mass $1$)}\qquad\text{with}\qquad\eta=\epsilon/(2K),
\]
obtaining at most $2K/\epsilon+1$ intervals. Merging the two sets of partition points gives points
\[
-\infty=t_0<\dots<t_N=+\infty,\qquad N\le\frac{4K^2}{\epsilon^2}+\frac{2K}{\epsilon}+1,
\]
and every open interval $(t_i,t_{i+1})$ is contained in an interval of each of the two partitions, so that
\[
\E_P\bigl[v^2\ind\{t_i<S<t_{i+1}\}\bigr]\le\frac{\epsilon^2}{4},\qquad K\,Q(t_i<S<t_{i+1})\le\frac{\epsilon}{2}.
\]
For $t\in[t_i,t_{i+1})$ we have
\[
\ind\{s\le t_i\}\le\ind\{s\le t\}\le\ind\{s<t_{i+1}\},\qquad F_Q(t_i)\le F_Q(t)\le F_Q(t_{i+1}-),
\]
hence
\[
l_i:=v\bigl(\ind\{s\le t_i\}-F_Q(t_{i+1}-)\bigr)\ \le\ h_t\ \le\ v\bigl(\ind\{s<t_{i+1}\}-F_Q(t_i)\bigr)=:u_i ,
\]
and
\[
u_i-l_i=v\bigl(\ind\{t_i<s<t_{i+1}\}+Q(t_i<S<t_{i+1})\bigr),\qquad\|u_i-l_i\|_{P,2}\le\frac{\epsilon}{2}+K\,Q(t_i<S<t_{i+1})\le\epsilon .
\]
The functions $h_{t_i}$ at the finitely many partition points, and the two null functions at the infinite endpoints, are covered by brackets of size zero. Altogether
\[
N_{[\,]}\bigl(\epsilon,\Hcal,L_2(P)\bigr)\le2N+1\le\frac{8K^2}{\epsilon^2}+\frac{4K}{\epsilon}+3 .
\]
Since the logarithm of this bound is of order $\log(1/\epsilon)$ as $\epsilon$ tends to zero, the entropy integral converges.
\end{proof}

\begin{proof}[Proof of Theorem~\ref{thm:fclt}]
By Lemma~\ref{lem:bracketH} and Ossiander's theorem \citep[Theorem~2.5.6]{vandervaart1996weak}, $\Hcal$ is $P$-Donsker:
\[
\Gm h_{(\cdot)}\weak\mathbb G h_{(\cdot)}\qquad\text{in }\linf(\Rbar),
\]
where $\mathbb G$ is the $P$-Brownian bridge indexed by $\Hcal$, a tight centered Gaussian process with
\[
\Cov(\mathbb Gh_s,\mathbb Gh_t)=P(h_sh_t)-Ph_s\,Ph_t=P(h_sh_t)=C(s,t),
\]
because $Ph_t=0$. Expanding the product in \eqref{eq:C} and using
\[
\E_P\bigl[v^2\ind\{S\le s\}\ind\{S\le t\}\bigr]=G(s\wedge t),\qquad\E_P\bigl[v^2\ind\{S\le t\}\bigr]=G(t),\qquad\E_Pv^2=K^2
\]
gives the displayed formula for $C$. On the event that not all weights vanish, Lemma~\ref{lem:identity} gives
\[
\Zm=\frac{\Gm h_{(\cdot)}}{\Pn v},\qquad
\bigl\|\Zm-\Gm h_{(\cdot)}\bigr\|_\infty\le\bigl\|\Gm h_{(\cdot)}\bigr\|_\infty\,\Bigl|\frac{1}{\Pn v}-1\Bigr|\toP0,
\]
since the supremum of the empirical process is bounded in probability and $\Pn v\to1$ almost surely by the strong law; the exceptional event has probability $\Prob_P(v=0)^m$, which tends to zero. Hence
\[
\Zm\weak\GQ:=\mathbb Gh_{(\cdot)} .
\]

For \eqref{eq:brownian}, the process $t\mapsto\W(G(t))-F_Q(t)\W(K^2)$ is centered Gaussian with covariance
\[
G(s)\wedge G(t)-F_Q(s)G(t)-F_Q(t)G(s)+K^2F_Q(s)F_Q(t),
\]
using
\[
\Cov\bigl(\W(a),\W(b)\bigr)=a\wedge b,\qquad G\le K^2 ;
\]
as $G$ is nondecreasing,
\[
G(s)\wedge G(t)=G(s\wedge t),
\]
so the covariance is $C(s,t)$, and two centered Gaussian processes with the same covariance have the same law. For the remaining claims: the null functions at the infinite endpoints give $\GQ(\pm\infty)=0$;
\[
C(t,t)=\E_P\bigl[v^2(\ind\{S\le t\}-F_Q(t))^2\bigr]=\E_Q\bigl[v\,(\ind\{S\le t\}-F_Q(t))^2\bigr]\le\E_Qv=K^2,
\]
and when $v\le B$,
\[
C(t,t)\le B\,\E_Q\bigl[(\ind\{S\le t\}-F_Q(t))^2\bigr]=B\,F_Q(t)\bigl(1-F_Q(t)\bigr);
\]
if $F_Q$ is continuous then
\[
G(t)-G(t-)=\E_Q\bigl[v\ind\{S=t\}\bigr]=0,
\]
so $G$ is continuous and so are the paths in \eqref{eq:brownian}; and for $v\equiv1$ the clock is $F_Q$, $K^2=1$, and \eqref{eq:brownian} is the Brownian bridge composed with $F_Q$.
\end{proof}

\subsection{Proof of Proposition~\ref{prop:necessity}}\label{app:necessity}

\begin{lemma}\label{lem:converseCLT}
Let $\zeta,\zeta_1,\zeta_2,\dots$ be i.i.d.\ real random variables. If the normalized partial sums
\[
\frac{1}{\sqrt{m_k}}\sum_{i\le m_k}\zeta_i
\]
converge in distribution along some subsequence $m_k\to\infty$, then $\E\zeta^2<\infty$.
\end{lemma}

\begin{proof}
Let $\varphi$ be the characteristic function of $\zeta$ and $\hat\mu$ that of the limit, so that
\[
\varphi\bigl(u/\sqrt{m_k}\bigr)^{m_k}\longrightarrow\hat\mu(u)\qquad\text{for all }u .
\]
Fix $u_0>0$ with $|\hat\mu(u_0)|\ge1/2$, which exists by continuity of $\hat\mu$ at the origin. Then, for large $k$, using $1-x\le-\log x$ on the unit interval,
\[
m_k\Bigl(1-\bigl|\varphi(u_0/\sqrt{m_k})\bigr|^2\Bigr)\le-m_k\log\bigl|\varphi(u_0/\sqrt{m_k})\bigr|^2\longrightarrow-\log|\hat\mu(u_0)|^2\le\log4 .
\]
Let $\zeta'$ be an independent copy of $\zeta$, so that
\[
|\varphi(r)|^2=\E\cos\bigl(r(\zeta-\zeta')\bigr),\qquad\text{and}\qquad 1-\cos y\ge\frac{y^2}{3}\quad\text{for }|y|\le1 ;
\]
hence
\[
1-|\varphi(r)|^2=\E\bigl[1-\cos(r(\zeta-\zeta'))\bigr]\ge\frac{r^2}{3}\,\E\bigl[(\zeta-\zeta')^2\ind\{|\zeta-\zeta'|\le1/r\}\bigr].
\]
With $r=u_0/\sqrt{m_k}$ this gives
\[
\frac{u_0^2}{3}\,\E\bigl[(\zeta-\zeta')^2\ind\{|\zeta-\zeta'|\le\sqrt{m_k}/u_0\}\bigr]\le\log4+o(1),
\]
and monotone convergence yields
\[
\E(\zeta-\zeta')^2\le\frac{3\log4}{u_0^2}<\infty .
\]
Finally a finite second moment of the symmetrized variable implies one for $\zeta$ itself: otherwise the conditional second moment of $\zeta-\zeta'$ given $\zeta'$ would be infinite for every value of $\zeta'$.
\end{proof}

\begin{proof}[Proof of Proposition~\ref{prop:necessity}]
Suppose $\Zm(t)$ were tight. Then some subsequence converges in distribution, and since $\Pn v\to1$ almost surely, so does
\[
\sqrt m\,\Pn h_t=\Zm(t)\,\Pn v
\]
along it, by Slutsky's lemma. By Lemma~\ref{lem:converseCLT}, $Ph_t^2<\infty$. But
\begin{align*}
Ph_t^2&=\E_P\bigl[v^2(\ind\{S\le t\}-F_Q(t))^2\bigr]=\E_Q\bigl[v\,(\ind\{S\le t\}-F_Q(t))^2\bigr]\\
&\ge\min\{F_Q(t),1-F_Q(t)\}^2\,\E_Qv=\min\{F_Q(t),1-F_Q(t)\}^2K^2,
\end{align*}
because the squared centered indicator takes only the two values $F_Q(t)^2$ and $(1-F_Q(t))^2$. If $0<F_Q(t)<1$ and $K^2=\infty$ this is infinite, a contradiction.
\end{proof}

\begin{proof}[Proof of Corollary~\ref{cor:iff}]
The implication from (i) to (ii) is Theorem~\ref{thm:fclt}; from (ii) to (iii) is the continuous mapping theorem applied to a coordinate projection, which yields tightness at every $t$, in particular at an interior one, which exists because the score is nondegenerate under $Q$; from (iii) to (i) is Proposition~\ref{prop:necessity}.
\end{proof}

\subsection{Proof of Proposition~\ref{prop:bousquet}}\label{app:bousquet}

Let
\[
Z:=\sup_{t\in\Rbar}\Bigl|\sum_{i=1}^mh_t(Z_i)\Bigr|=\sup_{t\in\mathbb Q\cup\{\pm\infty\}}\Bigl|\sum_{i=1}^mh_t(Z_i)\Bigr|,
\]
the reduction to a countable index set being justified by right-continuity of $t\mapsto h_t(z)$. Each summand satisfies
\[
\E h_t(Z_i)=0,\qquad|h_t|\le v\le B,\qquad\E h_t(Z_i)^2\le\sigma_*^2 .
\]
Bousquet's inequality in the form of \citet[Theorem~12.5]{boucheron2013concentration}, applied to the class of the functions $h_t/B$ and $-h_t/B$, whose members are bounded above by one, states that for all $u>0$,
\[
\Prob\Bigl(Z\ge\E Z+\sqrt{2\bigl(m\sigma_*^2+2B\,\E Z\bigr)u}+\tfrac{Bu}{3}\Bigr)\le e^{-u}.
\]
By the bracketing maximal inequality \citep[Theorem~2.14.2]{vandervaart1996weak} with envelope $v$, whose $L_2(P)$-norm is $K$, and Lemma~\ref{lem:bracketH},
\begin{align*}
\E Z=\sqrt m\,\E\bigl\|\Gm h_{(\cdot)}\bigr\|_\infty&\le\sqrt m\,C\,K\int_0^1\sqrt{1+\log N_{[\,]}(\epsilon K,\Hcal,L_2(P))}\,\dd\epsilon\\
&\le\sqrt m\,C\,K\int_0^1\sqrt{1+\log\Bigl(\frac{8}{\epsilon^2}+\frac4\epsilon+3\Bigr)}\,\dd\epsilon=:C_1K\sqrt m,
\end{align*}
where $C$ is universal and the integral is a numerical constant (Lemma~\ref{lem:bracketH} gives a bound on the bracketing numbers at scale $\epsilon K$ that is free of $K$). Dividing by $m$, with probability at least $1-e^{-u}$,
\[
\sup_t|\Pn h_t|\le\frac{C_1K}{\sqrt m}+\sqrt{\frac{2(\sigma_*^2+2C_1BK/\sqrt m)u}{m}}+\frac{Bu}{3m}.
\]
By Hoeffding's inequality for the weights, which lie in $[0,B]$,
\[
\Prob\bigl(\Pn v\ge1-\eta\bigr)\ge1-e^{-2m\eta^2/B^2}.
\]
On the intersection of the two events, Lemma~\ref{lem:identity} gives the claim. \qed

\section{Proofs for Section~\ref{sec:coverage}}\label{app:coverage}

Throughout this section $K^2<\infty$ and $F_Q$ is continuous. Let
\[
\rho(s,t):=\|h_s-h_t\|_{P,2}.
\]
Since $\Hcal$ is Donsker, the empirical process is asymptotically uniformly $\rho$-equicontinuous in probability \citep[Theorem~1.5.7 and Section~2.1]{vandervaart1996weak}: for all $\epsilon,\eta>0$ there is $\zeta>0$ with
\begin{equation}\label{eq:equicont}
\limsup_{m\to\infty}\Prob^*\Bigl(\sup_{\rho(s,t)<\zeta}\bigl|\Gm h_s-\Gm h_t\bigr|>\epsilon\Bigr)<\eta ,
\end{equation}
and the same holds for $\Zm$, which is the empirical process divided by $\Pn v\to1$.

\begin{lemma}\label{lem:rho}
For $s<t$,
\[
\rho(s,t)^2\le\{G(t)-G(s)\}+K^2\{F_Q(t)-F_Q(s)\}^2 .
\]
Moreover, for every $\zeta>0$ there is $\kappa>0$ such that $F_Q(t)-F_Q(s)<\kappa$ implies $\rho(s,t)<\zeta$.
\end{lemma}

\begin{proof}
With $d=F_Q(t)-F_Q(s)$ we have
\[
h_t-h_s=v\bigl(\ind\{s<S\le t\}-d\bigr),
\]
so
\[
\rho(s,t)^2=(1-d)^2\{G(t)-G(s)\}+d^2\bigl\{K^2-(G(t)-G(s))\bigr\}\le\{G(t)-G(s)\}+d^2K^2 .
\]
The finite measure
\[
\mu_G(A)=\E_Q\bigl[v\ind\{S\in A\}\bigr]
\]
is absolutely continuous with respect to the law of $S$ under $Q$, hence uniformly so: for every $\eta>0$ there is $\kappa>0$ such that any set of $Q$-probability below $\kappa$ has $\mu_G$-mass below $\eta$. Apply this with $A=(s,t]$, whose $Q$-probability is $d$ by continuity of $F_Q$, and choose $\eta$ and $\kappa$ with $\eta+K^2\kappa^2<\zeta^2$.
\end{proof}

\begin{lemma}\label{lem:maxw}
If $K^2<\infty$ then $\sqrt m\max_{i\le m}\hat w_i\toP0$.
\end{lemma}

\begin{proof}
The largest normalized weight is $\max_iv_i/\sum_jv_j$, and $\sum_jv_j/m\to1$ almost surely. For $\epsilon>0$,
\[
\Prob\bigl(\max_iv_i>\epsilon\sqrt m\bigr)\le m\,\Prob(v>\epsilon\sqrt m)\le\frac{1}{\epsilon^2}\,\E\bigl[v^2\ind\{v>\epsilon\sqrt m\}\bigr]\longrightarrow0
\]
by dominated convergence.
\end{proof}

\begin{proof}[Proof of Theorem~\ref{thm:coverage}]
Work on the event that not all weights vanish. For $\beta\le\beta_1<1$ the quantile is finite, since $\Fhat$ equals one at the largest score. The weighted empirical distribution function is a right-continuous step function whose jump at any point is the total weight of the scores located there, and two scores with positive weight tie with probability zero: for every $s$,
\[
Q(S=s)=\E_P\bigl[v\ind\{S=s\}\bigr]=0
\]
by continuity of $F_Q$, which forces $v=0$ $P$-almost surely on every $P$-atom of $S$. Hence, almost surely,
\[
r_m(\beta):=\Fhat\bigl(\tauhat(\beta)\bigr)-\beta\in\bigl[0,\max_i\hat w_i\bigr],
\]
and therefore
\[
\Cm(\beta)=\beta-F_Q(\tauhat(\beta))=\bigl\{\Fhat(\tauhat(\beta))-F_Q(\tauhat(\beta))\bigr\}-r_m(\beta),
\]
that is,
\[
\sqrt m\,\Cm(\beta)=\Zm\bigl(\tauhat(\beta)\bigr)-\sqrt m\,r_m(\beta).
\]
The remainder is uniformly negligible by Lemma~\ref{lem:maxw}. It remains to show that
\[
\sup_{\beta\in[\beta_0,\beta_1]}\bigl|\Zm(\tauhat(\beta))-\Zm(q_Q(\beta))\bigr|\toP0 .
\]
By continuity of $F_Q$ we have $F_Q(q_Q(\beta))=\beta$, so
\[
\bigl|F_Q(\tauhat(\beta))-F_Q(q_Q(\beta))\bigr|=\bigl|F_Q(\tauhat(\beta))-\beta\bigr|\le\|\Fhat-F_Q\|_\infty+\max_i\hat w_i\toP0
\]
uniformly in $\beta$, using
\[
\|\Fhat-F_Q\|_\infty=m^{-1/2}\|\Zm\|_\infty=O_p(m^{-1/2})
\]
from Theorem~\ref{thm:fclt}. By Lemma~\ref{lem:rho} the $\rho$-distance between $\tauhat(\beta)$ and $q_Q(\beta)$ tends to zero in probability uniformly in $\beta$, and \eqref{eq:equicont} gives the claim. The weak convergence of the coverage process then follows from
\[
\Zm\circ q_Q\weak\GQ\circ q_Q,
\]
composition with a fixed map being continuous on $\linf$, and Slutsky's lemma. The representation of $\Lb$ is \eqref{eq:brownian} evaluated at $q_Q(\beta)$, with
\[
\tilde G(\beta)=G(q_Q(\beta))=\E_Q\bigl[v\ind\{S\le q_Q(\beta)\}\bigr]=\E_Q\bigl[v\ind\{F_Q(S)\le\beta\}\bigr],
\]
the two events differing by a $Q$-null set because $F_Q$ is continuous; the variance formula is $C(q,q)$ with $F_Q(q)=\beta$.
\end{proof}

\begin{proof}[Proof of Proposition~\ref{prop:testweight}]
Let $V=\sum_jv_j$ and fix $\beta\in[\beta_0,\beta_1]$. Since $\tauhat$ is nondecreasing and the test weight only raises the level,
\[
0\le P_e(\tauhat(\beta))-P_e^{\rm split}(D;\beta)=Q\Bigl(\tauhat(\beta)<S\le\tauhat\bigl(\beta(1+v(X)/V)\bigr)\Bigr).
\]
Let
\[
r(a):=\E_Q\bigl[v\ind\{v>a\}\bigr]=\E_P\bigl[v^2\ind\{v>a\}\bigr]\longrightarrow0\qquad(a\to\infty),
\]
and choose
\[
a_m=\sqrt m\,\epsilon_m,\qquad\epsilon_m:=\sqrt{r(m^{1/4})}+m^{-1/8},
\]
so that, for large $m$,
\[
\epsilon_m\to0,\qquad a_m\ge m^{3/8}\to\infty,\qquad r(a_m)\le r(m^{1/4})\le\epsilon_m^2 .
\] Splitting according to whether $v(X)\le a_m$ and using the monotonicity of $\tauhat$ on the first part,
\[
P_e(\tauhat(\beta))-P_e^{\rm split}(D;\beta)\le\Bigl[F_Q\bigl(\tauhat(\beta_a)\bigr)-F_Q\bigl(\tauhat(\beta)\bigr)\Bigr]+Q\bigl(v(X)>a_m\bigr),\qquad\beta_a:=\beta\Bigl(1+\frac{a_m}{V}\Bigr),
\]
on the event that $\beta_a$ does not exceed $\beta_1':=(1+\beta_1)/2$, whose probability tends to one since $a_m/V\toP0$. The second term satisfies
\[
Q\bigl(v(X)>a_m\bigr)\le\frac{r(a_m)}{a_m}\le\frac{\epsilon_m^2}{\sqrt m\,\epsilon_m}=o(m^{-1/2}).
\]
For the first term, by \eqref{eq:Cm},
\[
F_Q(\tauhat(\beta_a))-F_Q(\tauhat(\beta))=(\beta_a-\beta)+\Cm(\beta)-\Cm(\beta_a)=\frac{\beta a_m}{V}+\Cm(\beta)-\Cm(\beta_a).
\]
Here
\[
\sqrt m\,\frac{\beta a_m}{V}=\beta\,\epsilon_m\,\frac{m}{V}\toP0 .
\]
By Theorem~\ref{thm:coverage} on $[\beta_0,\beta_1']$,
\[
\sqrt m\{\Cm(\beta)-\Cm(\beta_a)\}=\Zm(q_Q(\beta))-\Zm(q_Q(\beta_a))+o_p(1)
\]
uniformly, and since
\[
F_Q(q_Q(\beta_a))-F_Q(q_Q(\beta))=\beta_a-\beta\le\frac{\beta_1a_m}{V}\toP0
\]
uniformly, Lemma~\ref{lem:rho} and \eqref{eq:equicont} show that the difference is uniformly negligible. All bounds are uniform in $\beta\in[\beta_0,\beta_1]$, which proves the proposition.
\end{proof}

\begin{proof}[Proof of Corollary~\ref{cor:asympdkw}]
The convergence of $\|\Zm\|_\infty$ is the continuous mapping theorem. The Borell--TIS inequality \citep[Theorem~2.1.1]{adler2007random}, applied to the centered Gaussian process
\[
(t,\pm1)\longmapsto\pm\GQ(t)
\]
on the doubled index set (indexed through a countable dense subset and extended by right-continuity), whose supremum is $\|\GQ\|_\infty$ and whose variances are at most $\sigma_*^2$, gives
\[
\Prob\bigl(\|\GQ\|_\infty\ge\E\|\GQ\|_\infty+u\bigr)\le e^{-u^2/(2\sigma_*^2)}.
\]
For the expectation, the canonical metric is
\[
d(s,t)^2:=\E\bigl(\GQ(s)-\GQ(t)\bigr)^2=\rho(s,t)^2\le\{G(t)-G(s)\}+K^2\{F_Q(t)-F_Q(s)\}^2\qquad(s<t)
\]
by Lemma~\ref{lem:rho}. Apply Lemma~\ref{lem:partition} to $\mu_G$ (mass $K^2$) with $\eta=\epsilon^2/2$ and to the law of $S$ under $Q$ (mass $1$) with $\eta=\epsilon/(\sqrt2K)$, and merge the partition points: the extended line is covered by at most
\[
\frac{2K^2}{\epsilon^2}+\frac{\sqrt2K}{\epsilon}+1
\]
open intervals, on each of which $G$ increases by at most $\epsilon^2/2$ and the squared $F_Q$-increment times $K^2$ is at most $\epsilon^2/2$, so that each has $d$-diameter at most $\epsilon$; together with the endpoints, the covering numbers satisfy
\[
N(\epsilon)\le\frac{4K^2}{\epsilon^2}+\frac{2\sqrt2K}{\epsilon}+2\qquad\text{for every }\epsilon>0,
\]
and the $d$-diameter of the index set is at most $2\sigma_*\le2K$. Dudley's entropy bound \citep{dudley1967sizes}, in the form
\[
\E\sup_t\GQ(t)\le c\int_0^{\operatorname{diam}}\sqrt{\log N(\epsilon)}\,\dd\epsilon
\]
with a universal $c$ \citep[Theorem~1.3.3]{adler2007random}, gives
\[
\E\|\GQ\|_\infty\le c\int_0^{2K}\sqrt{\log\Bigl(\frac{4K^2}{\epsilon^2}+\frac{2\sqrt2K}{\epsilon}+2\Bigr)}\,\dd\epsilon
=cK\int_0^2\sqrt{\log\Bigl(\frac4{x^2}+\frac{2\sqrt2}{x}+2\Bigr)}\,\dd x=:c_0K .
\]
The limsup statement follows from the weak convergence of $\|\Zm\|_\infty$, the tail bound with
\[
u=K\sqrt{2\log(1/\delta)}\ge\sigma_*\sqrt{2\log(1/\delta)},
\]
and the portmanteau theorem. The statement for the coverage error follows from
\[
\sqrt m\sup_\beta|\Cm(\beta)|\le\|\Zm\|_\infty+\sqrt m\max_i\hat w_i
\]
and Lemma~\ref{lem:maxw}.
\end{proof}

\subsection{Proof of Theorem~\ref{thm:bootstrap}}\label{app:bootstrap}

\begin{proof}
\emph{Multiplier process.} Write
\[
\Gm^\xi f:=\frac{1}{\sqrt m}\sum_{i=1}^m\xi_if(Z_i).
\]
Since $\Hcal$ is Donsker with envelope $v$ in $L_2(P)$ and the multipliers satisfy the moment condition, the conditional multiplier central limit theorem \citep[Theorem~2.9.6]{vandervaart1996weak} gives
\[
\Gm^\xi h_{(\cdot)}\weak_\xi\GQ\qquad\text{in }\linf(\Rbar)
\]
(the centering by $Pf$ in that theorem is absent here because $Ph_t=0$). The process \eqref{eq:multiplier} equals
\[
\hat{\mathbb Z}^\xi_m(t)=\frac{1}{\Pn v}\Bigl\{\Gm^\xi h_t-\bigl(\Fhat(t)-F_Q(t)\bigr)\,\Gm^\xi v\Bigr\},
\]
because
\[
v_i\bigl(\ind\{S_i\le t\}-\Fhat(t)\bigr)=h_t(Z_i)-v_i\bigl(\Fhat(t)-F_Q(t)\bigr).
\]
The correction satisfies
\[
\sup_t\bigl|\Fhat(t)-F_Q(t)\bigr|\,\bigl|\Gm^\xi v\bigr|=O_p(m^{-1/2})\,O_p(1)\toP0,
\]
since $\E(\Gm^\xi v)^2=K^2$, and $\Pn v\to1$ almost surely. Conditional weak convergence is preserved under perturbations that tend to zero in (unconditional) outer probability: for $f$ in the unit ball of bounded Lipschitz functions,
\[
\bigl|\E_\xi f(A_m)-\E_\xi f(B_m)\bigr|\le\E_\xi\bigl(\|A_m-B_m\|_\infty\wedge2\bigr),
\]
whose unconditional expectation tends to zero by dominated convergence, so the conditional expectation tends to zero in probability. This gives the first claim.

\emph{Composition.} Unconditionally, $\Gm^\xi h_{(\cdot)}\weak\GQ$ as well \citep[Theorem~2.9.2]{vandervaart1996weak}, so the multiplier process is asymptotically uniformly $\rho$-equicontinuous in probability. As in the proof of Theorem~\ref{thm:coverage}, the $\rho$-distance between $\tauhat(\beta)$ and $q_Q(\beta)$ tends to zero in probability uniformly in $\beta$, hence
\[
\sup_{\beta\in[\beta_0,\beta_1]}\bigl|\hat{\mathbb Z}^\xi_m(\tauhat(\beta))-\hat{\mathbb Z}^\xi_m(q_Q(\beta))\bigr|\toP0
\]
unconditionally, and therefore the composition with $\tauhat$ has the same conditional limit as the composition with $q_Q$, namely $\GQ\circ q_Q=\Lb$ (composition with a fixed map is continuous on $\linf$).

\emph{Efron's bootstrap.} Let the pairs $(S_i^*,v_i^*)$, $i\le m$, be drawn with replacement from the calibration pairs, let $\Pn^*$ be their empirical measure, and let $\Fhat^*$ be the self-normalized weighted distribution function computed from them. The identity of Lemma~\ref{lem:identity} applied to the bootstrap sample, with $\Pn$ in the role of $P$ and with
\[
\hat h_t:=v\bigl(\ind\{S\le t\}-\Fhat(t)\bigr),\qquad\Pn\hat h_t=0,
\]
gives
\[
\sqrt m\bigl(\Fhat^*(t)-\Fhat(t)\bigr)=\frac{\sqrt m(\Pn^*-\Pn)\hat h_t}{\Pn^*v}=\frac{\Gm^*h_t-(\Fhat(t)-F_Q(t))\,\Gm^*v}{\Pn^*v},
\]
where $\Gm^*$ denotes the bootstrap empirical process $\sqrt m(\Pn^*-\Pn)$.
Since $\Hcal\cup\{v\}$ is Donsker with square-integrable envelope, the bootstrap empirical process converges conditionally in probability, and indeed outer almost surely \citep[Theorems~3.6.1 and~3.6.2]{vandervaart1996weak}; the correction term and the denominator are handled exactly as above, and so is the composition with $\tauhat$.
\end{proof}

\begin{lemma}[uniform consistency of the variance estimator]\label{lem:sigmahatunif}
Assume $K^2<\infty$ and that $F_Q$ is continuous. Then
\[
\sup_{\beta\in[\beta_0,\beta_1]}\bigl|\hat\sigma_m^2(\beta)-\sigma^2(\beta)\bigr|\toP0 .
\]
\end{lemma}

\begin{proof}
Write $\hat\sigma_m^2(\beta)=\Pn[v^2(\ind\{S\le\tauhat(\beta)\}-\beta)^2]/(\Pn v)^2$. The class
\[
\Fcal_2:=\bigl\{v^2\bigl(\ind\{S\le t\}-c\bigr)^2:\ t\in\Rbar,\ c\in[0,1]\bigr\}
\]
has envelope $v^2\in L_1(P)$ and finite $L_1(P)$-bracketing numbers: for fixed $c$ the functions $v^2\{(1-c)^2\ind\{S\le t\}+c^2\ind\{S>t\}\}$ are bracketed by the partition lemma applied to the measure $A\mapsto\E_P[v^2\ind\{S\in A\}]$, and the dependence on $c$ is Lipschitz with constant $2v^2$. Hence $\Fcal_2$ is $P$-Glivenko--Cantelli, and
\[
\sup_{\beta}\Bigl|\Pn\bigl[v^2(\ind\{S\le\tauhat(\beta)\}-\beta)^2\bigr]-P\bigl[v^2(\ind\{S\le\tauhat(\beta)\}-\beta)^2\bigr]\Bigr|\le\|\Pn-P\|_{\Fcal_2}\to0\quad\text{a.s.}
\]
Next,
\[
P\bigl[v^2(\ind\{S\le t\}-\beta)^2\bigr]=(1-\beta)^2G(t)+\beta^2\bigl(K^2-G(t)\bigr),
\]
and $G=\tilde G\circ F_Q$ because $G$ is constant on every interval on which $F_Q$ is constant (such an interval is $Q$-null). Since $\tilde G$ is continuous on $[0,1]$ (proof of Corollary~\ref{cor:simultaneous}) and
\[
\sup_{\beta\in[\beta_0,\beta_1]}\bigl|F_Q(\tauhat(\beta))-\beta\bigr|=\sup_\beta|\Cm(\beta)|\toP0
\]
by Theorem~\ref{thm:coverage}, we get $\sup_\beta|G(\tauhat(\beta))-\tilde G(\beta)|\toP0$, hence the numerator converges uniformly to $\sigma^2(\beta)$; finally $\Pn v\to1$ a.s.
\end{proof}

\begin{proof}[Proof of Corollary~\ref{cor:simultaneous}]
Let $M:=\sup_{\beta\in\Bcal}\Lb(\beta)$. The process $\Lb$ has continuous paths on $[\beta_0,\beta_1]$: indeed
\[
\Lb(\beta)=\W(\tilde G(\beta))-\beta\,\W(K^2),\qquad\tilde G(\beta)=\E_Q\bigl[v\ind\{F_Q(S)\le\beta\}\bigr],
\]
and $\tilde G$ is continuous because $F_Q(S)$ is uniform under $Q$ and the finite measure
\[
A\mapsto\E_Q\bigl[v\ind\{F_Q(S)\in A\}\bigr]
\]
is absolutely continuous with respect to Lebesgue measure. Its variance is bounded below on $\Bcal$,
\[
\sigma^2(\beta)\ge K^2\min\{\beta_0,1-\beta_1\}^2>0,
\]
by the argument in the proof of Proposition~\ref{prop:necessity}. By Ylvisaker's theorem \citep{ylvisaker1968note,tsirelson1975density}, see \citet{lifshits1995gaussian}, the law of the supremum of a continuous Gaussian process over a compact set on which the variance is bounded away from zero has no atoms (for a finite set $\Bcal$ this is elementary). Moreover the distribution function of $M$ is log-concave, since Gaussian measures are log-concave \citep{borell1974convex} and the event that $M$ is at most $c$ is a convex set of paths; a log-concave distribution function is strictly increasing wherever it lies strictly between zero and one. Hence the distribution function of $M$ is continuous and strictly increasing at $c_{1-\delta}$. Theorem~\ref{thm:bootstrap} and the continuous mapping theorem (the supremum over the closed set $\Bcal$ is continuous on $\linf[\beta_0,\beta_1]$) give
\[
\sup_{\beta\in\Bcal}\hat\Lb^\xi_m(\beta)\weak_\xi M,
\]
and conditional weak convergence to a law whose distribution function is continuous and strictly increasing at $c_{1-\delta}$ implies convergence of the conditional quantiles, $\hat c\toP c_{1-\delta}$ \citep[Lemma~21.2 and Lemma~23.3]{vandervaart1998asymptotic}. For the coverage statement, note that for $\beta\in\Bcal$ and large $m$, so that the inflated level stays below $\beta_1'<1$,
\[
P_e\bigl(\tauhat(\beta+\hat c/\sqrt m)\bigr)-(1-\beta)=\Cm(\beta+\hat c/\sqrt m)-\frac{\hat c}{\sqrt m},
\]
so the failure event is
\[
\Bigl\{\sup_{\beta\in\Bcal}\sqrt m\,\Cm(\beta+\hat c/\sqrt m)>\hat c\Bigr\}.
\]
By Theorem~\ref{thm:coverage} on $[\beta_0,\beta_1']$ and Lemma~\ref{lem:rho} with \eqref{eq:equicont},
\[
\sup_{\beta\in\Bcal}\bigl|\sqrt m\,\Cm(\beta+\hat c/\sqrt m)-\sqrt m\,\Cm(\beta)\bigr|\toP0,
\]
hence the failure probability converges to $\Prob(M>c_{1-\delta})=\delta$ by Slutsky's lemma and the continuity of the law of $M$ at $c_{1-\delta}$. This proves (a). The two-sided and test-weight versions are identical, using Proposition~\ref{prop:testweight} for the latter.

\emph{Part (b).} Since $\sigma$ is continuous and bounded below on $[\beta_0,\beta_1]$, division by $\sigma$ is a continuous linear map of $\linf[\beta_0,\beta_1]$ into itself, so Theorem~\ref{thm:bootstrap} and the continuous mapping theorem give $\sup_\Bcal\hat\Lb^\xi_m/\sigma\weak_\xi M^s:=\sup_\Bcal\Lb/\sigma$, a supremum of a continuous Gaussian process with unit variance, whose law is atomless with a log-concave, hence strictly increasing, distribution function, as above. By Lemma~\ref{lem:sigmahatunif},
\[
\Delta_m:=\sup_{\beta\in\Bcal}\Bigl|\frac{\hat\Lb^\xi_m(\beta)}{\hat\sigma_m(\beta)}-\frac{\hat\Lb^\xi_m(\beta)}{\sigma(\beta)}\Bigr|\le\sup_\Bcal|\hat\Lb^\xi_m|\,\sup_\Bcal\Bigl|\frac{1}{\hat\sigma_m}-\frac1\sigma\Bigr|\toP0,
\]
and for every $\epsilon>0$ the conditional probability $\Prob_\xi(\Delta_m>\epsilon)$ tends to zero in probability, by Markov's inequality applied to $\E\,\Prob_\xi(\Delta_m>\epsilon)=\Prob(\Delta_m>\epsilon)$. Hence, for every $x$,
\begin{multline*}
\Prob_\xi\Bigl(\sup_\Bcal\hat\Lb^\xi_m/\sigma\le x-\epsilon\Bigr)-\Prob_\xi(\Delta_m>\epsilon)\le\Prob_\xi\Bigl(\sup_\Bcal\hat\Lb^\xi_m/\hat\sigma_m\le x\Bigr)\\
\le\Prob_\xi\Bigl(\sup_\Bcal\hat\Lb^\xi_m/\sigma\le x+\epsilon\Bigr)+\Prob_\xi(\Delta_m>\epsilon),
\end{multline*}
so the conditional distribution function of the studentized supremum converges in probability to that of $M^s$ at every point, and $\hat c^{\,s}\toP c^s_{1-\delta}$, the $(1-\delta)$-quantile of $M^s$, as in (a). For the coverage statement put $\beta'_\beta:=\beta+\hat c^{\,s}\hat\sigma_m(\beta)/\sqrt m$; then $\sup_\Bcal|\beta'_\beta-\beta|\toP0$, the failure event is
\[
\Bigl\{\sup_{\beta\in\Bcal}\bigl[\sqrt m\,\Cm(\beta'_\beta)-\hat c^{\,s}\hat\sigma_m(\beta)\bigr]>0\Bigr\},
\]
and $\sup_\Bcal|\sqrt m\,\Cm(\beta'_\beta)-\sqrt m\,\Cm(\beta)|\toP0$ as in (a), while $\sup_\Bcal|\hat c^{\,s}\hat\sigma_m(\beta)-c^s_{1-\delta}\sigma(\beta)|\toP0$ by Lemma~\ref{lem:sigmahatunif}. Hence the failure probability converges to $\Prob(\sup_\Bcal(\Lb-c^s_{1-\delta}\sigma)>0)=\Prob(M^s>c^s_{1-\delta})=\delta$, the law of $M^s$ being continuous at $c^s_{1-\delta}$.

\emph{Part (c).} Let $\pi_m(\beta)$ be a nearest point of $\Bcal_m$ to $\beta$, so that $\sup_\Bcal|\beta-\pi_m(\beta)|\to0$. Then $0\le\sup_\Bcal\hat\Lb^\xi_m-\sup_{\Bcal_m}\hat\Lb^\xi_m\le\sup_{\beta\in\Bcal}|\hat\Lb^\xi_m(\beta)-\hat\Lb^\xi_m(\pi_m(\beta))|$, which tends to zero in probability by the asymptotic equicontinuity of the multiplier process $\hat{\mathbb Z}^\xi_m$ over $\Hcal$ (which holds in outer probability, Theorem~2.9.6 of \citet{vandervaart1996weak}), the uniform continuity of $\tilde G$ and Lemma~\ref{lem:rho}, exactly as in the proof of Theorem~\ref{thm:coverage}. The argument of (b) then shows that the conditional quantiles computed over $\Bcal_m$ converge to the same limits $c_{1-\delta}$, respectively $c^s_{1-\delta}$, and the coverage statements over all of $\Bcal$ follow as before, since they only use these limits.
\end{proof}

\section{Proofs for Section~\ref{sec:estimated}}\label{app:estimated}

\begin{proof}[Proof of Theorem~\ref{thm:estimated}]
Work on the event that $\hat v$ belongs to $\Vcal$ and $\Pn\hat v>0$, whose probability tends to one. Write
\[
\tilde k_t:=\hat v\,\bigl(\ind\{S\le t\}-F_{\hat Q}(t)\bigr)\in\Hcal_\Vcal .
\]
By definition of the pseudo-target,
\[
P\tilde k_t=P\hat v\,F_{\hat Q}(t)-F_{\hat Q}(t)\,P\hat v=0,
\]
with $P$ integrating over $Z$ for fixed $\hat v$. Hence, exactly as in Lemma~\ref{lem:identity},
\[
\sqrt m\bigl(\Fchk(t)-F_{\hat Q}(t)\bigr)=\frac{\sqrt m\,\Pn\tilde k_t}{\Pn\hat v}=\frac{\Gm\tilde k_t}{\Pn\hat v},
\]
where $\Gm\tilde k_t$ denotes the centered and scaled empirical measure of the realized function. Now
\[
\tilde k_t-h_t=(\hat v-v)\bigl(\ind\{S\le t\}-F_Q(t)\bigr)-\hat v\bigl(F_{\hat Q}(t)-F_Q(t)\bigr),
\]
and, since $Pv=1$,
\begin{align*}
\sup_t\bigl|F_{\hat Q}(t)-F_Q(t)\bigr|&=\sup_t\Bigl|P\Bigl[\Bigl(\frac{\hat v}{P\hat v}-v\Bigr)\ind\{S\le t\}\Bigr]\Bigr|\le\Bigl\|\frac{\hat v}{P\hat v}-v\Bigr\|_{P,1}\\
&\le\frac{\|\hat v-v\|_{P,1}+|P\hat v-1|}{P\hat v}\le\frac{2\|\hat v-v\|_{P,2}}{P\hat v}\toP0,
\end{align*}
using
\[
|P\hat v-1|=|P(\hat v-v)|\le\|\hat v-v\|_{P,1},\qquad P\hat v\toP1 .
\]
Therefore
\[
\sup_t\|\tilde k_t-h_t\|_{P,2}\le\|\hat v-v\|_{P,2}+\|V\|_{P,2}\,\sup_t\bigl|F_{\hat Q}(t)-F_Q(t)\bigr|\toP0 .
\]
Both $\tilde k_t$ and $h_t$ lie in the Donsker class $\Hcal_\Vcal$, whose empirical process is asymptotically equicontinuous with respect to the $L_2(P)$ seminorm, which dominates the centered seminorm; as this equicontinuity is a statement about
\[
\sup\bigl\{|\Gm(f-g)|:\ f,g\in\Hcal_\Vcal,\ \|f-g\|_{P,2}<\zeta\bigr\},
\]
it applies to the random pair $(\tilde k_t,h_t)$ regardless of how $\hat v$ depends on the sample, and gives
\[
\sup_t\bigl|\Gm\tilde k_t-\Gm h_t\bigr|\toP0 .
\]
Finally
\[
\Pn\hat v-1=(\Pn-P)(\hat v-v)+(\Pn-P)v+(P\hat v-1)\toP0,
\]
the first term by the same equicontinuity (take $t=+\infty$ and $c=0$ in $\Hcal_\Vcal$). Combining with
\[
\Zm=\frac{\Gm h_{(\cdot)}}{\Pn v}
\]
and the boundedness in probability of $\|\Gm h_{(\cdot)}\|_\infty$ proves the theorem.
\end{proof}

\begin{proof}[Proof of Lemma~\ref{lem:parametric}]
Let $\epsilon\in(0,1]$. \emph{Brackets for $\Vcal$.} By \citet[Theorem~2.7.11]{vandervaart1996weak}, $\Vcal$ is covered by
\[
N_1\le C_d\Bigl(\frac{\operatorname{diam}\Theta\,\|\dot V\|_{P,2}}{\epsilon}\Bigr)^d
\]
brackets $[l,u']$ of the form $[v_{\theta_j}-\epsilon'\dot V,\ v_{\theta_j}+\epsilon'\dot V]$ with $\|u'-l\|_{P,2}\le\epsilon$; since every member of $\Vcal$ is nonnegative, we may replace $l$ by its positive part, so that
\[
0\le l\le u\le u',\qquad l\le V .
\]
\emph{Brackets for the indicator part.} Let $\mathcal G$ be the class of the functions
\[
g_{t,c}=\ind\{s\le t\}-c,\qquad t\in\Rbar,\ c\in[0,1],
\]
and let $\mu:=V^2\cdot P$, a finite measure of mass $\|V\|_{P,2}^2$. Apply Lemma~\ref{lem:partition} to the image of $\mu$ under $S$ with $\eta=\epsilon^2/4$, and take a grid $0=c_0<\dots<c_L=1$ of mesh $\epsilon/(2\|V\|_{P,2})$. For $t\in[t_i,t_{i+1})$ and $c\in[c_k,c_{k+1}]$,
\[
\ind\{s\le t_i\}-c_{k+1}\ \le\ g_{t,c}\ \le\ \ind\{s<t_{i+1}\}-c_k,
\]
a bracket of $L_2(\mu)$-size at most
\[
\sqrt{2\mu(t_i<S<t_{i+1})+2(c_{k+1}-c_k)^2\|V\|_{P,2}^2}\le\epsilon ;
\]
with the partition points this gives
\[
N_2\le2\Bigl(\frac{4\|V\|_{P,2}^2}{\epsilon^2}+1\Bigr)\Bigl(\frac{2\|V\|_{P,2}}{\epsilon}+1\Bigr)
\]
brackets $[l_g,u_g]$, all with values in $[-1,1]$. \emph{Product brackets.} For $u\in[l,u']$ and $g\in[l_g,u_g]$ as above, write $ug=lg+(u-l)g$ with $l\ge0$ and $|g|\le1$, so that
\[
l\,l_g-(u'-l)\ \le\ ug\ \le\ l\,u_g+(u'-l),
\]
a bracket of $L_2(P)$-size at most
\[
\|l(u_g-l_g)\|_{P,2}+2\|u'-l\|_{P,2}\le\|V(u_g-l_g)\|_{P,2}+2\epsilon=\|u_g-l_g\|_{L_2(\mu)}+2\epsilon\le3\epsilon .
\]
Hence
\[
N_{[\,]}\bigl(3\epsilon,\Hcal_\Vcal,L_2(P)\bigr)\le N_1N_2\le C\epsilon^{-(d+3)},
\]
and the bracketing integral is finite.
\end{proof}

\begin{proof}[Proof of Theorem~\ref{thm:parametric}]
By Theorem~\ref{thm:estimated} it suffices to show that
\[
\sqrt m\,(F_{\hat Q}-F_Q)=\sqrt m\,(\hat\theta-\theta_0)^\top\gamma+o_p(1)\qquad\text{uniformly in }t,
\]
and to identify the joint limit. Since $Pv=1$,
\[
F_{\hat Q}(t)-F_Q(t)=\frac{P\bigl[(v_{\hat\theta}-v)\bigl(\ind\{S\le t\}-F_Q(t)\bigr)\bigr]}{Pv_{\hat\theta}},
\]
and $L_1(P)$-differentiability gives
\[
v_{\hat\theta}-v=(\hat\theta-\theta_0)^\top\dot v+R,\qquad\|R\|_{P,1}=o_p(\|\hat\theta-\theta_0\|)=o_p(m^{-1/2}).
\]
Since the centered indicator is bounded by one, the numerator equals
\[
(\hat\theta-\theta_0)^\top\gamma(t)+O_p(\|R\|_{P,1})
\]
uniformly in $t$, and $Pv_{\hat\theta}\toP1$. This gives the expansion, with
\[
\sup_t\|\gamma(t)\|\le\|\dot v\|_{P,1}.
\]
The class $\Hcal\cup\{\psi_1,\dots,\psi_d\}$ is Donsker (a finite union of Donsker classes), so
\[
(\Gm h_{(\cdot)},\Gm\psi)\weak(\GQ,\mathbb G\psi)\qquad\text{jointly, with}\qquad
\Cov\bigl(\GQ(t),\mathbb G\psi\bigr)=P[h_t\psi]-Ph_t\,P\psi=P[h_t\psi],
\]
and conditioning on $X$ gives
\[
P[h_t\psi]=\E_Q\bigl[(\ind\{S\le t\}-F_Q(t))\,\psi(X)\bigr]=\E_Q\bigl[\{F_{S|X}(t)-F_Q(t)\}\,\psi(X)\bigr].
\]
Since $\xi_m$ is independent of the calibration sample, the triple $(\Gm h_{(\cdot)},\Gm\psi,\xi_m)$ converges jointly to $(\GQ,\mathbb G\psi,\mathcal N)$ with $\mathcal N$ independent of the first two. The continuous mapping theorem applied to the map
\[
(z,a,b)\longmapsto z+\gamma^\top(a+b)
\]
completes the proof.
\end{proof}

\begin{proof}[Proof of Corollary~\ref{cor:estcoverage}]
Ties among the scores have probability zero because $F_P$ is continuous, so, almost surely,
\[
\check r_m(\beta):=\Fchk\bigl(\tauchk(\beta)\bigr)-\beta\in\Bigl[0,\ \frac{\max_i\hat v(X_i)}{\sum_j\hat v(X_j)}\Bigr],
\]
and, as in the proof of Theorem~\ref{thm:coverage},
\[
\sqrt m\bigl(P_e(\tauchk(\beta))-(1-\beta)\bigr)=\sqrt m\bigl(\Fchk-F_Q\bigr)\bigl(\tauchk(\beta)\bigr)-\sqrt m\,\check r_m(\beta).
\]
The remainder is negligible by the argument of Lemma~\ref{lem:maxw} with the envelope $V$ in $L_2(P)$ and $\Pn\hat v\toP1$. By Theorem~\ref{thm:parametric} and its proof, uniformly in $t$,
\[
\sqrt m\,(\Fchk-F_Q)(t)=\Zm(t)+\sqrt m\,(\hat\theta-\theta_0)^\top\gamma(t)+o_p(1).
\]
Theorem~\ref{thm:parametric} also gives $\|\Fchk-F_Q\|_\infty\toP0$, which together with continuity and strict monotonicity of $F_Q$ near $q$ yields $\tauchk(\beta)\toP q$ by the standard argument:
\[
\Fchk\bigl(\tauchk(\beta)\bigr)\ge\beta>\Fchk\bigl(\tauchk(\beta)-\bigr),\qquad F_Q(q+\epsilon)>\beta>F_Q(q-\epsilon).
\]
Hence $\gamma(\tauchk(\beta))\toP\gamma(q)$ by continuity, and
\[
\Zm\bigl(\tauchk(\beta)\bigr)-\Zm(q)\toP0
\]
by Lemma~\ref{lem:rho} and \eqref{eq:equicont}, since $F_Q(\tauchk(\beta))\toP F_Q(q)$. Therefore
\[
\sqrt m\bigl(P_e(\tauchk(\beta))-(1-\beta)\bigr)=\Zm(q)+\sqrt m\,(\hat\theta-\theta_0)^\top\gamma(q)+o_p(1)\weak\GQ(q)+\gamma(q)^\top\bigl(\mathbb G\psi+\mathcal N\bigr),
\]
a centered normal variable with variance
\[
C(q,q)+2\,\gamma(q)^\top\Cov\bigl(\GQ(q),\mathbb G\psi\bigr)+\gamma(q)^\top\bigl(P\psi\psi^\top+\Sigma_{\rm aux}\bigr)\gamma(q),
\]
and, because $F_Q(q)=\beta$,
\[
\Cov\bigl(\GQ(q),\mathbb G\psi\bigr)=\E_Q\bigl[\{F_{S|X}(q)-\beta\}\psi(X)\bigr]=c(q).
\]
\end{proof}

\begin{proof}[Verification of Example~\ref{ex:meanshift}]
Differentiating the tilt at $\theta_0$,
\begin{align*}
\dot v(x)&=v(x)\,(x-\theta_0),\\
\gamma(t)&=\E_P\bigl[v(X)(X-\theta_0)\{F_{S|X}(t)-F_Q(t)\}\bigr]=\E_Q\bigl[(X-\theta_0)\{F_{S|X}(t)-F_Q(t)\}\bigr]\\
&=\Cov_Q\bigl(X,F_{S|X}(t)\bigr),
\end{align*}
since $\E_QX=\theta_0$; $L_1(P)$-differentiability holds because $P_X$ is Gaussian, and Lemma~\ref{lem:parametric} applies on any bounded parameter set, which the estimate enters with probability tending to one. (a) The estimation error
\[
\sqrt m\,(\hat\theta-\theta_0)=\sqrt{\frac m{N_t}}\,\sqrt{N_t}\,(\bar X^{\rm tgt}-\theta_0)-\sqrt{\frac m{N_s}}\,\sqrt{N_s}\,\bar X^{\rm src}
\]
is independent of the calibration sample and converges to $N(0,\lambda_t+\lambda_s)$. (b) Here
\[
\sqrt m\,(\hat\theta-\theta_0)=\xi_m+\Gm\psi,\qquad\xi_m=\sqrt m\,(\bar X^{\rm tgt}-\theta_0)\weak N(0,\lambda_t),
\]
with influence function $\psi(x)=-x$, so that $P\psi=0$ and $P\psi^2=1$,
and
\[
c(q)=-\E_Q\bigl[\{F_{S|X}(q)-\beta\}X\bigr]=-\gamma(q).
\]
Corollary~\ref{cor:estcoverage} gives
\[
\tau^2=\sigma^2-2\gamma^2+\gamma^2(1+\lambda_t)=\sigma^2+\gamma^2(\lambda_t-1).
\]
\end{proof}

\begin{proof}[Proof of Theorem~\ref{thm:estboot}]
Throughout, $\Prob^*$, $\E^*$ and $o_{p^*}(1)$ refer to the bootstrap distribution conditionally on the data, with convergence in (outer) probability over the data; $\weak_*$ is conditional weak convergence in probability in the sense of \citet[Section~2.9]{vandervaart1996weak}. Write $\Hcal_\Vcal$ for the Donsker class of Assumption~\ref{ass:donsker}, which contains the classes $\{v_\theta\ind\{S\le t\}\}$ and $\{v_\theta\}$, and note that its envelope $V$ is square integrable.

\emph{Step 1: the coverage process with estimated weights.} By Theorem~\ref{thm:parametric}, $\sqrt m(\Fchk-F_Q)\weak\GQ+\gamma^\top(\mathbb G\psi+\mathcal N)$ in $\linf(\Rbar)$. The weighted empirical distribution function $\Fchk$ has jumps at the scores $S_i$, which are almost surely distinct because $F_P$ is continuous, of sizes $\hat v_i/\sum_j\hat v_j$. Exactly as in the proof of Theorem~\ref{thm:coverage},
\begin{gather*}
\sqrt m\bigl(\beta-F_Q(\tauchk(\beta))\bigr)=\sqrt m\bigl(\Fchk-F_Q\bigr)(\tauchk(\beta))-\sqrt m\,\check r_m(\beta),\\
0\le\check r_m(\beta)\le\frac{\max_i\hat v_i}{\sum_j\hat v_j}\le\frac{\max_iV(X_i)}{\sum_j\hat v_j},
\end{gather*}
and $\sqrt m\max_iV(X_i)/\sum_j\hat v_j\toP0$ because $\E_PV^2<\infty$ (Lemma~\ref{lem:maxw}) and $\Pn\hat v\toP1$. Moreover $|F_Q(\tauchk(\beta))-\beta|\le\|\Fchk-F_Q\|_\infty+\max_i\hat v_i/\sum_j\hat v_j\toP0$ uniformly in $\beta$. We claim that $\gamma$ is a continuous function of $F_Q$: $\gamma=\tilde\gamma\circ F_Q$ with $\tilde\gamma(u):=\gamma(q_Q(u))$ continuous on $(0,1)$. Indeed, on a flat piece $[c,d]$ of $F_Q$ we have $0=Q(c<S\le d)=\E_P[v\,\Prob(c<S\le d\mid X)]$, so $\Prob(c<S\le d\mid X)=0$ almost surely on $\{v>0\}$; and on $\{v=0\}$, since $\theta_0$ is interior and $\theta\mapsto v_\theta(x)\ge0$ is differentiable there with $v_{\theta_0}(x)=0$, the derivative $\dot v_{\theta_0}(x)$ vanishes; hence $\gamma(d)-\gamma(c)=\E_P[\dot v\,\Prob(c<S\le d\mid X)]=0$, and the continuity of $\gamma$ on $\R$ gives the claim. Now, by the proofs of Theorems~\ref{thm:estimated} and~\ref{thm:parametric}, uniformly in $t$,
\[
\sqrt m(\Fchk-F_Q)(t)=\frac{\Gm\tilde k_t}{\Pn\hat v}+\gamma(t)^\top\sqrt m(\hat\theta-\theta_0)+o_p(1),
\]
where the first term is asymptotically $\rho$-equicontinuous and $\rho$-closeness follows from $F_Q$-closeness by Lemma~\ref{lem:rho}, while for the second, with $\tilde\gamma$ uniformly continuous on $[\beta_0/2,(1+\beta_1)/2]$,
\[
\sup_{\beta\in[\beta_0,\beta_1']}\bigl\|\gamma(\tauchk(\beta))-\gamma(q_Q(\beta))\bigr\|=\sup_\beta\bigl\|\tilde\gamma(F_Q(\tauchk(\beta)))-\tilde\gamma(\beta)\bigr\|\toP0
\]
because $\sup_\beta|F_Q(\tauchk(\beta))-\beta|\toP0$. Hence
\[
\sup_{\beta\in[\beta_0,\beta_1]}\bigl|\sqrt m(\Fchk-F_Q)(\tauchk(\beta))-\sqrt m(\Fchk-F_Q)(q_Q(\beta))\bigr|\toP0,
\]
and the first claim follows by composition with the fixed map $q_Q$.

\emph{Step 2: exact centering of the bootstrap statistic.} For $\theta\in\Theta$ let
\[
\hat F_{m,\theta}(t):=\frac{\Pn[v_\theta\ind\{S\le t\}]}{\Pn v_\theta},\qquad
\tilde h_{t,\theta}(x,s):=v_\theta(x)\bigl(\ind\{s\le t\}-\hat F_{m,\theta}(t)\bigr),
\]
so that $\Pn\tilde h_{t,\theta}=0$ for every $\theta$ and $t$, and $\hat F_{m,\hat\theta}=\Fchk$. Since the bootstrap sample is drawn from $\Pn$, the identity of Lemma~\ref{lem:identity} holds with $\Pn$ as the population:
\[
\sqrt m\bigl(\Fchk^*(t)-\hat F_{m,\hat\theta^*}(t)\bigr)=\frac{\sqrt m\,\Pn^*\tilde h_{t,\hat\theta^*}}{\Pn^*v_{\hat\theta^*}}=\frac{\mathbb G^*_m\tilde h_{t,\hat\theta^*}}{\Pn^*v_{\hat\theta^*}} .
\]
Therefore
\begin{equation}\label{eq:bootdecomp}
\sqrt m\bigl(\Fchk^*(t)-\Fchk(t)\bigr)=\frac{\mathbb G^*_m\tilde h_{t,\hat\theta^*}}{\Pn^*v_{\hat\theta^*}}+\sqrt m\bigl(\hat F_{m,\hat\theta^*}(t)-\hat F_{m,\hat\theta}(t)\bigr).
\end{equation}

\emph{Step 3: the bootstrap empirical process term.} Work on the event $\{\hat\theta,\hat\theta^*\in U\}$, whose conditional probability tends to one in probability because $\theta_0$ is interior, $\hat\theta\toP\theta_0$ and $\hat\theta^*-\hat\theta=O_{p^*}(m^{-1/2})$ by \ref{ass:boot-est}. By Theorem~3.6.1 of \citet{vandervaart1996weak}, since $\Hcal_\Vcal$ is Donsker, $\mathbb G^*_m\weak_*\mathbb G$ in $\linf(\Hcal_\Vcal)$ conditionally in probability; as in the proof of Theorem~\ref{thm:bootstrap} this implies unconditional weak convergence of $\mathbb G^*_m$ in $\linf(\Hcal_\Vcal)$, hence asymptotic equicontinuity of $\mathbb G^*_m$ with respect to the $L_2(P)$ seminorm, a statement about the supremum of $|\mathbb G^*_m(f-g)|$ over all pairs in the class at small distance, which therefore applies to pairs that depend on the data and on the bootstrap draws. Now
\[
\mathbb G^*_m\tilde h_{t,\hat\theta^*}=\mathbb G^*_m\bigl(v_{\hat\theta^*}\ind\{S\le t\}\bigr)-\hat F_{m,\hat\theta^*}(t)\,\mathbb G^*_mv_{\hat\theta^*},
\]
and, by the Lipschitz condition of Lemma~\ref{lem:parametric}, $\|v_{\hat\theta^*}\ind\{S\le t\}-v\ind\{S\le t\}\|_{P,2}\le\|\dot V\|_{P,2}\|\hat\theta^*-\theta_0\|$ uniformly in $t$, which tends to zero in probability because $\hat\theta^*-\theta_0=(\hat\theta^*-\hat\theta)+(\hat\theta-\theta_0)=O_{p^*}(m^{-1/2})+o_p(1)$ by \ref{ass:boot-est} and \eqref{eq:thetaexp}. The Glivenko--Cantelli property of $\Hcal_\Vcal$ and the same Lipschitz bound give $\sup_t|\hat F_{m,\hat\theta^*}(t)-F_Q(t)|=o_{p^*}(1)$ and $\Pn^*v_{\hat\theta^*}=1+o_{p^*}(1)$ (recall $Pv=1$). Hence, uniformly in $t$,
\[
\frac{\mathbb G^*_m\tilde h_{t,\hat\theta^*}}{\Pn^*v_{\hat\theta^*}}=\mathbb G^*_m\bigl(v\ind\{S\le t\}\bigr)-F_Q(t)\,\mathbb G^*_mv+o_{p^*}(1)=\mathbb G^*_mh_t+o_{p^*}(1).
\]

\emph{Step 4: the tilt term.} Since $\Pn\tilde h_{t,\hat\theta}=0$,
\[
\hat F_{m,\hat\theta^*}(t)-\hat F_{m,\hat\theta}(t)=\frac{\Pn\bigl[(v_{\hat\theta^*}-v_{\hat\theta})(\ind\{S\le t\}-\hat F_{m,\hat\theta}(t))\bigr]}{\Pn v_{\hat\theta^*}} .
\]
By \ref{ass:boot-smooth}, $v_{\hat\theta^*}-v_{\hat\theta}=(\hat\theta^*-\hat\theta)^\top\dot v_{\hat\theta}+R$ with $\Pn|R|\le\Pn\ddot H\,\|\hat\theta^*-\hat\theta\|^2=O_p(1)\,O_{p^*}(m^{-1})$. Hence, uniformly in $t$,
\[
\sqrt m\bigl(\hat F_{m,\hat\theta^*}(t)-\hat F_{m,\hat\theta}(t)\bigr)=\sqrt m(\hat\theta^*-\hat\theta)^\top\hat\gamma_m(t)+o_{p^*}(1),\qquad
\hat\gamma_m(t):=\frac{\Pn\bigl[\dot v_{\hat\theta}(\ind\{S\le t\}-\hat F_{m,\hat\theta}(t))\bigr]}{\Pn v_{\hat\theta}} ,
\]
where $\Pn v_{\hat\theta^*}-\Pn v_{\hat\theta}=O_{p^*}(m^{-1/2})$ has been absorbed into the remainder. The Glivenko--Cantelli property of $\{\dot v_\theta\ind\{S\le t\}\}$ in \ref{ass:boot-smooth}, the $L_1(P)$-continuity of $\theta\mapsto\dot v_\theta$ at $\theta_0$, $\hat\theta\toP\theta_0$, $\sup_t|\hat F_{m,\hat\theta}-F_Q|\toP0$ (Theorem~\ref{thm:parametric}) and $\Pn v_{\hat\theta}\toP1$ give $\sup_t\|\hat\gamma_m(t)-\gamma(t)\|\toP0$, with $\gamma$ as in \eqref{eq:gamma} (recall $Pv=1$). With \ref{ass:boot-est},
\[
\sqrt m\bigl(\hat F_{m,\hat\theta^*}(t)-\hat F_{m,\hat\theta}(t)\bigr)=\gamma(t)^\top\bigl(\mathbb G^*_m\psi+\xi^*_m\bigr)+o_{p^*}(1)\qquad\text{uniformly in }t .
\]

\emph{Step 5: conclusion.} Combining \eqref{eq:bootdecomp} with Steps 3 and 4,
\[
\sqrt m\bigl(\Fchk^*-\Fchk\bigr)(t)=\mathbb G^*_mh_t+\gamma(t)^\top\bigl(\mathbb G^*_m\psi+\xi^*_m\bigr)+o_{p^*}(1)\qquad\text{uniformly in }t .
\]
The process $(\mathbb G^*_mh_{(\cdot)},\mathbb G^*_m\psi)$ converges conditionally to $(\GQ,\mathbb G\psi)$, jointly, by Theorem~3.6.1 of \citet{vandervaart1996weak} applied to the Donsker class $\Hcal_\Vcal\cup\{\psi_1,\dots,\psi_d\}$; and $\xi^*_m$, which is computed from $\mathcal A^*$ alone, converges conditionally to $\mathcal N\sim N(0,\Sigma_{\rm aux})$ independently of the calibration resampling. Hence $\sqrt m(\Fchk^*-\Fchk)\weak_*\GQ+\gamma^\top(\mathbb G\psi+\mathcal N)$ in $\linf(\Rbar)$. The composition with $\tauchk(\beta)$ is handled as in Step~1: the first component $\mathbb G^*_mh_t$ is asymptotically $\rho$-equicontinuous and Lemma~\ref{lem:rho} converts $F_Q$-closeness into $\rho$-closeness, the second is $\gamma(t)^\top O_{p^*}(1)$ with $\gamma=\tilde\gamma\circ F_Q$, and $\sup_\beta|F_Q(\tauchk(\beta))-\beta|\toP0$; so $\Lambda^*_m(\beta)=\sqrt m(\Fchk^*-\Fchk)(q_Q(\beta))+o_{p^*}(1)$ uniformly in $\beta$, which gives $\Lambda^*_m\weak_*\Lb^{\rm est}$.

\emph{Inflation.} The variance $\tau^2(\beta)$ of $\Lb^{\rm est}(\beta)$ is continuous in $\beta$ (each of $\sigma^2$, $\tilde\gamma$ and $\beta\mapsto P[h_{q_Q(\beta)}\psi]$ is, the last by $|P[(h_s-h_t)\psi]|\le\rho(s,t)\|\psi\|_{P,2}$ and Lemma~\ref{lem:rho}), so under the hypothesis $\tau^2>0$ on $\Bcal$ it is bounded below there; the sufficient condition in the statement follows from $\tau^2(\beta)=P(h_q+\gamma(q)^\top\psi)^2+\gamma(q)^\top\Sigma_{\rm aux}\gamma(q)\ge\E_P[\Var(h_q+\gamma(q)^\top\psi\mid X)]=\E_P[v^2F_{S|X}(q)\{1-F_{S|X}(q)\}]$, because $\psi$ depends on $X$ only. The supremum over a closed $\Bcal$ is a continuous functional, so $\hat c^*\toP c^{\rm est}_{1-\delta}$, the $(1-\delta)$-quantile of $\sup_\Bcal\Lb^{\rm est}$, which is a continuity point of its distribution because the supremum of a centered Gaussian process with variances bounded away from zero has an atomless law (as in the proof of Corollary~\ref{cor:simultaneous}). The failure event $\{\exists\beta\in\Bcal:P_e(\tauchk(\beta+\hat c^*/\sqrt m))>1-\beta\}$ equals $\{\sup_{\beta\in\Bcal}\sqrt m(\beta'-F_Q(\tauchk(\beta')))>\hat c^*\}$ with $\beta'=\beta+\hat c^*/\sqrt m$, and the argument of the proof of Corollary~\ref{cor:simultaneous}, with the first claim of the theorem (on $[\beta_0,\beta_1']$) in place of Theorem~\ref{thm:coverage}, gives probability tending to $\delta$. Finally, the multiplier bootstrap with $\hat v$ held fixed is $\hat{\mathbb Z}^\xi_m$ of \eqref{eq:multiplier} with $\hat v$ in place of $v$. Writing $\tilde k_t=\hat v(\ind\{S\le t\}-F_{\hat Q}(t))$ as in the proof of Theorem~\ref{thm:estimated}, its summands are $\xi_i\{\tilde k_t(Z_i)-\hat v_i(\Fchk(t)-F_{\hat Q}(t))\}$, and the second part contributes $\|\Fchk-F_{\hat Q}\|_\infty\,|m^{-1/2}\sum_i\xi_i\hat v_i|=O_p(m^{-1/2})O_p(1)$ by Theorem~\ref{thm:estimated}; the first part satisfies $\sup_t|\Gm^\xi\tilde k_t-\Gm^\xi h_t|\toP0$ by the equicontinuity of the multiplier process over $\Hcal_\Vcal$ and $\sup_t\|\tilde k_t-h_t\|_{P,2}\toP0$. Hence the fixed-weight bootstrap converges conditionally to $\GQ$, and, by Step~1 and Lemma~\ref{lem:rho}, to $\Lb$ after composition with $\tauchk$; its conditional $(1-\delta)$-quantile at a single level converges to $z_{1-\delta}\sigma(\beta)$, so the failure probability of the set run at $\beta+\hat c/\sqrt m$ converges to $\Prob(\Lb^{\rm est}(\beta)>z_{1-\delta}\sigma(\beta))=1-\Phi(z_{1-\delta}\sigma(\beta)/\tau(\beta))$, which exceeds $\delta$ if and only if $\tau^2(\beta)>\sigma^2(\beta)$. With independent auxiliary data $\psi\equiv0$ and $\tau^2-\sigma^2=\gamma(q)^\top\Sigma_{\rm aux}\gamma(q)$; in case~(b) of Example~\ref{ex:meanshift}, $\tau^2-\sigma^2=\gamma(q)^2(\lambda_t-1)<0$ when $\lambda_t<1$ and $\gamma(q)\ne0$.
\end{proof}

\begin{proof}[Verification of Assumption~\ref{ass:boot} for Example~\ref{ex:meanshift}]
Take $v_\theta(x)=\exp(\theta x-\theta^2/2)$ on a bounded open neighbourhood $\Theta$ of $\theta_0$, so that $v_{\theta_0}=v$. Then $\dot v_\theta(x)=(x-\theta)v_\theta(x)$, $|v_\theta-v_{\theta'}-(\theta-\theta')\dot v_{\theta'}|\le\ddot H\,|\theta-\theta'|^2$ with $\ddot H(x)=C\,(1+x^2)\sup_{\theta\in\Theta}e^{\theta x}\in L_1(N(0,1))$, $\theta\mapsto\dot v_\theta$ is $L_1$-continuous by dominated convergence, and $\{\dot v_\theta\ind\{S\le t\}\}$ is Glivenko--Cantelli as the product of the Glivenko--Cantelli classes $\{\dot v_\theta\}$ (Lipschitz in $\theta$ with integrable Lipschitz constant, \citealp[Example~19.7]{vandervaart1998asymptotic}) and $\{\ind\{S\le t\}\}$, both with integrable envelopes \citep[Corollary~9.27]{kosorok2008introduction}. For the estimator, with $N_s=m$ (case (b)),
\begin{gather*}
\sqrt m(\hat\theta^*-\hat\theta)=-\sqrt m\bigl(\bar X^*-\bar X\bigr)+\sqrt m\bigl(\bar X'^*-\bar X'\bigr)=\mathbb G^*_m\psi+\xi^*_m,\\
\psi(x)=-x,\qquad\xi^*_m=\sqrt{m/N_t}\,\sqrt{N_t}\bigl(\bar X'^*-\bar X'\bigr),
\end{gather*}
exactly, and $\xi^*_m\weak_*N(0,\lambda_t)$ by the bootstrap central limit theorem for the mean of the target covariates; in case (a) the source mean is resampled as well, $\psi\equiv0$ and $\xi^*_m\weak_*N(0,\lambda_t+\lambda_s)$.
\end{proof}

\section{Proofs for Section~\ref{sec:stable}}\label{app:stable}

We use the classical criterion for convergence of row sums of a triangular array of i.i.d.\ summands to an infinitely divisible law without Gaussian part \citep[Chapter~4, \S25]{gnedenko1954limit}; see also \citet[Chapter~15]{kallenberg2002foundations} and \citet[Chapter~7]{resnick2007heavy}. Let $\xi_{m1},\dots,\xi_{mm}$ be i.i.d.\ random vectors in $\R^k$ with
\[
\E\xi_{m1}=0,\qquad\E\|\xi_{m1}\|<\infty .
\]
If
\begin{enumerate}[label=(C\arabic*),leftmargin=*,itemsep=1pt]
\item the measures $m\,\Prob(\xi_{m1}\in\cdot)$ converge vaguely on $\R^k\setminus\{0\}$ to a L\'evy measure $\nu$ with
\[
\int\bigl(\|y\|^2\wedge\|y\|\bigr)\,\nu(\dd y)<\infty,
\]
\item the small values carry no variance in the limit,
\[
\lim_{\epsilon\downarrow0}\limsup_{m}\,m\,\E\bigl[\|\xi_{m1}\|^2\ind\{\|\xi_{m1}\|\le\epsilon\}\bigr]=0,
\]
\item the large values are uniformly integrable,
\[
\lim_{M\to\infty}\limsup_m\,m\,\E\bigl[\|\xi_{m1}\|\ind\{\|\xi_{m1}\|>M\}\bigr]=0,
\]
\end{enumerate}
then the row sums converge in distribution to the random vector $\Scal$ with
\[
\E\,e^{i\langle u,\Scal\rangle}=\exp\Bigl\{\int\bigl(e^{i\langle u,y\rangle}-1-i\langle u,y\rangle\bigr)\,\nu(\dd y)\Bigr\}.
\]
Indeed, (C1) and (C2) are the L\'evy-measure and no-Gaussian-part conditions of the criterion with truncation at a level $\tau$, and (C3) together with the centering of the summands identifies the drift: for every $\tau$ with $\nu(\{\|y\|=\tau\})=0$,
\[
m\,\E\bigl[\xi_{m1}\ind\{\|\xi_{m1}\|\le\tau\}\bigr]=-m\,\E\bigl[\xi_{m1}\ind\{\|\xi_{m1}\|>\tau\}\bigr]\longrightarrow-\int_{\|y\|>\tau}y\,\nu(\dd y),
\]
by (C1) and the uniform integrability (C3).

Two consequences of Assumption~\ref{ass:heavy}(i) are used repeatedly \citep[Chapter~2]{resnick2007heavy}. By the uniform convergence theorem for regularly varying functions,
\[
m\,\Prob(v>xa_m)\longrightarrow x^{-p}\qquad\text{locally uniformly in }x>0,
\]
and by Karamata's theorem,
\begin{equation}\label{eq:karamata}
m\,a_m^{-2}\,\E\bigl[v^2\ind\{v\le xa_m\}\bigr]\longrightarrow\frac{p}{2-p}\,x^{2-p},\qquad
m\,a_m^{-1}\,\E\bigl[v\ind\{v>xa_m\}\bigr]\longrightarrow\frac{p}{p-1}\,x^{1-p}.
\end{equation}

\begin{proof}[Proof of Theorem~\ref{thm:stablefidi}]
Let
\[
\mathbf c(s)=\bigl(\mathbf c_{(-\infty,t_j]}(s)\bigr)_{j\le k},\qquad Y_i=v_i\,\mathbf c(S_i),\qquad\xi_{mi}=a_m^{-1}Y_i .
\]
Then
\[
\E Y_i=\E_P\bigl[v\,\mathbf c(S)\bigr]=\E_Q\mathbf c(S)=0,
\]
and, by Lemma~\ref{lem:identity} applied coordinatewise,
\[
\frac{m}{a_m}\Bigl(\Fhat(t_j)-F_Q(t_j)\Bigr)_{j\le k}=\frac{1}{\Pn v}\sum_{i=1}^m\xi_{mi}.
\]
Since $\Pn v\to1$ almost surely, it suffices to verify (C1)--(C3). Note that
\[
\underline c:=\min_j\min\{F_Q(t_j),1-F_Q(t_j)\}\le\|\mathbf c(s)\|\le\sqrt k
\]
(if some coordinate has $F_Q(t_j)$ equal to zero or one it vanishes $Q$-almost surely and can be dropped).

(C1) Let $A$ be a bounded Borel subset of $\R^k\setminus\{0\}$ at positive distance $\epsilon$ from the origin with $\nu(\partial A)=0$. Then $\xi_{m1}\in A$ forces $v>\epsilon a_m/\sqrt k$, and
\[
m\,\Prob(\xi_{m1}\in A)=m\,\Prob\bigl(v>\epsilon a_m/\sqrt k\bigr)\cdot\Prob\bigl(a_m^{-1}v\,\mathbf c(S)\in A\ \big|\ v>\epsilon a_m/\sqrt k\bigr).
\]
The first factor tends to $(\epsilon/\sqrt k)^{-p}$. By Assumption~\ref{ass:heavy}(ii), the conditional law of the pair
\[
\Bigl(\frac{v}{\epsilon a_m/\sqrt k},\ S\Bigr)\quad\text{given}\quad v>\epsilon a_m/\sqrt k
\]
converges weakly to the product of a Pareto law with index $p$ and $\Pi$ (the sets $(r,\infty)\times A$ form a convergence-determining class, the first marginal being continuous); the map $(r,s)\mapsto r\epsilon k^{-1/2}\mathbf c(s)$ is continuous except on the $\Pi$-null set $\{t_1,\dots,t_k\}$, so the continuous mapping theorem gives convergence of the second factor to the pushforward probability. Multiplying out,
\[
m\,\Prob(\xi_{m1}\in A)\longrightarrow\int_0^\infty\!\!\int\ind\{r\mathbf c(s)\in A\}\,\Pi(\dd s)\,p\,r^{-p-1}\dd r=\nu(A),
\]
where we used that the Pareto law rescaled by $\epsilon/\sqrt k$ and multiplied by $(\epsilon/\sqrt k)^{-p}$ is the measure $pr^{-p-1}\dd r$ on $(\epsilon/\sqrt k,\infty)$. This is vague convergence to $\nu$, and $\nu$ satisfies the integrability condition because
\[
\|r\,\mathbf c(s)\|\le\sqrt k\,r\qquad\text{and}\qquad
\int_0^\infty(r^2\wedge r)\,p\,r^{-p-1}\dd r<\infty\qquad\text{for }1<p<2 .
\]

(C2) The event $\|\xi_{m1}\|\le\epsilon$ implies $v\le\epsilon a_m/\underline c$, so by \eqref{eq:karamata}
\[
m\,\E\bigl[\|\xi_{m1}\|^2\ind\{\|\xi_{m1}\|\le\epsilon\}\bigr]\le k\,m\,a_m^{-2}\,\E\bigl[v^2\ind\{v\le\epsilon a_m/\underline c\}\bigr]\longrightarrow k\,\frac{p}{2-p}\Bigl(\frac{\epsilon}{\underline c}\Bigr)^{2-p},
\]
which tends to zero as $\epsilon$ decreases to zero.

(C3) The event $\|\xi_{m1}\|>M$ implies $v>Ma_m/\sqrt k$, so by \eqref{eq:karamata}
\[
m\,\E\bigl[\|\xi_{m1}\|\ind\{\|\xi_{m1}\|>M\}\bigr]\le\sqrt k\,m\,a_m^{-1}\,\E\bigl[v\ind\{v>Ma_m/\sqrt k\}\bigr]\longrightarrow k^{p/2}\,\frac{p}{p-1}\,M^{1-p},
\]
which tends to zero as $M$ tends to infinity.
\end{proof}

\begin{proof}[Proof of Theorem~\ref{thm:stablecoverage}]
Let $W_i=F_Q(S_i)$. Under $Q$ the variable $F_Q(S)$ is uniform on the unit interval because $F_Q$ is continuous, and
\[
\hat W(\beta):=F_Q\bigl(\tauhat(\beta)\bigr)
\]
is the self-normalized weighted $\beta$-quantile of the $W_i$: since $F_Q$ is nondecreasing and continuous, the weighted empirical distribution function $\Fhat^W$ of the $W_i$ satisfies
\[
\Fhat^W\bigl(F_Q(t)\bigr)=\Fhat(t)\qquad\text{almost surely},
\]
and the infimum defining $\tauhat$ maps to the infimum defining the weighted quantile of the $W_i$. Fix $z\in\R$ and put
\[
w_m:=\beta-\frac{za_m}{m}\longrightarrow\beta .
\]
For the right-continuous step function $\Fhat^W$ the duality
\[
\{\hat W(\beta)\le w_m\}=\{\Fhat^W(w_m)\ge\beta\}
\]
is exact, and on the event that not all weights vanish,
\[
\Fhat^W(w_m)\ge\beta\iff\sum_{i=1}^mv_i\bigl(\ind\{W_i\le w_m\}-w_m\bigr)\ge(\beta-w_m)\sum_{i=1}^mv_i\iff\sum_{i=1}^m\xi_{mi}\ge z\,\Pn v,
\]
where
\[
\xi_{mi}:=a_m^{-1}v_i\bigl(\ind\{W_i\le w_m\}-w_m\bigr).
\]
Since
\[
\frac{m}{a_m}\bigl\{P_e(\tauhat(\beta))-(1-\beta)\bigr\}=\frac{m}{a_m}\bigl\{\beta-\hat W(\beta)\bigr\},
\]
the event that this quantity is at least $z$ is exactly
\[
\{\hat W(\beta)\le w_m\}=\Bigl\{\sum_{i=1}^m\xi_{mi}\ge z\,\Pn v\Bigr\}.
\]
We verify (C1)--(C3) for the one-dimensional array $\xi_{mi}$, which is centered:
\[
\E\xi_{m1}=a_m^{-1}\bigl\{Q(W\le w_m)-w_m\bigr\}=0 .
\]

(C1) For $y>0$, the event $\xi_{m1}>y$ is the event that $W\le w_m$ and $v>t_m$, where
\[
t_m:=\frac{a_my}{1-w_m},\qquad\frac{t_m}{a_m}\longrightarrow\frac{y}{1-\beta}.
\]
Thus
\[
m\,\Prob(\xi_{m1}>y)=m\,\Prob(v>t_m)\cdot\Prob(W\le w_m\mid v>t_m),
\]
and the first factor tends to $(y/(1-\beta))^{-p}$. For the second, assume $z>0$, so that $w_m$ increases to $\beta$ (the case $z\le0$ is symmetric); for every $\epsilon>0$ and large $m$,
\[
\Prob(W\le\beta-\epsilon\mid v>t_m)\le\Prob(W\le w_m\mid v>t_m)\le\Prob(W\le\beta\mid v>t_m),
\]
and by Assumption~\ref{ass:heavy}(ii) with $r=1$ the two bounds converge to
\[
\Pi\bigl(\{s:F_Q(s)\le\beta-\epsilon\}\bigr)\qquad\text{and}\qquad\Pi\bigl(\{s:F_Q(s)\le\beta\}\bigr)
\]
respectively, provided these are $\Pi$-continuity sets, which holds for all but countably many $\epsilon$ and, for $\epsilon=0$, by the assumption that $\Pi$ does not charge the set where $F_Q$ equals $\beta$. Letting $\epsilon$ decrease to zero along admissible values,
\[
\Prob(W\le w_m\mid v>t_m)\longrightarrow\Pi\bigl(\{F_Q\le\beta\}\bigr)=\Pi\bigl((-\infty,q]\bigr)=\pi_\beta,
\]
where $\{F_Q\le\beta\}$ is the union of $(-\infty,q]$ and the $\Pi$-null set $\{F_Q=\beta\}$. Hence
\[
m\,\Prob(\xi_{m1}>y)\longrightarrow\pi_\beta(1-\beta)^p\,y^{-p}.
\]
Similarly, for $y<0$, the event $\xi_{m1}<y$ is the event that $W>w_m$ and $v>a_m|y|/w_m$, and
\[
m\,\Prob(\xi_{m1}<y)\longrightarrow(1-\pi_\beta)\beta^p\,|y|^{-p}.
\]
So the L\'evy measure $\nu_\beta$ has density
\[
p\,\pi_\beta(1-\beta)^p\,y^{-p-1}\ \text{ on }(0,\infty),\qquad p\,(1-\pi_\beta)\beta^p\,|y|^{-p-1}\ \text{ on }(-\infty,0).
\]

(C2)--(C3) As in the previous proof,
\[
a_m^{-1}v\,\min\{w_m,1-w_m\}\le|\xi_{m1}|\le a_m^{-1}v,
\]
and \eqref{eq:karamata} gives both conditions.

Hence the row sums converge in distribution to $\Scal_\beta$, the mean-zero infinitely divisible law with L\'evy measure $\nu_\beta$ and no Gaussian part, which is $p$-stable. Since $\Pn v\to1$ almost surely and stable laws are absolutely continuous,
\[
\Prob\Bigl(\sum_i\xi_{mi}\ge z\,\Pn v\Bigr)\longrightarrow\Prob(\Scal_\beta\ge z)\qquad\text{for every }z,
\]
which is the claimed convergence in distribution. For the parameters: for $1<p<2$ and $u>0$ \citep[Lemma~14.11]{sato1999levy},
\[
\int_0^\infty\bigl(e^{iuy}-1-iuy\bigr)y^{-p-1}\dd y=\Gamma(-p)\,e^{-i\pi p/2}\,u^p ,
\]
so with
\[
c_+=p\,\pi_\beta(1-\beta)^p,\qquad c_-=p\,(1-\pi_\beta)\beta^p
\]
the log-characteristic function is
\[
u^p\,\Gamma(-p)\Bigl[(c_++c_-)\cos\tfrac{\pi p}{2}-i\,(c_+-c_-)\sin\tfrac{\pi p}{2}\Bigr]
=-\varsigma^p u^p\Bigl(1-i\kappa\tan\tfrac{\pi p}{2}\Bigr),
\]
with
\[
\varsigma^p=-\Gamma(-p)\cos\tfrac{\pi p}{2}\,(c_++c_-)=\Gamma(1-p)\cos\tfrac{\pi p}{2}\bigl\{\pi_\beta(1-\beta)^p+(1-\pi_\beta)\beta^p\bigr\},
\qquad
\kappa=\frac{c_+-c_-}{c_++c_-},
\]
using $-p\,\Gamma(-p)=\Gamma(1-p)$. This is the law $S_p(\varsigma,\kappa,0)$ in the parametrization of \citet[Definition~1.1.6]{samorodnitsky1994stable}.
\end{proof}

\section{Proofs for Section~\ref{sec:rates}}\label{app:rates}

Throughout this appendix $v\le B$, $F_Q$ is continuous, and $c,c',c''$ denote absolute constants whose value may change from line to line, while $C$ denotes a constant depending only on $(\beta_0,\beta_1)$. We write
\[
U_m(t):=\Gm h_t,\qquad \Zm(t)=\frac{U_m(t)}{\Pn v},\qquad q_j:=q_Q(b_j),\qquad \sigma_0:=\min\{\beta_0,\,1-\beta_1\}.
\]
Since $K^2\ge1$ and $\tilde G(\beta)\in[0,K^2]$, the single-level variance is bounded below on the grid:
\begin{equation}\label{eq:sigmalower}
\sigma^2(\beta)=\beta^2K^2+(1-2\beta)\tilde G(\beta)\ \ge\ \min\{\beta^2K^2,(1-\beta)^2K^2\}\ \ge\ \sigma_0^2\qquad(\beta\in[\beta_0,\beta_1]),
\end{equation}
because $\sigma^2(\beta)$ is affine in $\tilde G(\beta)$ and equals $\beta^2K^2$ at $\tilde G=0$ and $(1-\beta)^2K^2$ at $\tilde G=K^2$. We also have $\sigma^2(\beta)\le K^2\le B$.

\subsection{Proof of Proposition~\ref{prop:berryesseen}}

\begin{proof}
Fix $\beta\in(0,1)$ and write $\sigma=\sigma(\beta)$, $q(w)=q_Q(w)$. Since $F_Q$ is continuous, $F_Q(q(w))=w$ for $w\in(0,1)$, and $F_Q(t)<w$ for $t<q(w)$. Hence
\[
\{F_Q(\tauhat(\beta))<w\}=\{\tauhat(\beta)<q(w)\}=\{\Fhat(q(w)-)\ge\beta\},
\]
the last equality because $\tauhat(\beta)=\inf\{t:\Fhat(t)\ge\beta\}$ and $\Fhat$ is a right-continuous step function, so that $\Fhat(q(w)-)=\sup_{t<q(w)}\Fhat(t)$ is attained. Given $x\in\R$ put
\[
w:=\beta-\frac{x\sigma}{\sqrt m},\qquad\text{so that}\qquad
\Prob\bigl(\sqrt m\,\Cm(\beta)\le x\sigma\bigr)=\Prob\bigl(F_Q(\tauhat(\beta))\ge w\bigr)=1-\Prob\bigl(\Fhat(q(w)-)\ge\beta\bigr)
\]
when $w\in(0,1)$. On the event $\Pn v>0$, whose complement has probability at most $(1-1/B)^m\le e^{-m/B}\le\sqrt{B/(2em)}\le B/(\sigma\sqrt{2em})$ (using $e^{-y}\le(2ey)^{-1/2}$ and $\sigma^2\le B$),
\[
\{\Fhat(q(w)-)\ge\beta\}=\Bigl\{\sum_{i=1}^m\zeta_i\ge0\Bigr\},\qquad \zeta_i:=v_i\bigl(\ind\{S_i<q(w)\}-\beta\bigr),
\]
and the $\zeta_i$ are i.i.d.\ with $|\zeta_i|\le B$ and
\[
\E\zeta_i=Q\bigl(S<q(w)\bigr)-\beta=w-\beta=-\frac{x\sigma}{\sqrt m},\qquad
s_w^2:=\Var\zeta_i=(1-2\beta)\,G\bigl(q(w)\bigr)+\beta^2K^2-(w-\beta)^2 .
\]
(We used $E_P[v\ind\{S<q(w)\}]=Q(S<q(w))=w$ and that $G$ is continuous because $G\ll F_Q$.) At $w=\beta$ the first two terms equal $\sigma^2$, and $|G(q(w))-G(q(\beta))|=\E_Q[v\ind\{S\text{ between }q(w),q(\beta)\}]\le B|w-\beta|$, so
\begin{equation}\label{eq:sw}
|s_w^2-\sigma^2|\le\frac{B|x|\sigma}{\sqrt m}+\frac{x^2\sigma^2}{m}.
\end{equation}
\emph{Central range.} Let $|x|\le\sigma\sqrt m/(4B)$. Since $\sigma^2=\E_P[v^2(\ind\{S\le q\}-\beta)^2]\le B\,\E_Q[(\ind\{S\le q\}-\beta)^2]=B\beta(1-\beta)$, we have $|w-\beta|=|x|\sigma/\sqrt m\le\sigma^2/(4B)\le\beta(1-\beta)/4$, so $w\in(0,1)$. Moreover \eqref{eq:sw} and $\sigma^2\le B\le B^2$ give $|s_w^2-\sigma^2|\le5\sigma^2/16$, so $\sigma^2/2\le s_w^2\le2\sigma^2$. The Berry--Esseen theorem for the standardized sum, with $\E|\zeta-\E\zeta|^3\le2B\,s_w^2$, gives
\[
\Bigl|\Prob\Bigl(\sum_i\zeta_i\ge0\Bigr)-\bar\Phi\Bigl(\frac{x\sigma}{s_w}\Bigr)\Bigr|
=\Bigl|\Prob\Bigl(\frac{\sum_i(\zeta_i-\E\zeta_i)}{s_w\sqrt m}\ge\frac{x\sigma}{s_w}\Bigr)-\bar\Phi\Bigl(\frac{x\sigma}{s_w}\Bigr)\Bigr|
\le\frac{2C_{\rm BE}B}{s_w\sqrt m}\le\frac{2\sqrt2\,C_{\rm BE}B}{\sigma\sqrt m}.
\]
By the mean value theorem, with $\xi$ between $x$ and $x\sigma/s_w$, so that $|\xi|\ge|x|/\sqrt2$,
\[
\Bigl|\Phi\Bigl(\frac{x\sigma}{s_w}\Bigr)-\Phi(x)\Bigr|=\phi(\xi)\,|x|\,\Bigl|\frac{\sigma}{s_w}-1\Bigr|
\le\phi\bigl(|x|/\sqrt2\bigr)\,|x|\,\frac{|\sigma^2-s_w^2|}{s_w(\sigma+s_w)}
\le\sqrt2\,\phi\bigl(|x|/\sqrt2\bigr)\Bigl(\frac{Bx^2}{\sigma\sqrt m}+\frac{|x|^3}{m}\Bigr),
\]
and $\sup_x x^2\phi(x/\sqrt2)+\sup_x|x|^3\phi(x/\sqrt2)\le c$. Since $1/m\le B/(\sigma\sqrt m)$ (as $\sigma\le\sqrt B\le B\sqrt m$), the total error in the central range is at most $cB/(\sigma\sqrt m)$.

\emph{Tails.} Let $x>\sigma\sqrt m/(4B)$; the case $x<-\sigma\sqrt m/(4B)$ is symmetric. Then $\bar\Phi(x)\le e^{-x^2/2}\le e^{-m\sigma^2/(32B^2)}$. If $w\in(0,1)$, Bernstein's inequality for the $\zeta_i$, with $|\zeta_i-\E\zeta_i|\le2B$ and variance $s_w^2\le\sigma^2+Bx\sigma/\sqrt m+x^2\sigma^2/m$ by \eqref{eq:sw}, gives
\begin{align*}
\Prob\bigl(\sqrt m\,\Cm(\beta)>x\sigma\bigr)&=\Prob\Bigl(\sum_i(\zeta_i-\E\zeta_i)\ge x\sigma\sqrt m\Bigr)
\le\exp\Bigl(-\frac{x^2\sigma^2m}{2m\sigma^2+\frac{10}{3}Bx\sigma\sqrt m+2x^2\sigma^2}\Bigr)\\
&\le\exp\Bigl(-\frac13\min\Bigl\{\frac{x^2}{2},\frac{3x\sigma\sqrt m}{10B},\frac m2\Bigr\}\Bigr),
\end{align*}
and on the range $x>\sigma\sqrt m/(4B)$, using $\sigma^2\le B\le B^2$, each of the three terms is at least $\sigma^2m/(32B^2)$, so the probability is at most $e^{-\sigma^2m/(96B^2)}$. If $w\le0$ the left-hand side is $0$. In both cases the error is at most $2e^{-m\sigma^2/(96B^2)}\le cB/(\sigma\sqrt m)$, by $e^{-y}\le(2ey)^{-1/2}$. (For $x<0$ and $w\ge1$ one uses that $F_Q(\tauhat(\beta))<1$ almost surely: $\tauhat(\beta)$ is a score $S_i$ with $v_i>0$, and $\E_P[v\ind\{F_Q(S)=1\}]=Q(F_Q(S)=1)=0$ forces $v\ind\{F_Q(S)=1\}=0$ $P$-a.s., so $\Prob(F_Q(S)=1,\,v>0)=0$.) Collecting the bounds proves the proposition with a universal $C_0$.
\end{proof}

\subsection{Deviation and oscillation bounds for bounded ratios}

\begin{lemma}[deviation bounds]\label{lem:devbounded}
There is an absolute constant $c_1$ such that for every $m\ge2$ and every finite set $q_1,\dots,q_L$, with probability at least $1-4/m$, simultaneously,
\begin{align}
|\Pn v-1|&\le c_1\sqrt{\frac{B\log m}{m}}+c_1\frac{B\log m}{m},\label{eq:dev1}\\
\|\Fhat-F_Q\|_\infty&\le c_1\sqrt{\frac{B\log m}{m}}\qquad\text{provided }m\ge4B^2\log m,\label{eq:dev2}\\
\max_{i\le m}\hat w_i&\le\frac{B}{m\,\Pn v},\qquad \hat w_i:=\frac{v_i}{\sum_jv_j},\label{eq:dev3}\\
\max_{j\le L}|U_m(q_j)|&\le c_1\sqrt{B\log(Lm)}+c_1\frac{B\log(Lm)}{\sqrt m}.\label{eq:dev4}
\end{align}
\end{lemma}

\begin{proof}
\eqref{eq:dev1}: Bernstein's inequality for $v_i-1\in[-1,B]$ with variance $K^2-1\le B$. \eqref{eq:dev2}: Proposition~\ref{prop:bousquet} with $u=\log m$ and $\eta=1/2$ gives, with probability at least $1-1/m-e^{-m/(2B^2)}\ge1-2/m$ when $m\ge4B^2\log m$,
\[
\|\Fhat-F_Q\|_\infty\le2\Bigl[\frac{C_1K}{\sqrt m}+\sqrt{\frac{2(B/4+2C_1B^{3/2}/\sqrt m)\log m}{m}}+\frac{B\log m}{3m}\Bigr]\le c\sqrt{\frac{B\log m}{m}},
\]
using $K\le\sqrt B$, $\sigma_*^2\le B/4$, $B^{3/2}\log m/m^{3/2}\le B\log m/m$ and $B\log m/m\le\sqrt{B\log m/m}$ under $m\ge4B^2\log m$. \eqref{eq:dev3} is $v_i\le B$. \eqref{eq:dev4}: Bernstein's inequality for $h_{q_j}$, bounded by $B$ with variance at most $K^2\le B$, at level $1/(Lm)$, and a union bound.
\end{proof}

\begin{lemma}[modulus of continuity of the weighted score process]\label{lem:osc}
For $\delta\in(0,1]$ let
\begin{align*}
\omega_m(\delta)&:=\sup\bigl\{|U_m(s)-U_m(t)|:\ s<t,\ F_Q(t)-F_Q(s)\le\delta\bigr\},\\
\omega_m'(\delta)&:=\sup\bigl\{|\Gm(v\ind\{s<S\le t\})|:\ s<t,\ F_Q(t)-F_Q(s)\le\delta\bigr\}.
\end{align*}
There is an absolute constant $c_2$ such that for every $m\ge2$, with probability at least $1-2/m$,
\[
\max\{\omega_m(\delta),\omega'_m(\delta)\}\le c_2\Bigl[\sqrt{B\delta\log(8m/\delta)}+\frac{B\log(8m/\delta)}{\sqrt m}\Bigr].
\]
\end{lemma}

\begin{proof}
Consider the class $\Hcal_\delta:=\{h_s-h_t:\ s<t,\ F_Q(t)-F_Q(s)\le\delta\}$. Each $f\in\Hcal_\delta$ has the form $f=v(\ind\{s<S\le t\}-(F_Q(t)-F_Q(s)))$, so that
\[
Pf=0,\qquad \|f\|_\infty\le B,\qquad Pf^2=(1-2d)\,\E_P\bigl[v^2\ind\{s<S\le t\}\bigr]+d^2K^2\le B\delta+B\delta^2\le2B\delta ,
\]
where $d=F_Q(t)-F_Q(s)\le\delta$ and we used $\E_P[v^2\ind\{s<S\le t\}]=\E_Q[v\ind\{s<S\le t\}]\le Bd$ and $K^2\le B$.
Since $h_t$ is right-continuous in $t$ for every fixed argument, the suprema may be taken over rational $s,t$, so the class is countable. Brackets for $h_s$ and $h_t$ of $L_2(P)$-size $\epsilon$ combine into a bracket for $h_s-h_t$ of size $2\epsilon$, and Lemma~\ref{lem:bracketH} with $K\le\sqrt B$ gives, for $\epsilon\le\sqrt B$,
\[
N_{[\,]}(2\epsilon,\Hcal_\delta,L_2(P))\le N_{[\,]}(\epsilon,\Hcal,L_2(P))^2\le\Bigl(\frac{15B}{\epsilon^2}\Bigr)^2 .
\]
Hence for $\tau\le\sqrt B$, using $\int_0^\tau\sqrt{\log(\tau/\epsilon)}\,\dd\epsilon=\tau\sqrt\pi/2$,
\[
J_{[\,]}(\tau,\Hcal_\delta,L_2(P))=\int_0^\tau\sqrt{1+\log N_{[\,]}(\epsilon,\Hcal_\delta,L_2(P))}\,\dd\epsilon\le c\,\tau\sqrt{\log(8\sqrt B/\tau)} .
\]
Lemma~3.4.2 of \citet{vandervaart1996weak}, for classes of functions bounded by $M=B$ with $Pf^2<\tau^2$, $\tau^2=2B\delta$, yields
\[
\E\,\omega_m(\delta)\le c\,J_{[\,]}(\tau)\Bigl(1+\frac{J_{[\,]}(\tau)\,B}{\tau^2\sqrt m}\Bigr)\le c'\Bigl[\sqrt{B\delta\log(8/\delta)}+\frac{B\log(8/\delta)}{\sqrt m}\Bigr].
\]
Bousquet's inequality \citep[Theorem~12.5]{boucheron2013concentration}, applied to the centered class $\Hcal_\delta\cup(-\Hcal_\delta)$ with $x=\log m$, gives with probability at least $1-1/m$
\[
\omega_m(\delta)\le\E\omega_m(\delta)+\sqrt{2\Bigl(2B\delta+\frac{2B\,\E\omega_m(\delta)}{\sqrt m}\Bigr)\log m}+\frac{B\log m}{3\sqrt m},
\]
and $\sqrt{2B\E\omega_m(\delta)\log m/\sqrt m}\le\E\omega_m(\delta)+B\log m/\sqrt m$ by the inequality between the geometric and arithmetic means. Collecting terms gives the bound for $\omega_m(\delta)$. The centered class $\{v\ind\{s<S\le t\}-Q(s<S\le t)\}$, on which $\Gm$ takes the same values as on $\{v\ind\{s<S\le t\}\}$, has the same bounds on its envelope, variance and bracketing numbers, which gives the bound for $\omega_m'(\delta)$.
\end{proof}

\subsection{Reduction to a vector sum}

In the rest of the appendix we fix the grid $b_1<\dots<b_L$ in $[\beta_0,\beta_1]$ and assume
\begin{equation}\label{eq:standing}
C_\star\,B\,m^{-1/4}\log^{5/4}(Lm)\le1
\end{equation}
for a constant $C_\star=C_\star(\beta_0,\beta_1)$ large enough; if \eqref{eq:standing} fails, the statements of Theorem~\ref{thm:rates} and Corollary~\ref{cor:finitesimul} are trivial with $C=C_\star$. Under \eqref{eq:standing}, $B\le m^{1/4}$, hence $m\ge4B^2\log m$ for $m\ge m_0$ absolute and the smallness of all the error terms below is guaranteed: in particular \eqref{eq:dev1} gives $\Pn v\in[1/2,2]$ and $|\Pn v-1|\le c\sqrt{B\log m/m}$. Define
\begin{equation}\label{eq:deltam}
\delta_m:=c_1\sqrt{\frac{B\log m}{m}}+\frac{2B}{m}\le1,
\end{equation}
and let $\mathcal E_m$ be the intersection of the events of Lemma~\ref{lem:devbounded} (with the points $q_1,\dots,q_L$) and of Lemma~\ref{lem:osc} with $\delta=\delta_m$, so that $\Prob(\mathcal E_m^c)\le6/m$.

\begin{lemma}[weighted Bahadur--Kiefer bound]\label{lem:bk}
On $\mathcal E_m$, for every $j\le L$,
\begin{align*}
\bigl|\sqrt m\,\Cm(b_j)-U_m(q_j)\bigr|&\le\epsilon_m:=c\Bigl[B^{3/4}\Bigl(\frac{\log m}{m}\Bigr)^{1/4}\sqrt{\log(8m/\delta_m)}+\frac{B\log(Lm)}{\sqrt m}\Bigr]\\
&\le c'\,B\,m^{-1/4}\log(Lm).
\end{align*}
\end{lemma}

\begin{proof}
The exact representation in the proof of Theorem~\ref{thm:coverage} is
\[
\sqrt m\,\Cm(b_j)=\Zm(\tauhat(b_j))-\sqrt m\,r_m(b_j),\qquad 0\le r_m(b_j)=\Fhat(\tauhat(b_j))-b_j\le\max_i\hat w_i\le\frac{2B}{m}.
\]
Moreover $|F_Q(\tauhat(b_j))-b_j|\le\|\Fhat-F_Q\|_\infty+\max_i\hat w_i\le\delta_m$ on $\mathcal E_m$, and $F_Q(q_j)=b_j$, so the pair $(\tauhat(b_j),q_j)$ is admissible in the definition of $\omega_m(\delta_m)$. Therefore
\begin{align*}
\bigl|\Zm(\tauhat(b_j))-U_m(q_j)\bigr|&\le\frac{|U_m(\tauhat(b_j))-U_m(q_j)|}{\Pn v}+|U_m(q_j)|\,\frac{|1-\Pn v|}{\Pn v}\\
&\le2\,\omega_m(\delta_m)+2\,|U_m(q_j)|\,|1-\Pn v| .
\end{align*}
Inserting the bounds of Lemmas~\ref{lem:devbounded} and~\ref{lem:osc},
\begin{gather*}
\omega_m(\delta_m)\le c\Bigl[\sqrt{B\delta_m\log(8m/\delta_m)}+\frac{B\log(8m/\delta_m)}{\sqrt m}\Bigr],\\
|U_m(q_j)|\,|1-\Pn v|\le c\sqrt{B\log(Lm)}\sqrt{\frac{B\log m}{m}},
\end{gather*}
and $\sqrt m\,r_m\le2B/\sqrt m$, gives the first inequality, with $B\delta_m\le c B^{3/2}\sqrt{\log m/m}$. The second follows from $B^{3/4}\le B$ and $B\log(Lm)/\sqrt m\le B m^{-1/4}\log(Lm)$.
\end{proof}

\begin{lemma}[variance estimation]\label{lem:sigmahat}
On $\mathcal E_m$ intersected with an event of probability at least $1-2/m$, for every $j\le L$,
\[
|\hat\sigma_m^2(b_j)-\sigma^2(b_j)|\le c\,B^{3/2}\sqrt{\frac{\log(Lm)}{m}},\qquad\text{hence}\qquad
\Bigl|\frac{\sigma(b_j)}{\hat\sigma_m(b_j)}-1\Bigr|\le\frac{c\,B^{3/2}}{\sigma_0^2}\sqrt{\frac{\log(Lm)}{m}}\le\frac12 .
\]
\end{lemma}

\begin{proof}
Write $g_{t,b}:=v^2(\ind\{S\le t\}-b)^2$, so that $\hat\sigma^2_m(b_j)=\Pn g_{\tauhat(b_j),b_j}/(\Pn v)^2$ and $\sigma^2(b_j)=Pg_{q_j,b_j}$. Since $(\ind\{S\le t\}-b)^2-(\ind\{S\le q\}-b)^2=(\ind\{S\le t\}-\ind\{S\le q\})(\ind\{S\le t\}+\ind\{S\le q\}-2b)$ and the second factor has absolute value at most one when the first is nonzero,
\[
|\Pn g_{\tauhat(b_j),b_j}-\Pn g_{q_j,b_j}|\le B\,\Pn\bigl[v\ind\{S\text{ between }\tauhat(b_j),q_j\}\bigr]
\le B\Bigl[\frac{\omega_m'(\delta_m)}{\sqrt m}+\delta_m\Bigr]\le c\,B\delta_m ,
\]
because $P[v\ind\{s<S\le t\}]=Q(s<S\le t)\le\delta_m$ for admissible pairs. Bernstein's inequality for $g_{q_j,b_j}$, bounded by $B^2$ with variance at most $\E_Pv^4\le B^2K^2\le B^3$, at level $1/(Lm)$ with a union bound, gives $|\Pn g_{q_j,b_j}-Pg_{q_j,b_j}|\le c\sqrt{B^3\log(Lm)/m}+cB^2\log(Lm)/m$ for all $j$ with probability at least $1-1/m$. Finally $|(\Pn v)^{-2}-1|\le c|\Pn v-1|\le c\sqrt{B\log m/m}$ and $\Pn g\le 2B$. Collecting terms and using \eqref{eq:standing} gives the first bound; the second follows from $\sigma(b_j)\ge\sigma_0$ and \eqref{eq:standing} with $C_\star$ large enough.
\end{proof}

\subsection{Gaussian approximation and anti-concentration}

\begin{lemma}[high-dimensional CLT]\label{lem:cckk}
Let $X_i:=(h_{q_j}(Z_i)/\sigma_j)_{j\le L}\in\R^L$ and let
\[
T^U_m:=\max_{j\le L}\frac{\Gm h_{q_j}}{\sigma_j}=\max_{j\le L}\frac{1}{\sqrt m}\sum_{i=1}^mX_{ij}.
\]
There is an absolute constant $c_3$ such that for $m\ge3$ (the case $m=2$ is covered by enlarging the constant)
\[
\rho_m:=\sup_{x\in\R}\bigl|\Prob(T^U_m\le x)-\Prob(T\le x)\bigr|\le c_3\Bigl(\frac{B^2\log^5(Lm)}{\sigma_0^2\,m}\Bigr)^{1/4}.
\]
The same bound holds for the two-sided statistics $\max_j|\Gm h_{q_j}|/\sigma_j$ and $T^{|\cdot|}$.
\end{lemma}

\begin{proof}
The vectors $X_i$ are i.i.d.\ with $\E X_{ij}=0$, $\E X_{ij}^2=1$ and $|X_{ij}|\le B/\sigma_0$ by \eqref{eq:sigmalower}. With $B_n:=B/(\sigma_0\log2)\ge1$, the exponential-moment Condition~E of \citet{chernozhukov2022improved} holds in either of its forms, $\E\exp(|X_{ij}|/B_n)\le2$ and $\E\exp(X_{ij}^2/B_n^2)\le2$, because $|X_{ij}|/B_n\le\log2<1$; their Condition~M holds with $b_1=b_2=1$, since $\E X_{ij}^4\le(B/\sigma_0)^2\le B_n^2$. (If $L=1$ we duplicate the coordinate, which changes nothing.) Theorem~2.1 of that paper gives, for every $\alpha\in(0,1)$,
\[
\bigl|\Prob\bigl(T^U_m>c^G_{1-\alpha}\bigr)-\alpha\bigr|\le c\Bigl(\frac{B_n^2\log^5(Lm)}{m}\Bigr)^{1/4},
\]
where $c^G_{1-\alpha}$ is the $(1-\alpha)$-quantile of $T=\max_j\Lb(b_j)/\sigma_j$, whose covariance matrix is the common covariance of the $X_i$. The law of $T$ is absolutely continuous, so its quantiles sweep its support; outside the support both probabilities are within the same bound of $0$ or $1$ by monotonicity. Hence the bound holds uniformly in $x$. For the two-sided version apply the same theorem to the $2L$-dimensional vectors $(X_i,-X_i)$.
\end{proof}

\begin{lemma}[anti-concentration]\label{lem:nazarov}
Let $Y$ be a centered Gaussian vector in $\R^L$ with $\Var Y_j\ge b>0$ for all $j$. Then for every $x\in\R$ and $\epsilon>0$,
\[
\Prob\Bigl(x<\max_{j\le L}Y_j\le x+\epsilon\Bigr)\le\frac{\epsilon}{\sqrt b}\bigl(\sqrt{2\log L}+2\bigr).
\]
\end{lemma}

\begin{proof}
This is Nazarov's inequality \citep{nazarov2003maximal} in the form of \citet{chernozhukov2017detailed}: $\Prob(Y\le y+\epsilon\mathbf1)-\Prob(Y\le y)\le\epsilon(\sqrt{2\log L}+2)/\sqrt b$ for all $y\in\R^L$, applied at $y=x\mathbf1$.
\end{proof}

\subsection{Proof of Theorem~\ref{thm:rates}}

\begin{proof}[Proof of \eqref{eq:rate-gauss}]
Assume \eqref{eq:standing}. On $\mathcal E_m$, Lemma~\ref{lem:bk} gives $|T_m-T^U_m|\le\epsilon_m/\sigma_0$. Hence for every $x$, with $\eta_m:=\epsilon_m/\sigma_0$,
\begin{align*}
\Prob(T_m\le x)&\le\Prob(T^U_m\le x+\eta_m)+\Prob(\mathcal E_m^c)\le\Prob(T\le x+\eta_m)+\rho_m+\frac6m\\
&\le\Prob(T\le x)+\eta_m\bigl(\sqrt{2\log L}+2\bigr)+\rho_m+\frac6m,
\end{align*}
by Lemma~\ref{lem:cckk} and Lemma~\ref{lem:nazarov} applied to $Y_j=\Lb(b_j)/\sigma_j$, which has unit variances. The reverse inequality is proved in the same way. Since
\begin{gather*}
\eta_m\bigl(\sqrt{2\log L}+2\bigr)\le\frac{c}{\sigma_0}\Bigl[B^{3/4}m^{-1/4}\log^{3/4}m+\frac{B\log(Lm)}{\sqrt m}\Bigr]\sqrt{\log(Lm)}\le\frac{c\,B\,m^{-1/4}\log^{5/4}(Lm)}{\sigma_0},\\
\rho_m\le\frac{c_3\sqrt B\,m^{-1/4}\log^{5/4}(Lm)}{\sqrt{\sigma_0}},
\end{gather*}
where the first bound uses the two-term form of $\epsilon_m$ in Lemma~\ref{lem:bk}, $B^{3/4}\le B$, and $\log(Lm)\le m^{1/5}$, which follows from \eqref{eq:standing}; and $6/m\le6Bm^{-1/4}$. The bound \eqref{eq:rate-gauss} follows with $C$ depending only on $\sigma_0$. The two-sided case is identical, with the two-sided versions of Lemmas~\ref{lem:bk} (which holds coordinatewise) and~\ref{lem:cckk}.
\end{proof}

\begin{proof}[Proof of \eqref{eq:rate-boot}]
Let $X_i$ be as in Lemma~\ref{lem:cckk}, $\bar X:=\Pn X$, and let
\[
T^B_m:=\max_{j\le L}\frac{1}{\sqrt m}\sum_{i=1}^m\xi_i\bigl(X_{ij}-\bar X_j\bigr),\qquad \xi_i\sim N(0,1)\ \text{i.i.d., independent of the data},
\]
be the Gaussian multiplier bootstrap statistic of \citet{chernozhukov2022improved}, with conditional quantiles $c^B_{1-\delta}$. Theorem~2.2 of that paper, under the conditions verified in Lemma~\ref{lem:cckk}, gives
\begin{equation}\label{eq:cckkboot}
\sup_{\delta\in(0,1)}\bigl|\Prob\bigl(T^U_m>c^B_{1-\delta}\bigr)-\delta\bigr|\le\rho^B_m:=c_3'\Bigl(\frac{B^2\log^5(Lm)}{\sigma_0^2\,m}\Bigr)^{1/4},
\end{equation}
with another absolute constant $c_3'$.
\emph{Step 1: the bootstrap statistics are close.} Expanding the definitions,
\[
\frac{\hat\Lb^\xi_m(b_j)}{\hat\sigma_j}-\frac{1}{\sqrt m}\sum_i\xi_i(X_{ij}-\bar X_j)=\frac{A_j}{\Pn v\,\hat\sigma_j}+\Bigl(\frac{\sigma_j}{\Pn v\,\hat\sigma_j}-1\Bigr)\frac{1}{\sqrt m}\sum_i\xi_i(X_{ij}-\bar X_j),
\]
where
\begin{align*}
A_j:={}&\frac{1}{\sqrt m}\sum_i\xi_iv_i\bigl(\ind\{S_i\le\tauhat(b_j)\}-\ind\{S_i\le q_j\}\bigr)\\
&-\bigl(\Fhat(\tauhat(b_j))-b_j\bigr)\frac{1}{\sqrt m}\sum_i\xi_iv_i
+\bigl(\Pn h_{q_j}\bigr)\frac{1}{\sqrt m}\sum_i\xi_i .
\end{align*}
Conditionally on the data, the three terms of $A_j$ are centered Gaussian with variances
\begin{gather*}
\Pn\bigl[v^2(\ind\{S\le\tauhat(b_j)\}-\ind\{S\le q_j\})^2\bigr]\le B\,\Pn\bigl[v\ind\{S\text{ between }\tauhat(b_j),q_j\}\bigr]\le cB\delta_m,\\
\bigl(\Fhat(\tauhat(b_j))-b_j\bigr)^2\Pn v^2\le\frac{4B^2}{m^2}\cdot2B,\qquad
(\Pn h_{q_j})^2\le\frac{c\,B\log(Lm)}{m},
\end{gather*}
on $\mathcal E_m$, by the argument of Lemma~\ref{lem:sigmahat}, by \eqref{eq:dev3} and by \eqref{eq:dev4}. The conditional variance of $m^{-1/2}\sum_i\xi_i(X_{ij}-\bar X_j)$ is $\Pn X_j^2-\bar X_j^2$, which lies in $[1/2,2]$ on $\mathcal E_m$ intersected with an event of probability $1-1/m$ (Bernstein for $X_{ij}^2$, bounded by $B^2/\sigma_0^2$ with variance at most $B^3/\sigma_0^4$, and \eqref{eq:dev4}). Gaussian tail bounds and a union bound over $j$ therefore give, with conditional probability at least $1-1/m$,
\begin{gather*}
\max_j|A_j|\le c\sqrt{B\delta_m\log(Lm)}+\frac{cB^{3/2}\sqrt{\log m}}{m}+\frac{c\sqrt{B}\log(Lm)}{\sqrt m},\\
\max_j\Bigl|\frac{1}{\sqrt m}\sum_i\xi_i(X_{ij}-\bar X_j)\Bigr|\le c\sqrt{\log(Lm)} .
\end{gather*}
Together with Lemma~\ref{lem:sigmahat} and $|\Pn v-1|\le c\sqrt{B\log m/m}$ this yields, on an event $\mathcal E'_m$ of the data with $\Prob(\mathcal E_m'^c)\le c/m$ and with conditional probability at least $1-1/m$,
\[
\bigl|\hat T^\xi_m-T^B_m\bigr|\le\epsilon'_m:=c\Bigl[\frac{B^{3/4}\log^{3/4}(Lm)}{\sigma_0\,m^{1/4}}+\frac{B^{3/2}\log(Lm)}{\sigma_0^2\sqrt m}\Bigr]\le\frac{c'\,B\,m^{-1/4}\log(Lm)}{\sigma_0^2},
\]
where the last step uses $B\le m^{1/4}$.

\emph{Step 2: the conditional laws are close.} On $\mathcal E'_m$, $T^B_m$ is conditionally the maximum of a Gaussian vector with variances at least $1/2$, so Lemma~\ref{lem:nazarov} and Step~1 give
\[
\sup_x\bigl|\Prob_\xi(\hat T^\xi_m\le x)-\Prob_\xi(T^B_m\le x)\bigr|\le\rho'_m:=\sqrt2\,\epsilon'_m\bigl(\sqrt{2\log L}+2\bigr)+\frac1m .
\]
Consequently the conditional quantiles satisfy $c^B_{1-\delta-\rho'_m}\le\hat c_{1-\delta}\le c^B_{1-\delta+\rho'_m}$ on $\mathcal E'_m$, whenever $\delta\pm\rho'_m\in(0,1)$.

\emph{Step 3: the statistics are close.} By Lemmas~\ref{lem:bk} and~\ref{lem:sigmahat}, on $\mathcal E_m\cap\mathcal E'_m$,
\begin{align*}
|\hat T_m-T^U_m|&\le\frac{\epsilon_m}{\sigma_0}+\max_j\frac{|\sqrt m\,\Cm(b_j)|}{\sigma_j}\Bigl|\frac{\sigma_j}{\hat\sigma_j}-1\Bigr|\\
&\le\frac{\epsilon_m}{\sigma_0}+\frac{c\sqrt{B\log(Lm)}}{\sigma_0}\cdot\frac{cB^{3/2}}{\sigma_0^2}\sqrt{\frac{\log(Lm)}{m}}=:\epsilon''_m\le\frac{c\,B\,m^{-1/4}\log(Lm)}{\sigma_0^3},
\end{align*}
using $|\sqrt m\,\Cm(b_j)|\le|U_m(q_j)|+\epsilon_m\le c\sqrt{B\log(Lm)}$ and $B^2/\sqrt m\le B m^{-1/4}$.

\emph{Step 4: conclusion.} Let
\[
\theta_m:=\sqrt2\,\epsilon''_m\bigl(\sqrt{2\log L}+2\bigr),\qquad\nu_m:=\rho'_m+\theta_m,\qquad r_m:=\nu_m+\rho^B_m+\frac cm .
\]
Under \eqref{eq:standing}, $B^2\log(Lm)\le m$ and $B^4\log(Lm)\le m$, so the two-term forms of $\epsilon'_m$ and $\epsilon''_m$ give, exactly as in the proof of \eqref{eq:rate-gauss},
\begin{gather*}
\nu_m\le C\,B\,m^{-1/4}\log^{5/4}(Lm),\qquad\rho^B_m\le\frac{c\sqrt B\,m^{-1/4}\log^{5/4}(Lm)}{\sqrt{\sigma_0}},\\
\text{hence}\qquad 2\nu_m+r_m\le C\,B\,m^{-1/4}\log^{5/4}(Lm),
\end{gather*}
and $\nu_m\le1/4$ once $C_\star$ is large enough.

\emph{Interior levels.} Let $2\nu_m\le\delta\le1-2\nu_m$, so that all the quantile levels below lie in $(0,1)$. On $\mathcal E_m\cap\mathcal E'_m$, by Steps~2 and~3,
\[
\{\hat T_m>\hat c_{1-\delta}\}\subset\{T^U_m>c^B_{1-\delta-\rho'_m}-\epsilon''_m\}.
\]
Let $G^B$ denote, conditionally on the data, the Gaussian vector $m^{-1/2}\sum_i\xi_i(X_i-\bar X)$, so that $T^B_m=\max_jG^B_j$. By Lemma~\ref{lem:nazarov}, conditionally on the data,
\[
\Prob_\xi\bigl(T^B_m\le c^B_{1-\delta-\rho'_m}-\epsilon''_m\bigr)\ge1-\delta-\nu_m,\qquad\text{that is,}\qquad c^B_{1-\delta-\rho'_m}-\epsilon''_m\ge c^B_{1-\delta-\nu_m}.
\]
Hence, by \eqref{eq:cckkboot},
\[
\Prob(\hat T_m>\hat c_{1-\delta})\le\Prob\bigl(T^U_m>c^B_{1-\delta-\nu_m}\bigr)+\frac cm\le\delta+r_m ,
\]
and the lower bound is obtained symmetrically, from $\hat c_{1-\delta}\le c^B_{1-\delta+\rho'_m}$ and $c^B_{1-\delta+\rho'_m}+\epsilon''_m\le c^B_{1-\delta+\nu_m}$:
\begin{equation}\label{eq:interior}
\bigl|\Prob(\hat T_m>\hat c_{1-\delta})-\delta\bigr|\le r_m\qquad\text{for }2\nu_m\le\delta\le1-2\nu_m .
\end{equation}

\emph{Levels near the endpoints.} For every realization of the data the quantile $\hat c_{1-\delta}$ is nonincreasing in $\delta$, so the events $\{\hat T_m>\hat c_{1-\delta}\}$ increase with $\delta$ and
\[
\delta\mapsto\Prob(\hat T_m>\hat c_{1-\delta})\quad\text{is nondecreasing on }(0,1).
\]
For $\delta<2\nu_m$, monotonicity and \eqref{eq:interior} at the level $2\nu_m$ give
\[
0\le\Prob(\hat T_m>\hat c_{1-\delta})\le\Prob(\hat T_m>\hat c_{1-2\nu_m})\le2\nu_m+r_m ,
\]
and for $\delta>1-2\nu_m$, monotonicity and \eqref{eq:interior} at the level $1-2\nu_m$ give
\[
1\ge\Prob(\hat T_m>\hat c_{1-\delta})\ge\Prob(\hat T_m>\hat c_{2\nu_m})\ge1-2\nu_m-r_m .
\]
In the first case $\delta$ and the probability both lie in $[0,2\nu_m+r_m]$, in the second both lie in $[1-2\nu_m-r_m,1]$, so that
\[
\sup_{\delta\in(0,1)}\bigl|\Prob(\hat T_m>\hat c_{1-\delta})-\delta\bigr|\le2\nu_m+r_m\le C\,B\,m^{-1/4}\log^{5/4}(Lm),
\]
which is \eqref{eq:rate-boot}. The two-sided case is identical.
\end{proof}

\subsection{Proof of Corollary~\ref{cor:finitesimul}}

\begin{proof}
Assume \eqref{eq:standing} and $\delta\ge\max\{1/m,2\rho'_m\}$; for smaller $\delta$ the inflated levels are larger than those for $\delta_0:=\max\{1/m,2\rho_m'\}$, so the failure probability is at most that for $\delta_0$, which is bounded by $\delta_0+CBm^{-1/4}\log^{5/4}(Lm)$ by the case treated below, and $\delta_0\le CBm^{-1/4}\log^{5/4}(Lm)$. 

\emph{The critical value is nonnegative.} This holds for every $\delta<1/2$, with no condition on the approximation errors. On $\mathcal E_m\cap\mathcal E'_m$ the estimates $\hat\sigma_j$ are positive (Lemma~\ref{lem:sigmahat}) and, the multipliers being Gaussian, each coordinate $\hat\Lb^\xi_m(b_j)/\hat\sigma_j$ of $\hat T^\xi_m$ is, conditionally on the data, a centered Gaussian variable, possibly degenerate at zero. Hence
\[
\Prob_\xi\bigl(\hat T^\xi_m<0\bigr)\le\Prob_\xi\bigl(\hat\Lb^\xi_m(b_1)<0\bigr)\le\frac12<1-\delta ,
\]
so that no $x<0$ satisfies $\Prob_\xi(\hat T^\xi_m\le x)\ge1-\delta$, and by the definition \eqref{eq:bootmax} of the conditional quantile
\[
\hat c_{1-\delta}\ge0\qquad\text{for every }\delta\in(0,1/2)
\]
(with $\hat c_{1-\delta}=0$ if all coordinates vanish). In particular $b_j'\ge b_j$.

\emph{The failure event.} Since $P_e(\tauhat(b'_j))>1-b_j$ means $F_Q(\tauhat(b_j'))<b_j=b_j'-\hat c_{1-\delta}\hat\sigma_j/\sqrt m$,
\[
\Bigl\{\exists j:\ P_e\bigl(\tauhat(b_j')\bigr)>1-b_j\Bigr\}=\Bigl\{\max_{j\le L}\frac{\sqrt m\,\Cm(b_j')}{\hat\sigma_j}>\hat c_{1-\delta}\Bigr\}.
\]
\emph{The inflated levels are close to the grid.} On $\mathcal E'_m$ the conditional law of $\hat T^\xi_m$ is within $\rho'_m$ of that of $T^B_m$ in Kolmogorov distance, and the coordinates of $T^B_m$ are conditionally Gaussian with variances at most $2$, so the union bound gives $\Prob_\xi(\hat T^\xi_m>x)\le Le^{-x^2/4}+\rho'_m\le\delta$ at $x=2\sqrt{\log(L/(\delta-\rho'_m))}$, whence $\hat c_{1-\delta}\le2\sqrt{\log(2Lm)}$. With Lemma~\ref{lem:sigmahat},
\[
|b_j'-b_j|=\frac{\hat c_{1-\delta}\hat\sigma_j}{\sqrt m}\le\delta_m':=c\sqrt{\frac{B\log(Lm)}{m}}\qquad\text{on }\mathcal E_m\cap\mathcal E'_m ,
\]
and under \eqref{eq:standing} with $C_\star$ large the levels $b_j'$ lie in $[\beta_0/2,(1+\beta_1)/2]$.

\emph{The statistic at the inflated levels is close to $T^U_m$.} The levels $b_j'$ are random, but the bounds of Lemmas~\ref{lem:devbounded} and~\ref{lem:osc} are uniform: $|F_Q(\tauhat(b))-b|\le\delta_m$ for every $b$ on $\mathcal E_m$, and $\omega_m(\cdot)$ is a supremum over all admissible pairs. Hence, on $\mathcal E_m\cap\mathcal E'_m$ intersected with the event of Lemma~\ref{lem:osc} for $\delta=\delta'_m$, the argument of Lemma~\ref{lem:bk} gives
\begin{align*}
\bigl|\sqrt m\,\Cm(b_j')-U_m(q_j)\bigr|&\le2\omega_m(\delta_m)+\omega_m(\delta'_m)+2\bigl(|U_m(q_j)|+\omega_m(\delta_m')\bigr)|1-\Pn v|+\frac{2B}{\sqrt m}\\
&\le c\Bigl[B^{3/4}m^{-1/4}\log^{3/4}(Lm)+\frac{B\log(Lm)}{\sqrt m}\Bigr],
\end{align*}
where $|U_m(q_Q(b_j'))-U_m(q_j)|\le\omega_m(\delta'_m)$ was used because $|F_Q(q_Q(b_j'))-F_Q(q_j)|=|b_j'-b_j|\le\delta'_m$. Dividing by $\hat\sigma_j\ge\sigma_0/2$ and using Lemma~\ref{lem:sigmahat} as in Step~3 of the proof of \eqref{eq:rate-boot},
\[
\Bigl|\max_j\frac{\sqrt m\,\Cm(b_j')}{\hat\sigma_j}-T^U_m\Bigr|\le\frac{c}{\sigma_0^3}\Bigl[B^{3/4}m^{-1/4}\log^{3/4}(Lm)+\frac{B^2\log(Lm)}{\sqrt m}\Bigr],
\]
which has the two-term form of $\epsilon''_m$, so that under \eqref{eq:standing} its product with the factor $\sqrt{\log(Lm)}$ of the anti-concentration step is at most $C\,B\,m^{-1/4}\log^{5/4}(Lm)$. The upper bound in the interior case of Step~4 of that proof uses only that the quantile level $1-\delta-\nu_m$ lies in $(0,1)$. With $\epsilon''_m$ replaced by the quantity above, $\nu_m\le1/4$ under \eqref{eq:standing} with $C_\star$ large, and $\delta<1/2$, so that bound applies and gives
\[
\Prob\Bigl(\max_{j}\frac{\sqrt m\,\Cm(b_j')}{\hat\sigma_j}>\hat c_{1-\delta}\Bigr)\le\delta+C\,B\,m^{-1/4}\log^{5/4}(Lm).
\]
\emph{The sets of \citet{tibshirani2019conformal}.} The atom at $+\infty$ can only raise the level of the weighted quantile, so the coverage of the set at level $b_j'$ is at least $F_Q(\tauhat(b_j'))$, and the failure event for these sets is contained in the failure event treated above.
\end{proof}%

\section{Proofs for Section~\ref{sec:efficiency}}\label{app:efficiency}

We use the framework of \citet[Chapter~25]{vandervaart1998asymptotic}: for a model $\mathcal P$ and $P\in\mathcal P$, the tangent set $\dot{\mathcal P}_P$ is the set of score functions of differentiable paths $\eta\mapsto P_\eta$ in $\mathcal P$ with $P_0=P$; a functional $\theta$ is pathwise differentiable at $P$ if there is $\dot\theta_P\in L_2(P)$ with $\partial_\eta\theta(P_\eta)|_{\eta=0}=\E_P[\dot\theta_Pg]$ for every path with score $g$; the efficient influence function is the projection of $\dot\theta_P$ onto the closed linear span of the tangent set, and the efficiency bound is its second moment. An estimator that is asymptotically linear with influence function equal to the efficient influence function is regular and efficient \citep[Lemma~25.23]{vandervaart1998asymptotic}.

\subsection{Proof of Theorem~\ref{thm:eff}}

\begin{proof}
\emph{Tangent space.} A differentiable path in $\mathcal M_1$ has the form $P_\eta=P_{X,\eta}\otimes P_{S|X,\eta}$ with score $g(x,s)=a(x)+b(s,x)$, where $a$ is the score of the marginal path and $b$ that of the conditional path, so that $\E_P[a]=0$ and $\E_P[b\mid X]=0$. Differentiating the constraint $\int v\,\dd P_{X,\eta}=1$ gives $\E_P[v\,a(X)]=0$. Conversely, for bounded $a$ with $\E a=\E[va]=0$ and bounded $b$ with $\E[b\mid X]=0$, the path $\dd P_\eta=(1+\eta a)(1+\eta b)\,\dd P$ stays in $\mathcal M_1$ for small $|\eta|$. Hence the closed linear span of the tangent set is
\begin{gather*}
\overline{\dot{\mathcal M}}_1=\mathcal T_X\oplus\mathcal T_{S|X},\\
\mathcal T_X:=\{a(X)\in L_2(P):\E a=0,\ \E[va]=0\},\qquad
\mathcal T_{S|X}:=\{b(S,X)\in L_2(P):\E[b\mid X]=0\},
\end{gather*}
the two spaces being orthogonal. Bounded functions are dense in each: for $b\in\mathcal T_{S|X}$ the truncations $b\ind\{|b|\le n\}-\E[b\ind\{|b|\le n\}\mid X]$ are bounded, conditionally centered and converge to $b$ in $L_2$; for $a\in\mathcal T_X$, fix $c\ge1$ with $\Prob(v>c)>0$ (possible since $\E v=1$ and $v$ is not a.s.\ constant), so that the bounded functions $e_1=1$ and $e_2=\ind\{v>c\}$ satisfy $\E[(v-1)e_2]>0=\E[(v-1)e_1]$, and correct the truncation $a\ind\{|a|\le n\}$ by a combination $\alpha_ne_1+\beta_ne_2$ that restores $\E[\cdot]=\E[v\,\cdot]=0$; the coefficients tend to zero.

\emph{Pathwise derivative.} Along a path with score $g=a+b$,
\begin{align*}
\partial_\eta\theta_t(P_\eta)\big|_{\eta=0}&=\partial_\eta\int\!\!\int v(x)\ind\{s\le t\}\,\dd P_{S|X,\eta}(s\mid x)\,\dd P_{X,\eta}(x)\Big|_{\eta=0}\\
&=\E_P\bigl[v\ind\{S\le t\}\{a(X)+b(S,X)\}\bigr]=\E_P\bigl[g_t\,g\bigr],
\end{align*}
with $g_t:=v\ind\{S\le t\}-F_Q(t)\in L_2(P)$, using $\E_P g=0$; the interchange of derivative and integral is justified by $v\in L_2(P)$ and boundedness of the scores of the generating paths.

\emph{Projection.} Decompose $g_t=\{vm_t(X)-F_Q(t)\}+v\{\ind\{S\le t\}-m_t(X)\}$. The second term lies in $\mathcal T_{S|X}$. The first is a function of $X$ with mean zero, and its projection onto $\mathcal T_X$, the orthogonal complement of $\mathrm{span}\{1,v\}$ in $L_2(P_X)$, is obtained by removing its component along $v-1$:
\[
\Pi_{\mathcal T_X}\bigl(vm_t-F_Q(t)\bigr)=vm_t-F_Q(t)-c_t(v-1),\qquad
c_t=\frac{\E_P[(vm_t-F_Q(t))(v-1)]}{\E_P(v-1)^2}=\frac{G(t)-F_Q(t)}{K^2-1},
\]
because $\E_P[v^2m_t]=G(t)$, $\E_P[vm_t]=F_Q(t)$, $\E_Pv=1$ and $\E_P(v-1)^2=K^2-1$. Adding the two components gives $\psi_t=v\ind\{S\le t\}-F_Q(t)-c_t(v-1)$, and
\[
\E_P\psi_t^2=\Var_P(g_t)-2c_t\Cov_P(g_t,v-1)+c_t^2(K^2-1)=\bigl(G(t)-F_Q(t)^2\bigr)-\frac{(G(t)-F_Q(t))^2}{K^2-1},
\]
since $\Cov_P(g_t,v-1)=G(t)-F_Q(t)$ and $\Var_P(g_t)=\E_P[v^2\ind\{S\le t\}]-F_Q(t)^2$.

\emph{The self-normalized estimator.} By Lemma~\ref{lem:identity}, $\sqrt m(\Fhat(t)-F_Q(t))=\Gm h_t/\Pn v=\Gm h_t+o_p(1)$, so $\Fhat(t)$ is asymptotically linear with influence function $h_t=v(\ind\{S\le t\}-F_Q(t))=g_t-F_Q(t)(v-1)$. An asymptotically linear estimator is regular if and only if its influence function has the same inner product as $\dot\theta_P$ with every tangent vector \citep[Section~25.3]{vandervaart1998asymptotic}; here $h_t-g_t=-F_Q(t)(v-1)$ is orthogonal to $\mathcal T_X$ by definition of $\mathcal T_X$ and to $\mathcal T_{S|X}$ because it is a function of $X$, so $\Fhat(t)$ is regular. Its asymptotic variance is
\[
\E_Ph_t^2=\E_P\psi_t^2+(c_t-F_Q(t))^2(K^2-1)=V_1(t)+\frac{(K^2F_Q(t)-G(t))^2}{K^2-1},
\]
by orthogonality of $\psi_t$ and $v-1$, and it equals $V_1(t)$ if and only if $G(t)=K^2F_Q(t)$, that is, $\Cov_Q(v,\ind\{S\le t\})=\E_Q[v\ind\{S\le t\}]-\E_Qv\,Q(S\le t)=0$.
\end{proof}

\subsection{Proof of Theorem~\ref{thm:effthreshold}}

\begin{proof}
\emph{Uniform asymptotic linearity.} Since $\E_P[v^2\ind\{S\le t\}]=G(t)$ and the classes $\{v^2\ind\{S\le t\}\}$, $\{v\ind\{S\le t\}\}$ have integrable envelopes, the Glivenko--Cantelli theorem and $\Pn v^2\to K^2>1$ give $\sup_t|\hat c_m(t)-c_t|\toP0$. Expanding,
\begin{align*}
\sqrt m\bigl(\tilde F_m(t)-F_Q(t)\bigr)&=\Gm\bigl(v\ind\{S\le t\}\bigr)-\hat c_m(t)\,\Gm v\\
&=\Gm h_t+\bigl(F_Q(t)-\hat c_m(t)\bigr)\Gm v
=\Gm\psi_t+\bigl(c_t-\hat c_m(t)\bigr)\Gm v,
\end{align*}
where we used $\Gm(v\ind\{S\le t\})=\Gm(v\ind\{S\le t\}-F_Q(t))=\Gm(h_t+F_Q(t)(v-1))$ and $\Gm(v-1)=\Gm v$. Since $\Gm v=O_p(1)$, the remainder is $o_p(1)$ uniformly in $t$, which is the first claim; by Lemma~\ref{lem:identity} and Theorem~\ref{thm:fclt} the process $\Gm\psi_t=\Gm h_t+(F_Q(t)-c_t)\Gm v$ converges jointly with $\Gm v$ to $(\GQ(t)+(F_Q(t)-c_t)\mathcal N,\mathcal N)$, where $\mathcal N$ is the limit of $\Gm v$, with $\Var\mathcal N=K^2-1$ and $\Cov(\GQ(t),\mathcal N)=\Cov_P(h_t,v)=G(t)-K^2F_Q(t)$.

\emph{Quantile.} $\tilde F_m$ is a right-continuous step function with jumps at the $S_i$ of size
\[
\frac{v_i}{m}-\frac{(\Pn v-1)}{\Pn v^2-1}\Bigl(\frac{v_i^2}{m}-\hat w_i\Bigr),
\]
so that $\max_i|\text{jump}_i|\le\max_iv_i/m+O_p(m^{-1/2})\max_i(v_i^2/m+\hat w_i)=o_p(m^{-1/2})$, by Lemma~\ref{lem:maxw} and $\max_iv_i^2/m\toP0$ (which follows from $\E v^2<\infty$ as in that lemma). By definition of $\tilde\tau(\beta)$ as an infimum, $\tilde F_m(t)<\beta$ for $t<\tilde\tau(\beta)$ and $\tilde F_m(\tilde\tau(\beta))\ge\beta$, so
\begin{gather*}
0\le\tilde r_m(\beta):=\tilde F_m(\tilde\tau(\beta))-\beta\le\max_i|\text{jump}_i|=o_p(m^{-1/2}),\\
\sqrt m\bigl(\beta-F_Q(\tilde\tau(\beta))\bigr)=\sqrt m\bigl(\tilde F_m-F_Q\bigr)(\tilde\tau(\beta))-\sqrt m\,\tilde r_m(\beta).
\end{gather*}
Moreover $\sup_t|\tilde F_m-F_Q|\toP0$ by the Glivenko--Cantelli theorem, so $|F_Q(\tilde\tau(\beta))-\beta|\le\sup_t|\tilde F_m-F_Q|+\max_i|\text{jump}_i|\toP0$ uniformly in $\beta$, and $F_Q(q_Q(\beta))=\beta$.

\emph{Equicontinuity.} Write $\sqrt m(\tilde F_m-F_Q)(t)=\Gm h_t+(F_Q(t)-\hat c_m(t))\Gm v$. The first term is handled by Lemma~\ref{lem:rho} and \eqref{eq:equicont} exactly as in the proof of Theorem~\ref{thm:coverage}: $\sup_\beta|\Gm h_{\tilde\tau(\beta)}-\Gm h_{q_Q(\beta)}|\toP0$. For the second, $\sup_t|\hat c_m(t)-c_t|\toP0$, $\Gm v=O_p(1)$, and $t\mapsto c_t$ is continuous in the $F_Q$-scale: $c_{q_Q(\beta)}=(\tilde G(\beta)-\beta)/(K^2-1)$ and $\tilde G$ is continuous because $G\ll F_Q$ (Lemma~\ref{lem:rho}). Hence $(F_Q-\hat c_m)(\tilde\tau(\beta))-(F_Q-c)(q_Q(\beta))\toP0$ uniformly in $\beta$. Therefore
\[
\sqrt m\bigl(\beta-F_Q(\tilde\tau(\beta))\bigr)=\Gm h_{q_Q(\beta)}+\bigl(\beta-\tilde c(\beta)\bigr)\Gm v+o_p(1)\weak\Lb(\beta)+\bigl(\beta-\tilde c(\beta)\bigr)\mathcal N
\]
in $\linf[\beta_0,\beta_1]$, by Theorem~\ref{thm:coverage}, the joint convergence of $(\Gm h_{q_Q(\cdot)},\Gm v)$, and the continuous mapping theorem. The variance is
\[
\sigma^2(\beta)+2(\beta-\tilde c)(\tilde G-K^2\beta)+(\beta-\tilde c)^2(K^2-1)=\sigma^2(\beta)-\frac{(\beta K^2-\tilde G(\beta))^2}{K^2-1},
\]
because $(\beta-\tilde c(\beta))(K^2-1)=\beta K^2-\tilde G(\beta)$.

\emph{Lower bound.} Let $\tau_m(\beta)$ be a regular estimator of $q_Q(\beta)$ in $\mathcal M_1$ and suppose $F_Q$ is differentiable at $q=q_Q(\beta)$ with $f_Q(q)>0$. The functional $P\mapsto q_Q(\beta)=F_Q^{-1}(\beta)$ is pathwise differentiable with efficient influence function $-\psi_q/f_Q(q)$, by the chain rule for the inverse map \citep[Lemma~21.3 and Section~25.7]{vandervaart1998asymptotic}; the inverse map is Hadamard differentiable at $F_Q$ tangentially to functions continuous at $q$, and the derivative $t\mapsto\E_P[g_tg]$ of $t\mapsto\theta_t$ is continuous at $q$ because $v\ind\{S=q\}=0$ almost surely ($\E_P[v\ind\{S=q\}]=Q(S=q)=0$). By the convolution theorem \citep[Theorem~25.20]{vandervaart1998asymptotic} the limit law of $\sqrt m(\tau_m(\beta)-q)$ is the convolution of $N(0,V_1(q)/f_Q(q)^2)$ with another law, so its variance, when finite, is at least $V_1(q)/f_Q(q)^2$. Since $\sqrt m(\beta-F_Q(\tau_m(\beta)))=-f_Q(q)\sqrt m(\tau_m(\beta)-q)+o_p(1)$ by the delta method, the asymptotic variance of the coverage error is at least $V_1(q)=\sigma^2_{\rm eff}(\beta)$; the identity $V_1(q)=\sigma^2_{\rm eff}(\beta)$ is the variance formula of Theorem~\ref{thm:eff} with $F_Q(q)=\beta$ and $G(q)=\tilde G(\beta)$.
\end{proof}

\subsection{The variance estimator and degenerate levels}\label{app:degenerate}

The inflation of the efficient threshold described after Theorem~\ref{thm:effthreshold} rests on the following statement; part~(c) requires a nondegenerate limit, and Example~\ref{ex:degenerate} shows what can happen without it.

\begin{proposition}[inflation of the efficient threshold]\label{prop:sigmaeff}
Under the assumptions of Theorem~\ref{thm:effthreshold}, let $\beta\in(\beta_0,\beta_1)$, $\delta\in(0,1)$ and
\[
\hat\sigma^2_{\rm eff}(\beta):=\Pn\hat\psi^2,\qquad\hat\psi:=v\,\ind\{S\le\tilde\tau(\beta)\}-\beta-\hat c_m(\tilde\tau(\beta))\,(v-1).
\]
\begin{enumerate}[label=(\alph*),leftmargin=*,itemsep=1pt]
\item $\hat\sigma^2_{\rm eff}(\beta)\toP\sigma^2_{\rm eff}(\beta)$.
\item $\sigma^2_{\rm eff}(\beta)=0$ if and only if $v$ is $P$-almost surely constant on $\{S\le q_Q(\beta)\}$ and $P$-almost surely constant on $\{S>q_Q(\beta)\}$.
\item If $\sigma^2_{\rm eff}(\beta)>0$, then
\[
\Prob\Bigl(P_e\Bigl(\tilde\tau\Bigl(\beta+\frac{z_{1-\delta}\,\hat\sigma_{\rm eff}(\beta)}{\sqrt m}\Bigr)\Bigr)>1-\beta\Bigr)\longrightarrow\delta .
\]
\end{enumerate}
\end{proposition}

\begin{proof}
Write $q:=q_Q(\beta)$, $A:=\{S\le q\}$, $\tilde c:=\tilde c(\beta)$ and
\[
\psi_q=v\,\ind_A-\beta-\tilde c\,(v-1),\qquad\E_P\psi_q^2=V_1(q)=\sigma^2_{\rm eff}(\beta),
\]
the last identity being the one at the end of the previous proof.

(a) Let $\hat c:=\hat c_m(\tilde\tau(\beta))$. Expanding the square,
\begin{align*}
\Pn\hat\psi^2&=\Pn\bigl(v^2\ind\{S\le\tilde\tau(\beta)\}\bigr)-2\beta\,\Pn\bigl(v\ind\{S\le\tilde\tau(\beta)\}\bigr)+\beta^2\\
&\quad-2\hat c\,\Bigl[\Pn\bigl(v^2\ind\{S\le\tilde\tau(\beta)\}\bigr)-\Pn\bigl(v\ind\{S\le\tilde\tau(\beta)\}\bigr)\Bigr]+2\beta\hat c\,(\Pn v-1)+\hat c^{\,2}\,\Pn(v-1)^2 .
\end{align*}
By the Glivenko--Cantelli theorem for the classes $\{v\ind\{S\le t\}\}$ and $\{v^2\ind\{S\le t\}\}$, by $F_Q(\tilde\tau(\beta))\toP\beta$ and $\sup_t|\hat c_m(t)-c_t|\toP0$ (previous proof), and by the continuity of $\tilde G$ and of $\tilde c$,
\[
\Pn\bigl(v\ind\{S\le\tilde\tau(\beta)\}\bigr)\toP\beta,\qquad\Pn\bigl(v^2\ind\{S\le\tilde\tau(\beta)\}\bigr)\toP\tilde G(\beta),\qquad\hat c\toP\tilde c,
\]
where we used $G(t)=\tilde G(F_Q(t))$, valid because $G\ll F_Q$. With $\Pn v\toP1$, $\Pn(v-1)^2\toP K^2-1$ and $\tilde c\,(K^2-1)=\tilde G(\beta)-\beta$,
\[
\Pn\hat\psi^2\toP\tilde G(\beta)-\beta^2-2\tilde c\,\bigl(\tilde G(\beta)-\beta\bigr)+\tilde c^{\,2}(K^2-1)=\tilde G(\beta)-\beta^2-\frac{\bigl(\tilde G(\beta)-\beta\bigr)^2}{K^2-1}=V_1(q).
\]

(b) Since
\[
\E_P\psi_q=0,\qquad\E_P\bigl[\psi_q\,(v-1)\bigr]=\tilde G(\beta)-\beta-\tilde c\,(K^2-1)=0,
\]
$\psi_q$ is the residual of the projection of $v\ind_A$ onto the span of $1$ and $v$ in $L^2(P)$, so that
\[
\sigma^2_{\rm eff}(\beta)=0\iff v\,\ind_A=a_0+a_1v\quad P\text{-almost surely, for some constants }a_0,a_1 .
\]
The right-hand side reads
\[
(1-a_1)\,v=a_0\ \text{ on }A,\qquad a_1v=-a_0\ \text{ on }A^c .
\]
If $v$ were not almost surely constant on $A$, the first relation would force $a_1=1$ and $a_0=0$, and the second would give $v=0$ on $A^c$, contradicting
\[
\E_P\bigl[v\,\ind_{A^c}\bigr]=1-F_Q(q)=1-\beta>0 ;
\]
if $v$ were not almost surely constant on $A^c$, the second relation would force $a_1=a_0=0$ and the first $v=0$ on $A$, contradicting $\E_P[v\ind_A]=\beta>0$. Conversely, if $v=a$ on $A$ and $v=b$ on $A^c$, then $a\ne b$ because $K^2>1$, and
\[
v\,\ind_A=\frac{a\,(v-b)}{a-b}
\]
is an affine function of $v$.

(c) Let
\[
Z_m(u):=\sqrt m\,\bigl(u-F_Q(\tilde\tau(u))\bigr),\qquad\beta_m':=\beta+\frac{z_{1-\delta}\,\hat\sigma_{\rm eff}(\beta)}{\sqrt m},
\]
so that, since $P_e=1-F_Q$,
\[
\Bigl\{P_e\bigl(\tilde\tau(\beta_m')\bigr)>1-\beta\Bigr\}=\bigl\{F_Q(\tilde\tau(\beta_m'))<\beta\bigr\}=\bigl\{Z_m(\beta_m')>z_{1-\delta}\,\hat\sigma_{\rm eff}(\beta)\bigr\}.
\]
By Theorem~\ref{thm:effthreshold}, $Z_m\weak\tilde\Lb$ in $\linf[\beta_0,\beta_1]$, and $\tilde\Lb$ has continuous sample paths ($\Lb$ has, as shown in the proof of Corollary~\ref{cor:simultaneous}, and $\tilde c$ is continuous); hence, by the continuous mapping theorem applied to the modulus of continuity,
\[
\lim_{\eta\downarrow0}\limsup_{m\to\infty}\Prob\Bigl(\sup_{|u-u'|\le\eta}\bigl|Z_m(u)-Z_m(u')\bigr|>\varepsilon\Bigr)=0\qquad\text{for every }\varepsilon>0 .
\]
By (a), $\beta_m'-\beta=O_p(m^{-1/2})$, so $\beta_m'\in[\beta_0,\beta_1]$ with probability tending to one and
\[
Z_m(\beta_m')=Z_m(\beta)+o_p(1)\weak N\bigl(0,\sigma^2_{\rm eff}(\beta)\bigr),\qquad z_{1-\delta}\,\hat\sigma_{\rm eff}(\beta)\toP z_{1-\delta}\,\sigma_{\rm eff}(\beta).
\]
When $\sigma_{\rm eff}(\beta)>0$ the limit law is continuous, and Slutsky's lemma gives
\[
\Prob\bigl(Z_m(\beta_m')>z_{1-\delta}\,\hat\sigma_{\rm eff}(\beta)\bigr)\longrightarrow1-\Phi(z_{1-\delta})=\delta .
\]
When $\sigma_{\rm eff}(\beta)=0$ the limit is a point mass at zero, at which the threshold $z_{1-\delta}\hat\sigma_{\rm eff}(\beta)\toP0$ also concentrates, and the argument gives nothing.
\end{proof}

\begin{example}[a degenerate level]\label{ex:degenerate}
Let $S=X$, $Q_X=\mathrm{Unif}(0,1)$ and let $P_X$ have the density
\begin{gather*}
p_X(x)=\tfrac12\,\ind\{0<x\le\tfrac12\}+\tfrac32\,\ind\{\tfrac12<x<1\},\\
\text{so that}\qquad v(x)=2\,\ind\{x\le\tfrac12\}+\tfrac23\,\ind\{x>\tfrac12\},\qquad K^2=\tfrac43 ,
\end{gather*}
and the assumptions of Theorem~\ref{thm:effthreshold} hold. At $\beta=1/2$,
\[
q_Q(\beta)=\tfrac12,\qquad\tilde G(\beta)=1,\qquad\sigma^2(\beta)=\tfrac13,\qquad\sigma^2_{\rm eff}(\beta)=\frac13-\frac{(2/3-1)^2}{4/3-1}=0,
\]
in agreement with Proposition~\ref{prop:sigmaeff}(b): the constraint $\int v\,\dd P_X=1$ alone gives $P(X\le1/2)=1/4$, hence $F_Q(1/2)=1/2$. Let
\[
r_m:=\Pn\ind\{X\le\tfrac12\}\toP\tfrac14 .
\]
Then
\begin{gather*}
\Pn v=\frac{4r_m+2}{3},\qquad\Pn v^2=\frac{32r_m+4}{9},\qquad\Fhat(\tfrac12)=\frac{3r_m}{2r_m+1},\\
\hat c_m(\tfrac12)=\frac{\Pn(v^2\ind\{S\le\frac12\})-\Fhat(\frac12)}{\Pn v^2-1}=\frac{9\,r_m(8r_m+1)}{(2r_m+1)(32r_m-5)},
\end{gather*}
and therefore
\[
\tilde F_m(\tfrac12)-\frac12=2r_m-\hat c_m(\tfrac12)\,\frac{4r_m-1}{3}-\frac12=\frac{(4r_m-1)^2\,(4r_m+5)}{2\,(2r_m+1)\,(32r_m-5)}\ \ge0\qquad\text{on }\{r_m>5/32\},
\]
an event of probability tending to one. On this event let $S_{(k)}$ be the largest calibration score not exceeding $1/2$; almost surely $S_{(k)}<1/2$, and $\tilde F_m$ is constant on $[S_{(k)},1/2]$, so that
\[
\tilde F_m(S_{(k)})\ge\tfrac12,\qquad\tilde\tau(\tfrac12)\le S_{(k)}<\tfrac12,\qquad F_Q\bigl(\tilde\tau(\tfrac12)\bigr)<\tfrac12 .
\]
Hence
\[
\Prob\Bigl(P_e\bigl(\tilde\tau(\tfrac12)\bigr)>\tfrac12\Bigr)\longrightarrow1 :
\]
the efficient threshold run at the nominal level under-covers with probability tending to one, by an amount that is $o_p(m^{-1/2})$ in accordance with Theorem~\ref{thm:effthreshold}, and the Gaussian limit, degenerate here, does not justify an inflation proportional to $\hat\sigma_{\rm eff}(1/2)\toP0$. Degenerate levels require a separate, higher-order analysis.
\end{example}

\subsection{Proof of Proposition~\ref{prop:bounds23}}

\begin{proof}
(a) In $\mathcal M_2$ the marginal $P_X$ is fixed, so the closed tangent space is $\mathcal T_{S|X}$ and the efficient influence function is the projection of $g_t$ onto it, $g_t-\E[g_t\mid X]=v(\ind\{S\le t\}-m_t(X))$, with second moment $\E_P[v^2m_t(1-m_t)]$. The inequality $V_2\le V_1$ is
\[
V_1(t)-V_2(t)=\Var_P\bigl(vm_t(X)\bigr)-\frac{\E_P[(vm_t-F_Q(t))(v-1)]^2}{\E_P(v-1)^2}\ge0,
\]
by the Cauchy--Schwarz inequality.

(b) Formulate $\mathcal M_3$ as a pooled i.i.d.\ sample of size $n=m+N$ from the law of $(R,X,RS)$, where $R\sim\mathrm{Bernoulli}(\pi)$, $X\mid R=1\sim P_X$, $X\mid R=0\sim Q_X$, and $S\mid X,R=1\sim P_{S|X}$; the functional is $\theta_t=\E[m_t(X)\mid R=0]$, which is identified because $S\perp R\mid X$. This is the missing-at-random mean with propensity $e(x)=\Prob(R=1\mid X=x)$ and $(1-e(x))/e(x)=\frac{1-\pi}{\pi}v(x)$. The efficient influence function in the nonparametric model is \citep{hahn1998role}, in the form of \citet[Lemma~1]{yang2024doubly},
\[
\varphi_t=\frac{R}{\pi}\,v(X)\bigl\{\ind\{S\le t\}-m_t(X)\bigr\}+\frac{1-R}{1-\pi}\bigl\{m_t(X)-\theta_t\bigr\},
\]
with second moment $\E\varphi_t^2=\pi^{-1}\E_P[v^2m_t(1-m_t)]+(1-\pi)^{-1}\Var_Q(m_t(X))$. This is the bound in the normalization by $\sqrt n$; multiplying by $\pi=m/n$ gives the normalization by $\sqrt m$ and $V_3(t)$, since $\pi/(1-\pi)=m/N$. The limit as $N/m\to\infty$ is immediate.
\end{proof}

\end{document}